\documentclass[a4paper]{article}
\usepackage{amscd,amsmath,amssymb,amsfonts,times}
\usepackage{mathrsfs}
\usepackage[dvips]{color} 
\usepackage[all]{xy}

\numberwithin{equation}{section}
    \newtheorem{thm}{Theorem}[section]
    \newtheorem{lem}[thm]{Lemma}
    
    \newtheorem{prop}[thm]{Proposition}
    \newtheorem{cor}[thm]{Corollary}

    \newtheorem{conj}[thm]{Conjecture}
    \newtheorem{defn}[thm]{Definition}
    
    \newtheorem{rem}[thm]{Remark}

\DeclareMathAlphabet{\mathpzc}{OT1}{pzc}{m}{it}
\newcommand{\qed}
{\mbox{}\nolinebreak$\square$\medbreak\par}
\newenvironment{pf}{\par\smallskip\noindent\emph{Proof.}}{\hfill\qed\par\smallskip}
\newenvironment{pf*}[1]{\par\smallskip\noindent\emph{#1.}}{\hfill\qed\par\smallskip}
\begin{document}
\title{New $p$-adic hypergeometric functions and syntomic regulators}
\author{M. Asakura}
\date\empty

\maketitle

\def\can{\mathrm{can}}
\def\cano{\mathrm{canonical}}
\def\ch{{\mathrm{ch}}}
\def\Coker{\mathrm{Coker}}
\def\crys{\mathrm{crys}}
\def\dlog{d{\mathrm{log}}}
\def\dR{{\mathrm{d\hspace{-0.2pt}R}}}            
\def\et{{\mathrm{\acute{e}t}}}  
\def\Frac{{\mathrm{Frac}}}
\def\phami{\phantom{-}}
\def\id{{\mathrm{id}}}              
\def\Image{{\mathrm{Im}}}        
\def\Hom{{\mathrm{Hom}}}  
\def\Ext{{\mathrm{Ext}}}
\def\MHS{{\mathrm{MHS}}}  
  
\def\ker{{\mathrm{Ker}}}          
\def\rig{{\mathrm{rig}}}
\def\Pic{{\mathrm{Pic}}}
\def\CH{{\mathrm{CH}}}
\def\NS{{\mathrm{NS}}}
\def\NF{{\mathrm{NF}}}
\def\End{{\mathrm{End}}}
\def\pr{{\mathrm{pr}}}
\def\Proj{{\mathrm{Proj}}}
\def\ord{{\mathrm{ord}}}
\def\reg{{\mathrm{reg}}}          %
\def\res{{\mathrm{res}}}          %
\def\Res{\mathrm{Res}}
\def\Spec{\operatorname{Spec}}     
\def\syn{{\mathrm{syn}}}
\def\unit{{\mathrm{unit}}}
\def\zar{{\mathrm{zar}}}
\def\bA{{\mathbb A}}
\def\bC{{\mathbb C}}
\def\C{{\mathbb C}}
\def\G{{\mathbb G}}
\def\bE{{\mathbb E}}
\def\bF{{\mathbb F}}
\def\F{{\mathbb F}}
\def\bG{{\mathbb G}}
\def\bH{{\mathbb H}}
\def\bJ{{\mathbb J}}
\def\bL{{\mathbb L}}
\def\cL{{\mathscr L}}
\def\bN{{\mathbb N}}
\def\bP{{\mathbb P}}
\def\P{{\mathbb P}}
\def\bQ{{\mathbb Q}}
\def\Q{{\mathbb Q}}
\def\bR{{\mathbb R}}
\def\R{{\mathbb R}}
\def\bZ{{\mathbb Z}}
\def\Z{{\mathbb Z}}
\def\cH{{\mathscr H}}
\def\cD{{\mathscr D}}
\def\cE{{\mathscr E}}
\def\cF{{\mathscr F}}
\def\O{{\mathscr O}}
\def\cR{{\mathscr R}}
\def\cS{{\mathscr S}}
\def\cX{{\mathscr X}}
\def\cY{{\mathscr Y}}
%
\def\ve{\varepsilon}
\def\vG{\varGamma}
\def\vg{\varGamma}
%
%
%
%
\def\lra{\longrightarrow}
\def\lla{\longleftarrow}
\def\Lra{\Longrightarrow}
\def\hra{\hookrightarrow}
\def\lmt{\longmapsto}
\def\ot{\otimes}
\def\op{\oplus}
\def\l{\lambda}
\def\Isoc{{\mathrm{Isoc}}}
\def\Fil{{\mathrm{Fil}}}
\def\FIsoc{{F\text{-Isoc}}}
\def\FMIC{{F\text{-MIC}}}

\def\FilFMIC{{\mathrm{Fil}\text{-}F\text{-}\mathrm{MIC}}}

\def\Tr{\mathrm{Tr}}
\def\gauss{\text{\rm Gauss}}
\def\Log{{\mathscr{L}{og}}}
\def\wt#1{\widetilde{#1}}
\def\wh#1{\widehat{#1}}
\def\spt{\sptilde}
\def\ol#1{\overline{#1}}
\def\ul#1{\underline{#1}}
\def\us#1#2{\underset{#1}{#2}}
\def\os#1#2{\overset{#1}{#2}}
\def\wA{\wt A}

\begin{abstract}
We introduce new $p$-adic convergent functions, which we call the $p$-adic 
hypergeometric functions of logarithmic type.
The first main result is to prove the congruence relations that are similar to Dwork's.
The second main result
is that the special values of our new functions appear in the syntomic regulators for
hypergeometric curves, Fermat curves and some elliptic curves.
According to the $p$-adic Beilinson conjecture by Perrin-Riou, 
they are expected to be related with the special values of $p$-adic $L$-functions.
We provide one example for this.
\end{abstract}


\begin{flushleft}
\textbf{MSC} (2020) : 14F30, 19F27, 11S80 (primary), 19F15 (secondary)
\end{flushleft}

\section{Introduction}
Let $s\geq 1$ be an integer. 
For a $s$-tuple $\ul a=(a_1,\ldots,a_s)\in\Z^s_p$ of $p$-adic integers,
let
\[
F_{\ul a}(t)={}_sF_{s-1}\left({a_1,\ldots,a_s\atop 1,\ldots, 1}:t\right)
=\sum_{n=0}^\infty\frac{(a_1)_n}{n!}\cdots\frac{(a_s)_n}{n!}t^n
\]
be the hypergeometric power series where
$(\alpha)_n=\alpha(\alpha+1)\cdots(\alpha+n-1)$ denotes the Pochhammer symbol.
In his seminal paper \cite{Dwork-p-cycle}, B. Dwork discovered that
certain ratios of the hypergeometric power series are the uniform limit of rational functions.
We call his functions {\it Dwork's $p$-adic hypergeometric
functions}.
Let $\alpha'$ denote the Dwork prime, which is defined to be $(\alpha+l)/p$ where $l\in
\{0,1,\ldots,p-1\}$ is the unique integer such that $\alpha+l\equiv0$ mod $p$.
The $i$-th Dwork prime $a^{(i)}$ is defined by $a^{(i+1)}=(a^{(i)})'$ and $a^{(0)}=a$.
Write $\ul a'=(a'_1,\ldots,a'_s)$ and $\ul a^{(i)}=(a^{(i)}_1,\ldots,a^{(i)}_s)$. 
Then Dwork's $p$-adic hypergeometric function is defined
to be
\begin{equation}\label{Dwork}
\cF^{\mathrm{Dw}}_{\ul a}(t)=F_{\ul a}(t)/F_{\ul a'}(t^p).
\end{equation}
This is a convergent function in the sense of Krasner.
More precisely 
Dwork proved the congruence relations (\cite[p.41, Theorem 3]{Dwork-p-cycle})
\begin{equation}\label{Dwork-congruence}
\cF^{\mathrm{Dw}}_{\ul a}(t)
\equiv \frac{[F_{\ul a}(t)]_{<p^n}}{[F_{\ul a'}(t^p)]_{<p^n}}\mod p^n\Z_p[[t]]
\end{equation}
where for a power series $f(t)=\sum c_nt^n$, we write 
$[f(t)]_{<m}:=\sum_{n<m}c_nt^n$ the truncated polynomial.
This implies that 
$\cF^{\mathrm{Dw}}_{\ul a}(t)$ is a convergent function.
More precisely, for $f(t)\in\Z_p[t]$,
let $\ol{f(t)}\in \F_p[t]$ denote the reduction modulo $p$. 
Let $I\subset \Z_{\geq0}$ be a finite subset 
such that $\{\ol{[F_{\ul a^{(i)}}(t)]}_{<p}\}_{i\in \Z_{\geq0}}
=\{\ol{[F_{\ul a^{(i)}}(t)]}_{<p}\}_{i\in I}$
as sets. Put $h(t)=\prod_{i\in I}[F_{\ul a^{(i)}}(t)]_{<p}$.
Then \eqref{Dwork-congruence} implies
\[
\cF^{\mathrm{Dw}}_{\ul a}(t)
\in\Z_p\langle t,h(t)^{-1}\rangle:=\varprojlim_{n\geq1}(\Z_p/p^n\Z_p[t,h(t)^{-1}]),
\]
and hence that $\cF^{\mathrm{Dw}}_{\ul a}(t)$ is a convergent function on a domain
$\{z\in\C_p\mid |h(z)|_p=1\}$.

Dwork showed a geometric aspect of his $p$-adic hypergeometric functions.
Let $E$ be the elliptic curve over $\F_p$ defined by a Weierstrass equation
$y^2=x(1-x)(1- a x)$ with $a\in\F_p\setminus\{0,1\}$. 
Suppose that $E$ is ordinary, which means that $p\not{\hspace{-0.5mm}|}\,a_p$ 
where 
$T^2-a_pT+p$ is the characteristic polynomial of the Frobenius on $E$.
Let $\alpha_E$ be the root of $T^2-a_pT+p$ in $\Z_p$ such that $|\alpha_E|_p=1$,
which is often referred to as the unit root.
Then Dwork proved a formula
\[
\alpha_E=(-1)^{\frac{p-1}{2}}\cF^{\mathrm{Dw}}_{\frac{1}{2},\frac{1}{2}}(\wh a)
\]
where $\wh a\in \Z_p^\times$ is the Teichm\"uller lift of $a\in \F_p^\times$.
This is now called Dwork's unit root formula (cf. \cite[\S 7]{Put})

\medskip

In this paper, we introduce new $p$-adic hypergeometric functions, which we call
the {\it $p$-adic hypergeometric functions of logarithmic type}.
We shall first introduce the $p$-adic polygamma functions $\psi^{(r)}_p(z)$ 
in \S \ref{poly-sect}, which are slight modifications of Diamond's $p$-adic
polygamma functions \cite{Diamond}.
Let $W=W(\ol\F_p)$ be the Witt ring of $\ol\F_p$, and $K=\Frac W$ the fractional field.
Let $\sigma$ be a $p$-th Frobenius on $W[[t]]$ 
given by $\sigma(t)=ct^p$
with $c\in 1+pW$. Let $\log:\C^\times_p\to\C_p$ be the Iwasawa logarithmic function.
Let
\[
G_{\ul a}(t):=\psi_p(a_1)+\cdots+\psi_p(a_s)+s\gamma_p
-p^{-1}\log(c)+
\int_0^t(F_{\underline{a}}(t)-F_{\ul a'}(t^\sigma))\frac{dt}{t}
\]
be a power series
where 
$\int_0^t(-)\frac{dt}{t}$ means the 
operator which sends $\sum_{n\geq 1}a_nt^n$ to $\sum_{n\geq 1}\frac{a_n}{n}t^n$.
It is not hard to show $G_{\ul a}(t)\in W[[t]]$ (Lemma \ref{G-int}).
Then our new function is defined to be a ratio
\[
\cF^{(\sigma)}_{\underline{a}}(t):=
G_{\ul a}(t)/F_{\ul a}(t).
\]
Notice that $\cF^{(\sigma)}_{\ul a}(t)$ is also $p$-adically continuous with respect to
$\ul a$.
In case $a_1=\cdots=a_s=c=1$, one has $\cF^{(\sigma)}_{\ul a}(t)
=(1-t)\ln^{(p)}_1(t)$ the $p$-adic logarithm. In this sense, we can regard 
$\cF^{(\sigma)}_{\ul a}(t)$ as a deformation of the $p$-adic logarithm.

The first main result of this paper 
is the congruence relations for $\cF^{(\sigma)}_{\ul a}(t)$ that are 
similar to Dwork's. 
\begin{thm}[Theorem \ref{cong-thm}]\label{intro-thm0}
Suppose that $a_i\not\in\Z_{\leq 0}$ for all $i$. Then
\[
\cF^{(\sigma)}_{\underline{a}}(t)\equiv \frac{[G_{\ul a}(t)]_{<p^n}}{[F_{\ul a}(t)]_{<p^n}}
\mod p^nW[[t]]
\]
if $c\in 1+2pW$. If $p=2$ and $c\in 1+2W$, then the congruence holds modulo $p^{n-1}$.
\end{thm}
Thanks to this, 
$\cF^{(\sigma)}_{\ul a}(t)$ is a convergent function, 
\[
\cF^{(\sigma)}_{\ul a}(t)
\in\Z_p\langle t,h(t)^{-1}\rangle
\]
and then the special value at $t=\alpha$ is defined for $\alpha\in\C_p$ such that 
$|h(\alpha)|_p=1$.

\medskip

The second main result 
is to give a geometric aspect of our $\cF^{(\sigma)}_{\ul a}(t)$,
which concerns with the {\it syntomic regulator map}.
Let $X$ be a smooth variety over $W$.
Let $H^\bullet_{\syn}(X,\Q_p(j))$ be the rigid syntomic cohomology groups by
Besser \cite{Be1} (see also \cite[1B]{NN}), which agree with the syntomic cohomology
groups by Fontaine-Messing \cite{FM} (see also \cite{KaV}) when $X$ is projective. 
Let $K_i(X)$ be Quillen's algebraic $K$-groups.
Then there is the syntomic regulator map
\[
\reg_\syn^{i,j}:K_i(X)\lra H^{2j-i}_\syn(X,\Q_p(j))
\]
for each $i,j\geq0$ (\cite[Theorem 7.5]{Be1}, \cite[Theorem A]{NN}).
We shall concern ourselves with only $\reg_\syn^{2,2}$, which we 
abbreviate to $\reg_\syn$ in this paper.
Note that there is the natural isomorphism
$H^2_\syn(X,\Q_p(2))\cong H^1_\dR(X_K/K)$ where $K=\Frac W$ is the fractional
field and $X_K:=X\times_WK$.
Our second main result is to relate $\cF^{(\sigma)}_{a,b}(t)$ with
the syntomic regulator of a certain element of $K_2$ of a 
{\it hypergeometric curve}, which is defined 
in the following way.
Let $N,M\geq 2$ be integers, and
$p$ a prime such that $p\nmid NM$.
Let $\alpha\in W$ satisfy that $\alpha\not\equiv 0,1$ mod $p$.
Then we define a hypergeometric curve $X_\alpha$ to
be a projective smooth scheme over $W$ defined by
a bihomogeneous equation 
\[
(X_0^N-X_1^N)(Y_0^M-Y_1^M)=\alpha X^N_0Y_0^M
\]
in $\P^1_W(X_0,X_1)\times\P^1_W(Y_0,Y_1)$ (\S \ref{fermat-sect}).
We put $x=X_1/X_0$ and $y=Y_1/Y_0$.

\begin{thm}[Corollary \ref{main-thm4}]\label{intro-thm1}
Suppose $p>\max(N,M)$.
Let
\[
\xi=\left\{\frac{x-1}{x-\nu_1},\frac{y-1}{y-\nu_2}\right\}\in K_2(X_\alpha)\ot\Q
\]
for $(\nu_1,\nu_2)\in\mu_N(K)\times\mu_M(K)$ where $\mu_m(K)$ denotes
the group of $m$-th roots of unity in $K$ (cf. \eqref{m-fermat-eq1}).
Let
\[
Q: H^1_\dR(X_{\alpha,K}/K)\ot H^1_\dR(X_{\alpha,K}/K)\lra 
H^2_\dR(X_{\alpha,K}/K)\cong K
\]
be the cup-product pairing.
Suppose that
$h(\alpha)\not\equiv0$ mod $p$ where $h(t)$ is as above.
For a pair of integers $(i,j)$ such that $0<i<N$ and $0<j<M$, we put
$\omega_{i,j}:=Nx^{i-1}y^{j-M}/(1-x^N)dx$ a regular $1$-form (cf. \eqref{fermat-form}),
and
$e_{i,j}^{\text{\rm unit}}$ the unit root vectors which are explicitly given in Theorem \ref{uroot-thm}.
Then we have
\[
Q(\reg_\syn(\xi), e_{N-i,M-j}^{\text{\rm unit}})
=N^{-1}M^{-1}
(1-\nu^{-i}_1)(1-\nu^{-j}_2)
\cF_{a_i,b_j}^{(\sigma_\alpha)}(\alpha)
Q(\omega_{i,j}, e_{N-i,M-j}^{\text{\rm unit}}).
\]
\end{thm}
For some elliptic curves,
we also have similar results to the above (see \S \ref{gauss-sect} -- \S \ref{elliptic-sect}).
For the proof of Theorem \ref{intro-thm1}, we employ the main result in \cite{AM}
as a fundamental tool to compute the syntomic regulators.
We share a part of the technique
in the proof of \cite[Theorem 4.8]{AM} (see also Remark \ref{rem.miya} below).
A new ingredient is our function $\cF_{a,b}^{(\sigma)}(t)$
(indeed Theorem \ref{intro-thm0}
plays a key role in the proof of Theorem \ref{intro-thm1}).
This paper focuses on 
the cup-product of the regulator with a generator of the unit root
subspace, that appears in the $p$-adic Beilinson conjecture by Perrin-Riou. 

\medskip

A more striking application of our new function $\cF^{(\sigma)}_{\ul a}(t)$
is that one can describe
the syntomic regulators of the Ross symbols in $K_2$ of the Fermat curves.
\begin{thm}[Theorem \ref{fermat-main2}]
Suppose that $N|(p-1)$ and $M|(p-1)$.
Let $F$ be the Fermat curve over $K$ defined by an affine equation $z^N+w^M=1$.
Let $\{1-z,1-w\}\in K_2(F)$ be the Ross symbol \cite{ross2}.
Let 
\[
\reg_\syn(\{1-z,1-w\})=
\sum_{(i,j)\in I} A_{i,j}M^{-1}z^{i-1}w^{j-M}dz\in H^1_\dR(F/K)
\]
where the notation be as in the beginning of \S \ref{fermatcurve-sect}.
Then \[A_{i,j}=\cF_{a_i,b_j}^{(\sigma)}(1)\] for $(i,j)$ such that $i/N+j/M<1$.
\end{thm}
As long as the author knows, this is the first explicit description of the Ross symbol
in $p$-adic cohomology.

\medskip

In the last section, we discuss the $p$-adic Beilinson conjecture
by Perrin-Riou \cite[4.2.2]{Perrin-Riou} (see also \cite[Conj.2.7]{Colmez})
for $K_2$ of elliptic curves.
Let $E$ be an elliptic curve over $\Q$.
Let $p$ be a prime at which $E$
has a good ordinary reduction.
Let $\alpha_{E,p}$ be the unit root of the reduction $\ol E$ at $p$, 
and $e_\unit\in H^1_\dR(E/\Q_p)$ the eigenvector with eigenvalue $\alpha_{E,p}$. 
Let $\omega_E\in \vg(E,\Omega^1_{E/\Q})$ be a regular $1$-form.
Let $L_p(E,\chi,s)$ be the $p$-adic $L$-function
defined by Mazur and Swinnerton-Dyer \cite{MS}.
Then as a consequence of the $p$-adic Beilinson conjecture 
for elliptic curves,
one can expect
that there is an element $\xi\in K_2(E)$ which is integral in the sense of \cite{Scholl}
such that
\[
(1-p\alpha^{-1}_{E,p})\frac{Q(\reg_\syn(\xi),e_\unit)}{
Q(\omega_E,e_\unit)}\sim_{\Q^\times}L_p(E,\omega^{-1},0)
\]
where $x\sim_{\Q^\times}y$ means $xy\ne0$ and $x/y\in \Q^\times$.
We also refer to \cite[Conjecture 3.3]{ABC} for more precise statement.
According to our main results, we can replace the left hand side with the special values
of the $p$-adic hypergeometric functions of logarithmic type.
For example, let $E_a$ be the elliptic curve defined by $y^2=x(1-x)(1-(1-a)x)$ with $a\in \Q\setminus\{0,1\}$ and $p>3$ a prime where $E_a$ has a good ordinary reduction.
Then one predicts
\[
(1-p\alpha^{-1}_{E_a,p})
\cF_{\frac{1}{2},\frac{1}{2}}^{(\sigma_a)}(a)\sim_{\Q^\times}
L_p(E_a,\omega^{-1},0)
\]
if $a=-1,\pm2,\pm4,\pm8,\pm16,
\pm\frac{1}{2},  
\pm\frac{1}{8},\pm\frac{1}{4}, \pm\frac{1}{16}$ (Conjecture \ref{ell1-conj}).
See \S \ref{RZ-sect} for other cases.
The author has no idea how to attack the question in general, while
we have one example (the proof relies on Brunault's paper \cite{regexp}
and his appendix in \cite{ABC}).
\begin{thm}[Theorem \ref{brunault}]\label{F. Brunault}
$(1-p\alpha^{-1}_{E_4,p})\cF_{\frac{1}{2},\frac{1}{2}}^{(\sigma_4)}(4)=-
L_p(E_4,\omega^{-1},0).$
\end{thm}
We note that this is a $p$-adic counterpart of 
a formula of Rogers and Zudilin (\cite[Theorem 2, p.399 and (6), p.386]{RZ})
\[
2L'(E_4,0)=
\mathrm{Re}\left[\log 4-
{}_4F_3\left({\frac{3}{2},\frac{3}{2},1,1\atop 2,2,2};4\right)\right]
\left(={}_3F_2\left({\frac{1}{2},\frac{1}{2},\frac{1}{2}\atop \frac{3}{2},1};
\frac{1}{4}\right)\right).
\]
The conjectures in \S \ref{RZ-sect} give the first formulation of 
the $p$-adic counterparts
of Rogers-Zudilin type formulas.
We hope that our results will provide
a new direction of the study of the $p$-adic Beilinson conjecture.

\medskip

This paper is organized as follows.
\S \ref{poly-sect} is the preliminary section on 
Diamond's $p$-adic polygamma functions.
More precisely we shall give a slight modification of Diamond's polygamma functions
(though it might be known to the experts).
We give a self-contained exposition, because
the author does not find a suitable reference, especially 
concerning with our modified functions.
In \S \ref{pHGlog-sect}, we introduce the $p$-adic hypergeometric functions of
logarithmic type, and prove the congruence relations.
In \S \ref{reg-sect}, we show that our new $p$-adic hypergeometric functions appear in
the syntomic regulators of the hypergeometric curves.
Finally we discuss the $p$-adic Beilinson conjecture for $K_2$ of elliptic curves
in \S \ref{weakB-sect}.

\medskip

\noindent{\bf Acknowledgement.}
The author would like to express sincere gratitude
to Professor Masataka Chida for the stimulating discussion on 
the $p$-adic Beilinson conjecture. The discussion with him is the origin of this work.
He very much appreciate Professor Fran\c{c}ois Brunault 
for giving him lots of comments on his paper \cite{regexp}, and for the help of the proof
of Theorem \ref{F. Brunault}.

\medskip

\noindent{\bf Notation.}
For a field $K$, let $\mu_n(K)\subset K^\times$ denote the group of $n$-th roots of unity. We write $\mu_\infty(K)=\cup_{n\geq 1}\mu_n(K)$.
If there is no fear of confusion, we drop ``$K$'' and simply write $\mu_n$.
For a power series $f(t)=\sum_{i=0}^\infty a_it^i\in R[[t]]$ with coefficients in a 
commutative ring $R$, 
we write the truncated polynomial $\sum_{i=0}^{n-1}a_it^i$ by $[f(t)]_{<n}$.

\section{$p$-adic polygamma functions}\label{poly-sect}
The complex analytic polygamma functions are the $r$-th derivative
\[
\psi^{(r)}(z):=\frac{d^r}{dz^r}\left(\frac{\Gamma'(z)}{\Gamma(z)}\right), \quad r\in\Z_{\geq0}.
\]
In his paper \cite{Diamond}, Jack Diamond gave a $p$-adic counterpart of the polygamma functions $\psi^{(r)}_{D,p}(z)$ which are given in the following way.
\begin{equation}\label{diamond0}
\psi^{(0)}_{D,p}(z)
=\lim_{s\to\infty}\frac{1}{p^s}\sum_{n=0}^{p^s-1}\log(z+n),
\end{equation}
\begin{equation}\label{diamond1}
\psi^{(r)}_{D,p}(z)=(-1)^{r+1}r!\lim_{s\to\infty}\frac{1}{p^s}\sum_{n=0}^{p^s-1}\frac{1}{(z+n)^r},
\quad r\geq 1,
\end{equation}
where $\log:\C_p^\times\to\C_p$ is the Iwasawa logarithmic function which is characterized as a continuous 
homomorphism satisfying
$\log(p)=0$ and
\[
\log(z)=-\sum_{n=1}^\infty\frac{(1-z)^n}{n},\quad |z-1|_p<1.
\]
It should be noticed that the series \eqref{diamond0} and \eqref{diamond1}
converge only when $z\not\in \Z_p$, and hence $\psi^{(r)}_{D,p}(z)$
turn out to be locally analytic functions on $\C_p\setminus\Z_p$.
This causes inconvenience in our discussion.
In this section we give a continuous function $\psi_p^{(r)}(z)$ on $\Z_p$
which is a slight {\it modification} of $\psi^{(r)}_{D,p}(z)$.
See \S \ref{polygamma-sect} for the definition and also \S \ref{measure-sect} 
for an alternative definition in terms of $p$-adic measure.

\subsection{$p$-adic polylogarithmic functions}
Let $x$ be an indeterminate.
For an integer $r\in\Z$, the $r$-th $p$-adic polylogarithmic function $\ln_r^{(p)}(x)$
is defined as a formal power series
\begin{equation}\label{polylog-defn}
\ln_r^{(p)}(x):=\sum_{k\geq1,p\nmid k}\frac{x^k}{k^r}=\lim_{s\to\infty}
\left(\frac{1}{1-x^{p^s}}\sum_{1\leq k<p^s,\, p\nmid k}\frac{x^k}{k^r}\right)
\in \Z_p[[x]]
\end{equation}
which belongs to the ring 
\[
\Z_p\langle x,(1-x)^{-1}\rangle:=\varprojlim_{s}(\Z/p^s\Z[x,(1-x)^{-1}])
\]
of convergent power series.
If $r\leq 0$, this is a rational function, more precisely
\[
\ln_0^{(p)}(x)=\frac{1}{1-x}-\frac{1}{1-x^p},\quad 
\ln_{-r}^{(p)}(x)=\left(x\frac{d}{dx}\right)^r\ln_0^{(p)}(x).
\]
If $r>0$, this is known to be an {\it overconvergent function}, more precisely
it has a (unique) analytic continuation to the domain 
$|x-1|>|1-\zeta_p|$ where $\zeta_p\in \ol\Q_p$ 
is a primitive $p$-th root of unity.

Let $W(\ol\F_p)$ be the Witt ring of $\ol\F_p$ and $F$ the $p$-th Frobenius endomorphism.
Define the {\it $p$-adic logarithmic function}
\begin{equation}\label{log(p)}
\log^{(p)}(z):=\frac{1}{p}\log\left(\frac{z^p}{F(z)}\right)=-\sum_{n=1}^\infty
\frac{p^{-1}}{n}\left(1-\frac{z^p}{F(z)}\right)^n
\end{equation}
for $z\in W(\ol\F_p)^\times$, where $\log(z)$ is the Iwasawa logarithmic function.
\begin{lem}\label{lemma.log}
The function $\ln_1^{(p)}(x)$ agrees with 
\[
-\frac1p\log\left(\frac{(1-x)^p}{1-x^p}\right)=\sum_{n=1}^\infty
\frac{p^{n-1}w(x)^n}{n}\in \Z_p\langle x,(1-x)^{-1}\rangle
\]
where $w(x):=1-(1-x)^p/(1-x^p)$.
In particular, evaluating at $x=z$ for $z\in W(\ol\F_p)^\times$ such that $F(z)=z^p$
and $z\not\equiv 1$ mod $p$, one has $\ln^{(p)}_1(z)=-\log^{(p)}(1-z)$.
\end{lem}
\begin{pf}
We have the power series expression 
\[
-\frac1p\log\left(\frac{(1-x)^p}{1-x^p}\right)=\sum_{k\geq 1,\,p\nmid k}\frac{x^k}{k}
\]
in $\Z_p[[x]]$, and this agrees with the expression \eqref{polylog-defn} of $\ln_1^{(p)}(x)$.
Then the assertion is immediate 
as $\Z_p\langle x,(1-x)^{-1}\rangle\to\Z_p[[x]]$ is injective.
\end{pf}
\begin{prop}[cf. \cite{C-dlog} IV Prop.6.1, 6.2]
Let $r\in\Z$ be an integer. Then
\begin{equation}\label{polylog-diffeq}
\ln^{(p)}_r(x)=x\frac{d}{dx}\ln^{(p)}_{r+1}(x),
\end{equation}
\begin{equation}\label{polylog-refl}
\ln^{(p)}_r(x)=(-1)^{r+1}\ln^{(p)}_r(x^{-1}),
\end{equation}
\begin{equation}\label{polylog-dstr}
\sum_{\zeta\in\mu_N}\ln_r^{(p)}(\zeta x)=\frac{1}{N^{r-1}}\ln^{(p)}_r(x^N)
\quad \mbox{\rm(distribution formula)}.
\end{equation}
\end{prop}
\begin{pf}
\eqref{polylog-diffeq} and \eqref{polylog-dstr} are immediate from the power series
expansion $\ln_r^{(p)}(x)=\sum_{k\geq 1,p\nmid k} x^k/k^r$.
On the other hand \eqref{polylog-refl} follows from the fact
\[
\frac{1}{1-x^{-p^s}}\sum_{1\leq k<p^s,\, p\nmid k}\frac{x^{-k}}{k^r}
=
\frac{-1}{1-x^{p^s}}\sum_{1\leq k<p^s,\, p\nmid k}\frac{x^{p^s-k}}{k^r}
\equiv
\frac{(-1)^{r+1}}{1-x^{p^s}}\sum_{1\leq k<p^s,\, p\nmid k}
\frac{x^{p^s-k}}{(p^s-k)^r}
\]
modulo $p^s\Z[x,(1-x)^{-1}]$.
\end{pf}
\begin{lem}\label{ln-formula1}
Let $m, N\geq 2$ be integers prime to $p$.
Let $\ve\in\mu_m\setminus\{1\}$.
Then for any $n\in \{0,1,\ldots,N-1\}$, we have
\[
N^r\sum_{\nu^N=\ve}\nu^{-n}\ln^{(p)}_{r+1}(\nu)=\lim_{s\to\infty}
\frac{1}{1-\ve^{p^s}}\sum_{\us{k+n/N\not\equiv0\text{ \rm mod }p}{0\leq k<p^s}}
\frac{\ve^k}{(k+n/N)^{r+1}}.
\]  
\end{lem}
\begin{pf}
Note $\sum_{\nu^N=\ve}\nu^i=N\ve^{i/N}$ if $N|i$ and $=0$ otherwise.
We have 
\begin{align*}
N^r\sum_{\nu^N=\ve}\nu^{-n}\ln^{(p)}_{r+1}(\nu x)
&=N^r\sum_{k\geq1,p\nmid k}\sum_{\nu^N=\ve}\frac{\nu^{k-n}x^k}{k^{r+1}}\\
&=N^{r+1}\sum_{N|(k-n),p\nmid k}\frac{\ve^{(k-n)/N}x^k}{k^{r+1}}\\
&=\sum_{n+\ell N\not\equiv 0\text{ mod }p,\,\ell\geq 0}\frac{(\ve x^N)^\ell x^n}{(\ell+n/N)^{r+1}}
\quad(\ell=(k-n)/N)\\
&\equiv \frac{1}{1-(\ve x^N)^{p^s}}\sum_{\us{n+\ell N\not\equiv 0
\text{ mod }p}{0\leq \ell<p^s}}\frac{(\ve x^N)^\ell x^n}{(\ell+n/N)^{r+1}}\\
\end{align*}
modulo $p^s\Z[x,(1-\ve x^N)^{-1}]$.
Since $\ve\not\equiv 1$ mod $p$, the evaluation at $x=1$ makes sense, and then we have the desired equation.
\end{pf}
The following theorem  
is well-known to experts as Coleman's formula.
\begin{thm}[Coleman]\label{p-zeta-def}
Let $r\ne1$ be an integer.
Then for any integer $N\geq 2$ prime to $p$,
\begin{equation}\label{p-zeta-def-eq1}
\sum_{\ve\in \mu_N\setminus\{1\}}\ln_r^{(p)}(\ve)
=-(1-N^{1-r})L_p(r,\omega^{1-r})
\end{equation}
where $L_p(s,\chi)$ is the $p$-adic $L$-function 
and $\omega$ is the Teichm\"uller character.
\end{thm}
\begin{pf}
We give a self-contained and straightforward proof for convenience of the reader,
because the author does not find a suitable literature (note that \eqref{p-zeta-def-eq1} 
is not covered by \cite[I, (3)]{C-dlog}).

\medskip


We first show \eqref{p-zeta-def-eq1} in case $r=-m$ with $m\in \Z_{\geq1}$.
Note that $\ln_{-m}^{(p)}(x)$ is a rational function.
More precisely, let
\[
\ln_0(x):=\frac{x}{1-x},\quad 
\ln_{-m}(x):=\left(x\frac{d}{dx}\right)^m\ln_0(x),
\]
then $\ln_{-m}^{(p)}(x)=\ln_{-m}(x)-p^m\ln_{-m}(x^p)$.
Therefore
\[
\sum_{\ve\in \mu_N\setminus\{1\}}\ln_{-m}^{(p)}(\ve)
=(1-p^{m})
\sum_{\ve\in \mu_N\setminus\{1\}}\ln_{-m}(\ve).
\]
Since $L_p(-m,\omega^{1+m})=-(1-p^{m})B_{m+1}/(m+1)$ 
where $B_n$ are the Bernoulli numbers,
the equation
\eqref{p-zeta-def-eq1} for $r=-m$ is equivalent to
\begin{equation}\label{p-zeta-def-eq4}
(1-N^{m+1})\frac{B_{m+1}}{m+1}=
\sum_{\ve\in \mu_N\setminus\{1\}}\ln_{-m}(\ve).
\end{equation}
Put $\ell_r(x):=\ln_r(x)-N^{1-r}\ln_r(x^N)$.
By the distribution property
\[
\sum_{\ve\in \mu_N}\ln_r(\ve x)=
N^{1-r}\ln_r(x^N)
\]
which can be easily shown by a computation of power series expansions,
the right hand side of \eqref{p-zeta-def-eq4} equals to
the evaluation $-\ell_{-m}(x)|_{x=1}$ at $x=1$, and hence
\begin{equation}\label{p-zeta-def-eq5}
\sum_{\ve\in \mu_N\setminus\{1\}}\ln_{-m}(\ve)
=-\left(x\frac{d}{dx}\right)^m\ell_0(x)\bigg|_{x=1}.
\end{equation}
On the other hand, letting $x=e^z$, one has
\begin{align*}
\ell_0(e^z)=\frac{e^z}{1-e^z}-\frac{Ne^{Nz}}{1-e^{Nz}}
=-\sum_{n=1}^\infty \left(B_n\frac{z^{n-1}}{n!}-B_n\frac{N^nz^{n-1}}{n!}
\right)
=-\sum_{n=1}^\infty (1-N^n)B_n\frac{z^{n-1}}{n!}
\end{align*}
and hence
\begin{equation}\label{p-zeta-def-eq6}
\left(x\frac{d}{dx}\right)^m\ell_0(x)\bigg|_{x=1}=
\frac{d^m}{dz^m}\ell_0(e^z)\bigg|_{z=0}=
-(1-N^{m+1})\frac{B_{m+1}}{m+1}.
\end{equation}
Now \eqref{p-zeta-def-eq4} follows from \eqref{p-zeta-def-eq5} and
\eqref{p-zeta-def-eq6}.

\medskip

We have shown \eqref{p-zeta-def-eq1} for negative $r$.
Let $r\ne1$ be an arbitrary integer.
Since $\ln_r^{(p)}(x)=\sum_{p\nmid  k}x^k/k^r$,
one has that
for any integers $r,r'$ such that $r\equiv r'$ mod $(p-1)p^{s-1}$,
\[\ln_r^{(p)}(x)\equiv \ln_{r'}^{(p)}(x)
\mod p^s\Z_p[[x]]
\]
and hence modulo $p^s\Z_p\langle x,(1-x)^{-1}\rangle$.
This implies
\[\ln_r^{(p)}(\ve)\equiv \ln_{r'}^{(p)}(\ve)
\mod p^sW(\ol\F_p).
\]
Take $r'=r-p^{s+a-1}(p-1)$ with $a\gg 0$.
It follows that
\[
(1-N^{1-r'})L_p(r',\omega^{1-r'})
=(1-N^{1-r'})L_p(r',\omega^{1-r})\to
(1-N^{1-r})L_p(r,\omega^{1-r})
\]
as $a\to\infty$ by the continuity of the $p$-adic $L$-functions.
Since $r'<0$, one can apply \eqref{p-zeta-def-eq1} and then 
\[
-(1-N^{1-r})L_p(r,\omega^{1-r})
\equiv \sum_{\ve\in \mu_N\setminus\{1\}}\ln_r(\ve)\mod p^sW(\ol\F_p)
\]
for any $s>0$. This completes the proof.

\end{pf}

\subsection{$p$-adic polygamma functions}\label{polygamma-sect}
\begin{lem}\label{lemma.equiv}
\[
\sum_{1\leq k<p^s,p\nmid k}k^m\equiv 
\begin{cases}
-p^{s-1}&p\geq 3\mbox{ and }(p-1)|m\\
2^{s-1}&p=2\mbox{ and }2 |m\\
1&p=2\mbox{ and }s=1\\
0&\mbox{otherwise}
\end{cases}\mod p^s.
\]
\end{lem}
\begin{pf}
Let $p>2$. Let $\zeta\in \Z_p$ be a primitive $(p-1)$-th root of unity.
Then the set $\{\zeta^i(1+p)^i\mid 0\leq i<p^{s-1}(p-1)\}$ is a representative
of $(\Z/p^s\Z)^\times$.
Hence
\[
\sum_{1\leq k<p^s,p\nmid k}k^m\equiv 
\sum_{0\leq i<p^{s-1}(p-1)}(\zeta(1+p))^{mi}
=\frac{1-(\zeta(1+p))^{mp^{s-1}(p-1)}}{1-(\zeta(1+p))^{m}}
=\frac{1-(1+p)^{mp^{s-1}(p-1)}}{1-(\zeta(1+p))^{m}}
\mod p^s.
\] 
Note 
that $(1+p)^{m_0p^j}\equiv 1+m_0p^{j+1}$ mod $p^{j+2}$ for $p\nmid m_0$ and $j\geq0$.
Therefore, when $(p-1)\nmid m$, the last term vanishes.
If $(p-1)|m$,
then the last term is equivalent to $mp^s(p-1)/mp=p^{s-1}(p-1)\equiv-p^{s-1}$ modulo $p^s$.
This completes the proof in case $p>2$.
Let $p=2$.
When $s=1$, the proof is obvious. Suppose $s\geq 2$.
Then the set $\{\pm5^i\mid 0\leq i<2^{s-2}\}$ is a representative
of $(\Z/2^s\Z)^\times$.
Therefore
\[
\sum_{1\leq k<2^s,\,2\nmid k}k^m\equiv \sum_{0\leq i<2^{s-2}}5^{mi}+(-1)^m5^{mi}\mod 2^s.
\]
This vanishes when $m$ is odd. If $m$ is even, then the right hand side is
\[
2\sum_{0\leq i<2^{s-2}}5^{mi}
=2\frac{1-5^{2^{s-2}m}}{1-5^m}\equiv2\frac{2^sm}{4m}=2^{s-1}\mod 2^s
\]
as $(1+4)^{m_02^j}\equiv 1+2^{j+2}m_0$ mod $2^{j+3}$ for odd $m_0$ and $j\geq0$.
This completes the proof in case $p=2$.
\end{pf}
Let $r\in\Z$ be an integer. For $z\in \Z_p$, define
\begin{equation}\label{wt-polygamma-def}
\wt{\psi}_p^{(r)}(z):=\lim_{n\in\Z_{>0},n\to z}\sum_{1\leq k<n,p\nmid  k}
\frac{1}{k^{r+1}}.
\end{equation}
The limit exists by Lemma \ref{lemma.equiv}, and moreover it satisfies
\begin{equation}\label{equiv-psi}
z\equiv z'\hspace{-0.3cm}\mod p^s\Longrightarrow
\wt\psi^{(r)}_p(z)-\wt\psi^{(r)}_p(z') \equiv
\begin{cases}
0\hspace{-0.3cm}\mod p^s&p\geq 3\mbox{ and }(p-1)\not|(r+1)\\
0\hspace{-0.3cm}\mod p^s&p=2,\, s\geq 2\mbox{ and }2\not|(r+1)\\
0\hspace{-0.3cm}\mod p^{s-1}&\mbox{othewise.}
\end{cases}
\end{equation}
Thus, $\wt\psi^{(r)}_p(z)$ is a $p$-adic continuous function on $\Z_p$.
Define the {\it $p$-adic Euler constant} \footnote{This is different from
Diamond's $p$-adic Euler constant. His constant
is $p/(p-1)\gamma_p$, \cite[\S 7]{Diamond}.} by
\[
\gamma_p:=-\lim_{s\to\infty}\frac{1}{p^s}\sum_{0\leq j<p^s,p\nmid j}\log(j),\quad
(\log=\mbox{Iwasawa log}).
\]
where the convergence follows by
\begin{align}
\sum_{0\leq j<p^{s+1},p\nmid j}\log(j)-p
\sum_{0\leq j<p^s,p\nmid j}\log(j)
&=\sum_{k=0}^{p-1}\sum_{0\leq j<p^s,p\nmid j}
\log\left(1+\frac{kp^s}{j}\right)\notag\\
&\equiv\sum_{k=0}^{p-1}\sum_{0\leq j<p^s,p\nmid j}
\frac{kp^s}{j}\mod p^{2s-1}\notag\\
&\equiv 0\mod p^{2s-1}\label{equiv-log}\quad\text{(Lemma \ref{lemma.equiv})}.
\end{align}
We define the $r$-th {\it $p$-adic polygamma function} to be
\begin{equation}\label{polygamma-def}
\psi_p^{(r)}(z):=\begin{cases}
-\gamma_p+\wt{\psi}^{(0)}_p(z)&r=0\\
-L_p(1+r,\omega^{-r})+\wt{\psi}^{(r)}_p(z)&r\ne0.
\end{cases}
\end{equation}
If $r=0$, we also write $\psi_p(z)=\psi_p^{(0)}(z)$ and call it the {\it $p$-adic digamma function}.
\subsection{Formulas on $p$-adic polygamma functions}\label{formula-sect}
\begin{thm}\label{polygamma-thm1}
\begin{enumerate}
\item[\rm(1)] $\wt\psi^{(r)}_p(0)=\wt\psi^{(r)}_p(1)=0$ or equivalently
$\psi^{(r)}_p(0)=\psi^{(r)}_p(1)=-\gamma_p$ or $=-L_p(r+1,\omega^{-r})$.
\item[\rm(2)]
$\wt\psi^{(r)}_p(z)=(-1)^r\wt\psi^{(r)}_p(1-z)$ or equivalently
$\psi^{(r)}_p(z)=(-1)^r\psi^{(r)}_p(1-z)$ (note $L_p(1+r,\omega^{-r})=0$ for odd $r$).
\item[\rm(3)]
\[
\wt\psi^{(r)}_p(z+1)-\wt\psi^{(r)}_p(z)=\psi^{(r)}_p(z+1)-\psi^{(r)}_p(z)=\begin{cases}
z^{-r-1}&z\in \Z_p^\times\\
0&z\in p\Z_p.
\end{cases}
\]
\end{enumerate}
\end{thm}
Compare the above with the formulas on the complex analytic
polygamma functions, \cite[5.15.2, 5.15.5, 5.15.6]{NIST}.
\begin{pf}
(1) follows from Lemma \ref{lemma.equiv},
and (3) are immediate from definition.
We show (2). 
Since $\Z_{>0}$ is a dense subset in $\Z_p$, it is enough to show in case $z=n>0$ an integer.
Let $s>0$ be arbitrary such that $p^s>n$.
Then
\begin{align*}
\wt\psi^{(r)}_p(n)&\equiv
\sum_{1\leq k<n,p\nmid k}\frac{1}{k^{r+1}}\equiv 
(-1)^{r+1}\sum_{-n< k\leq -1,p\nmid k}\frac{1}{k^{r+1}}
\equiv 
(-1)^{r+1}\sum_{p^s-n+1\leq k< p^s,p\nmid k}\frac{1}{k^{r+1}}\\
&\equiv 
(-1)^{r+1}\sum_{0\leq k< p^s,p\nmid k}\frac{1}{k^{r+1}}-
(-1)^{r+1}\sum_{0\leq k< p^s-n+1,p\nmid k}\frac{1}{k^{r+1}}\\
&\equiv(-1)^{r}\sum_{0\leq k< p^s-n+1,p\nmid k}\frac{1}{k^{r+1}}\\
&\equiv (-1)^r\wt\psi^{(r)}_p(1-n)
\end{align*}
modulo $p^s$ or $p^{s-1}$. Since $s$ is an arbitrary large integer, this means 
$\wt\psi^{(r)}_p(n)=(-1)^r\wt\psi^{(r)}_p(1-n)$ as required.
\end{pf}
\begin{thm}\label{polygamma-thm2}
Let $0\leq n<N$ be integers and suppose $p\not|N$. Then
\begin{equation}\label{polygamma-thm2-eq}
\wt{\psi}_p^{(r)}\left(\frac{n}{N}\right)=N^r
\sum_{\ve\in \mu_N\setminus\{1\}}(1-\ve^{-n})\ln_{r+1}^{(p)}(\ve).
\end{equation}
\end{thm}
For example 
\[
{\psi}_p^{(r)}\left(\frac{1}{2}\right)=
-L_p(1+r,\omega^{-r})+2^{r+1}\ln_{r+1}^{(p)}(-1)=(1-2^{r+1})L_p(1+r,\omega^{-r}).
\]
Compare this with \cite[5.15.3]{NIST} the formula on the complex analytic polygamma functions.
\begin{pf}
We may assume $n>0$.
Let $s>0$ be an integer such that $p^s\equiv 1$ mod $N$.
Write $p^s-1=lN$.
\begin{align*}
S:=\sum_{\ve\in \mu_N\setminus\{1\}}(1-\ve^{-n})\ln_{r+1}^{(p)}(\ve)
&\os{\eqref{polylog-defn}}{\equiv}
\sum_{1\leq k<p^s,\, p\nmid k}\left(
\sum_{\ve\in \mu_N\setminus\{1\}}\frac{1-\ve^{-n}}{1-\ve^{p^s}}\frac{\ve^k}{k^{r+1}}\right)
\\
&=
\sum_{1\leq k<p^s,\, p\nmid k}\left(
\sum_{\ve\in \mu_N\setminus\{1\}}\frac{\ve^k+\cdots+\ve^{k+N-n-1}}{k^{r+1}}\right)
\end{align*}
modulo $p^s$.
Note $\sum_{\ve\in \mu_N\setminus\{1\}}\ve^i=N-1$ if $N|i$ and $=-1$ otherwise.
Let
$I$ be the set of integers $k$ satisfying that
$0\leq k<p^s$, $p\not{\hspace{-0.4mm}|}\,k$ and that
there is an integer $0\leq  i<N-n$ such 
that
$k+i\equiv 0$ mod $N$.
Then we have
\[
S\equiv 
\sum_{k\in I}
\frac{N}{k^{r+1}}\mod p^{s-1}
\]
by Lemma \ref{lemma.equiv}.
Hence 
\begin{align*}
&N^rS
\equiv 
\sum_{k\in I}\frac{1}{(k/N)^{r+1}}=
\sum_{k\equiv 0\text{ mod } N}
+\sum_{k\equiv -1\text{ mod } N}+\cdots+\sum_{k\equiv n-N+1\text{ mod }N}\\
&=
\sum_{\us{j\not\equiv 0\text{ mod }p}{1\leq j<p^s/N}}\frac{1}{j^{r+1}}
+\sum_{\us{j-1/N\not\equiv 0\text{ mod }p}{1\leq j<(p^s+1)/N}}\frac{1}{(j-1/N)^{r+1}}
+\cdots+
\sum_{\us{j-(N-n-1)/N\not\equiv 0\text{ mod }p}{1\leq j<(p^s+N-n-1)/N}}
\frac{1}{(j-(N-n-1)/N)^{r+1}}\\
&\equiv
\sum_{\us{j\not\equiv 0\text{ mod }p}{1\leq j\leq l}}\frac{1}{j^{r+1}}
+\sum_{\us{j+l\not\equiv 0\text{ mod }p}{1\leq j\leq l}}\frac{1}{(j+l)^{r+1}}
+\cdots+
\sum_{\us{j+l(N-n-1)\not\equiv 0\text{ mod }p}{1\leq j\leq l}}\frac{1}{(j+l(N-n-1))^{r+1}}\\
&=\sum_{\us{j\not\equiv 0\text{ mod }p}{1\leq j\leq l(N-n)}}\frac{1}{j^{r+1}}
=\sum_{\us{j\not\equiv 0\text{ mod }p}{0\leq j<l(N-n)+1}}\frac{1}{j^{r+1}}.
\end{align*}
Since $l(N-n)+1\equiv n/N$ mod $p^s$, the last summation is equivalent to $\wt\psi^{(r)}(n/N)$
modulo $p^{s-1}$ by definition. 
\end{pf}
\begin{rem}\label{definition-rem}
The complex analytic analogy of 
Theorem \ref{polygamma-thm2} is the following.
Let $\ln_r(z)=\ln^{an}_{r}(z)=\sum_{n=1}^\infty z^n/n^r$ be the analytic polylog.
Then
\begin{align*}
N^r\sum_{k=1}^{N-1}(1-e^{-2\pi ikn/N})\ln_{r+1}(e^{2\pi ik/N})&=
\sum_{m=1}^\infty\sum_{k=1}^{N-1}\frac{N^r}{m^{r+1}}(e^{2\pi ikm/N}-e^{2\pi ik(m-n)/N})\\
&=
\sum_{k=1}^\infty\frac{N^{r+1}}{(kN)^{r+1}}-\frac{N^{r+1}}{(kN-N+n)^{r+1}}\\
&=
\sum_{k=1}^\infty\frac{1}{k^{r+1}}-\frac{1}{(k-1+n/N)^{r+1}}.
\end{align*}
If $r=0$, then this is equal to $\psi(z)-\psi(1)$ (\cite[5.7.6]{NIST}).
If $r\geq 1$, then this is equal to $\zeta(r+1)+(-1)^r/r!\psi^{(r)}(n/N)$
(\cite[5.15.1]{NIST}).
\end{rem}
\begin{thm}\label{polygamma-thm4}
Let $m\geq 1$ be an positive integer prime to $p$. 
\begin{enumerate}
\item[\rm(1)]
Let $\psi_p(z)=\psi_p^{(0)}(z)$ be the $p$-adic digamma function.
Then
\[
\psi_p(mz)-\log^{(p)}(m)=\frac{1}{m}\sum_{i=0}^{m-1}\psi_p(z+\frac{i}{m}),
\]
(see \eqref{log(p)} for the definition of $\log^{(p)}(z)$).
\item[\rm(2)] If $r\ne 0$, we have
\[
\psi^{(r)}_p(mz)=\frac{1}{m^{r+1}}\sum_{i=0}^{m-1}\psi^{(r)}_p(z+\frac{i}{m}).
\]
\end{enumerate}
\end{thm}
See \cite[5.15.7]{NIST} for the corresponding formula on the complex analytic polygamma functions.
\begin{pf}
By Theorem \ref{p-zeta-def} (and Lemma \ref{lemma.log} in case $r=0$), 
the assertions are equivalent to
\begin{equation}\label{polygamma-thm4-eq}
\frac{1}{m^{r+1}}\sum_{i=0}^{m-1}\wt\psi^{(r)}_p(z+\frac{i}{m})
=\wt\psi^{(r)}_p(mz)
+\sum_{\ve\in\mu_N\setminus\{1\}}\ln^{(p)}_{r+1}(\ve)
\end{equation}
for all $r\in\Z$.
Since $\Z_{(p)}\cap[0,1)$ is a dense subset in $\Z_p$, it is enough to show the above
in case $z=n/N$ with $0\leq n<N$, $p\not|N$. 
By Theorem \ref{polygamma-thm2},
\begin{align*}
\frac{1}{m^{r+1}}\sum_{i=0}^{m-1}\wt\psi^{(r)}_p(z+\frac{i}{m})
&=\frac{1}{m^{r+1}}\sum_{i=0}^{m-1}\wt\psi^{(r)}_p(\frac{nm+iN}{mN})\\
&=\frac{N^r}{m}\sum_{i=0}^{m-1}\sum_{\nu\in\mu_{mN}\setminus\{1\}}
(1-\nu^{-nm-iN})\ln^{(p)}_{r+1}(\nu).
\end{align*}
The last summation is divided into the following two terms
\[
\sum_{i=0}^{m-1}
\sum_{\nu\in\mu_{N}\setminus\{1\}}
(1-\nu^{-nm})\ln^{(p)}_{r+1}(\nu)=
m\sum_{\nu\in\mu_{N}\setminus\{1\}}
(1-\nu^{-nm})\ln^{(p)}_{r+1}(\nu),
\]
\begin{align*}
\sum_{i=0}^{m-1}\sum_{\ve\in\mu_{m}\setminus\{1\}}\sum_{\nu^N=\ve}
(1-\nu^{-nm}\ve^{-i})\ln^{(p)}_{r+1}(\nu)
&=m\sum_{\ve\in\mu_{m}\setminus\{1\}}\sum_{\nu^N=\ve}
\ln^{(p)}_{r+1}(\nu)\\
&=\frac{m}{N^r}\sum_{\ve\in\mu_{m}\setminus\{1\}}
\ln^{(p)}_{r+1}(\ve)
\end{align*}
where the last equality follows from the distribution formula \eqref{polylog-dstr}.
Since the former is equal to $\wt\psi^{(r)}_p(nm/N)$ by Theorem \ref{polygamma-thm2},
the equality \eqref{polygamma-thm4-eq} follows.
\end{pf}
The relation with Morita's $p$-adic Gamma function $\Gamma_p(z)$ (e.g. \cite[11.6]{Cohen}) is as follows.
\begin{thm}\label{log.beta}
Let $B_p(x,y)=\Gamma_p(x)\Gamma_p(y)/\Gamma_p(x+y)$. Then
\[
\log B_p(z,q)=\sum_{i=1}^\infty\wt\psi^{(i-1)}_p(z)\frac{(-1)^{i}q^i}{i}
\]
for $z\in \Z_p$ and $q\in p\Z_p$
where $\log:1+p\Z_p\to p\Z_p$ is the Iwasawa logarithm.
\end{thm}
\begin{pf}
Fix $q\in p\Z_p$.
The functions
\[
\Z_p\lra p\Z_p,\quad z\longmapsto \log B_p(z,q)
\]
and
\[
\Z_p\lra p\Z_p,\quad z\longmapsto
\sum_{i=1}^\infty\wt\psi^{(i-1)}_p(z)\frac{(-q)^i}{i}
\]
are continuous.
Therefore it is enough to show
\begin{equation}\label{log.beta.eq1}
\log B_p(n,q)=\sum_{i=1}^\infty\wt\psi^{(i-1)}_p(n)\frac{(-q)^i}{i}
\end{equation}
for all $n\in\Z_{\geq0}$ as the set $\Z_{\geq0}$ is dense in $\Z_p$. 
We show it by induction on $n$. The case $n=0$ is trivial.
Suppose that \eqref{log.beta.eq1} is true for $n$.
If $p|n$, then $B_p(n+1,q)=B_p(n,q)$ (\cite[11.6.8. (3)]{Cohen})
and $\wt\psi^{(j)}_p(n+1)
=\wt\psi^{(j)}_p(n)$ (Theorem \ref{polygamma-thm1} (3)), so that
\eqref{log.beta.eq1} is true for $n+1$.
If $p\nmid n$, then
\begin{align*}
\log B_p(n+1,q)
&=\log \frac{\Gamma_p(n+1)\Gamma_p(q)}{\Gamma_p(n+1+q)}\\
&=\log \left(\frac{-n}{-(n+q)}\frac{\Gamma_p(n)\Gamma_p(q)}{\Gamma_p(n+q)}\right)
&\text{(\cite[11.6.8. (3)]{Cohen})}\\
&=-\log \left(1+\frac{q}{n}\right)+\log B_p(n,q)\\
&=\sum_{i=1}^\infty\frac{1}{n^i}\frac{(-q)^i}i+\sum_{i=1}^\infty
\wt\psi^{(i-1)}_p(n)\frac{(-q)^i}{i}\\
&=
\sum_{i=1}^\infty
\wt\psi^{(i-1)}_p(n+1)\frac{(-q)^i}{i},&\text{(Theorem \ref{polygamma-thm1} (3))}
\end{align*}
so that \eqref{log.beta.eq1} is true for $n+1$.
This completes the proof.
\end{pf}

\subsection{$p$-adic measure}\label{measure-sect}
For a function $g:\Z_p\to \C_p$, the Volkenborn integral is defined by
\[
\int_{\Z_p}g(t)dt=\lim_{s\to\infty}\frac{1}{p^s}\sum_{0\leq j<p^s}g(j)
\]
if the limit exists. We refer \cite[11.1.2]{Cohen} for a general theory on Volkenborn integrals
\begin{thm}\label{meas-thm2}
Let $\log:\C_p^\times\to \C_p$ be the Iwasawa logarithmic function.
Let 
\[
{\mathbf 1}_{\Z_p^\times}(z):=\begin{cases}
1&z\in\Z_p^\times\\
0&z\in p\Z_p
\end{cases}
\]
be the characteristic function.
Then
\[
\psi_p(z)=\int_{\Z_p}\log(z+t){\mathbf 1}_{\Z_p^\times}(z+t)dt.
\]
\end{thm}
\begin{pf}
Using the computation in \eqref{equiv-log}, one can easily show that
the Volkenborn integral 
$Q(z):=\int_{\Z_p}{\mathbf 1}_{\Z_p^\times}(z+t)\log(z+t)dt$ is defined, and it is continuous
with respect to $z\in\Z_p$. 
Moreover we have
\[
Q(z+1)-Q(z)\equiv\begin{cases}
p^{-s}(\log(z)-\log(z+p^s))& z\in \Z_p^\times\\
0&z\in p\Z_p
\end{cases}\mod p^{s'}
\]
where $s'=s-1$ if $p=2$ and $s'=s$ if $p\geq 3$.
For $z\in \Z_p^\times$, since
\[
p^{-s}(\log(z)-\log(z+p^s))=-p^{-s}\log(1+z^{-1}p^s)
\equiv z^{-1}\mod p^{s'},
\]
it follows from Theorem \ref{polygamma-thm1} (3) that
 $Q(z)$ differs from $\psi_p(z)$ by a constant.
Since
\[
Q(0)=\lim_{s\to\infty} \frac{1}{p^s}\sum_{0\leq j<p^s,p\nmid j}\log(j)= -\gamma_p=\psi_p(0),
\]
we obtain $Q(z)=\psi_p(z)$.
\end{pf}

\begin{thm}\label{meas-thm1}
If $r\ne0$, then
\[
\psi^{(r)}_p(z)=-\frac{1}{r}\int_{\Z_p}(z+t)^{-r}{\mathbf 1}_{\Z_p^\times}(z+t)dt
\]
where ${\mathbf 1}_{\Z_p^\times}(z)$ denotes the characteristic function as in Theorem \ref{meas-thm2}.
\end{thm}
\begin{pf}
Using Lemma \ref{lemma.equiv}, one sees that
the Volkenborn integral 
$Q(z)=-\frac{1}{r}\int_{\Z_p}(z+t)^{-r}{\mathbf 1}_{\Z_p^\times}(z+t)dt$ is defined.
Moreover 
if $z\in\Z_p^\times$, then
\[
Q(z+1)-Q(z)\equiv\frac{-1}{rp^s}\left(
\frac{1}{(z+p^s)^r}-\frac{1}{z^r}\right)\equiv z^{-1-r}\mod p^{s-\ord_p(r)},
\]
and if $z\in p\Z_p$, then $Q(z+1)\equiv Q(z)$.
This shows that $Q(z)-\psi^{(r)}_p(z)$ is a constant by Theorem \ref{polygamma-thm1} (3). 
We show $Q(0)=\psi^{(r)}_p(0)$.
By definition 
\[
Q(0)= \lim_{n\to\infty}\frac{-1}{rp^n}\sum_{0\leq k<p^n,p\nmid k}
\frac{1}{k^r}.
\]
Recall the original definition of
the $p$-adic $L$-function by Kubota-Leopoldt \cite{KL}
\[
L_p(s,\chi)=\frac{1}{s-1}\lim_{n\to\infty}\frac{1}{fp^n}
\sum_{0\leq k<fp^{n},p\nmid k}\chi(k)\langle k\rangle^{1-s},
\quad \langle k\rangle:=k/\omega(k)
\]
for a primitive Dirichlet character $\chi$ with conductor $f\geq1$.
This immediately implies $Q(0)=-L_p(1+r,\omega^{-r})=\psi^{(r)}_p(0)$.
\end{pf}
\section{$p$-adic hypergeometric functions of logarithmic type}\label{pHGlog-sect}
We write the Pochhammer symbol by $(a)_n$,
\[
(a)_0:=1,\quad (a)_n:=a(a+1)\cdots(a+n-1),\,n\geq 1.
\]
For $a\in \Z_p$, the {\it Dwork prime} $a'$ is defined to be $(a+l)/p$ 
where $l\in \{0,1,\ldots,p-1\}$ is the unique integer such that $a+l\equiv 0$ mod $p$.
The $i$-th Dwork prime is denoted by 
$a^{(i)}$ which is defined to be $(a^{(i-1)})'$ with $a^{(0)}=a$.
\subsection{Definition}\label{pHGlog-defn}
Let $s\geq 1$ be a positive integer.
Let $a_i,b_j\in\Q_p$ with $b_j\not\in \Z_{\leq0}$. Let
\[
{}_sF_{s-1}\left({a_1,\ldots,a_s\atop b_1,\ldots b_{s-1}}:t\right)
=\sum_{n=0}^\infty\frac{(a_1)_n\cdots(a_s)_n}{(b_1)_n\cdots(b_{s-1})_n}\frac{t^n}{n!}.
\]
be the {\it hypergeometric power series} with coefficients.
In what follows we only consider the cases $a_i\in \Z_p$ and $b_j=1$, and then 
we write
\[
F_{\ul a}(t):={}_sF_{s-1}\left({a_1,\ldots,a_s\atop 1,\ldots 1}:t\right)\in \Z_p[[t]]
\]
for 
$\ul a=(a_1,\ldots,a_s)\in \Z_p^s$.
\begin{defn}[$p$-adic hypergeometric functions of logarithmic type]
Write $\ul a^{(i)}=(a^{(i)}_1,\ldots,a^{(i)}_s)$ where $(-)^{(i)}$ denotes the $i$-th
Dwork prime.
Let $W=W(\ol\F_p)$ be the Witt ring of $\ol\F_p$.
Let $\sigma:W[[t]]\to W[[t]]$ be the $p$-th Frobenius endomorphism given by  
$\sigma(t)=ct^p$ with $c\in 1+pW$,
compatible with the Frobenius on $W$.
Put a power series
\[
G_{\ul a}(t):=\psi_p(a_1)+\cdots+\psi_p(a_s)+s\gamma_p
-p^{-1}\log(c)+
\int_0^t(F_{\underline{a}}(t)-F_{\ul a^{(1)}}(t^\sigma))\frac{dt}{t}
\]
where $\psi_p(z)$ is the $p$-adic digamma function defined in \S \ref{polygamma-sect},
and $\log(z)$ is the Iwasawa logarithmic function.
Then we define
\[
\cF^{(\sigma)}_{\underline{a}}(t)=G_{\ul a}(t)/F_{\underline{a}}(t),
\]
and call the {\rm $p$-adic hypergeometric functions of logarithmic type}.
\end{defn}
\begin{lem}\label{G-int}
$G_{\ul a}(t)\in W[[t]]$. Hence it follows $\cF^{(\sigma)}_{\underline{a}}(t)\in W[[t]]$.
\end{lem}
\begin{pf}
Let $G_{\underline{a}}(t)=\sum B_it^i$.
Let $F_{\underline{a}}(t)=\sum A_it^i$ and
$F_{\ul a^{(1)}}(t)=\sum \wA_it^i$.
If $p{\not|}i$, then
$B_i=A_i/i$ is obviously a $p$-adic integer.
For $i=mp^k$ with $k\geq 1$ and $p\not| m$, one has
\[
B_i=B_{mp^k}=\frac{A_{mp^k}-c^{mp^{k-1}}\wA_{mp^{k-1}}}{mp^k}.
\]
Since $c^{mp^{k-1}}\equiv 1$ mod $p^k$, it is enough to see 
$A_{mp^k}\equiv \wA_{mp^{k-1}}$ mod $p^k$. 
However this follows from \cite[p.36, Cor. 1]{Dwork-p-cycle}. 
\end{pf}
\subsection{Congruence relations}\label{cong-pf-sect1}
For a power series $f(t)=\sum_{n=0}^\infty A_nt^n$, we write $[f(t)]_{<m}:=\sum_{n<m}A_nt^n$
the truncated polynomial.
\begin{thm}\label{cong-thm}
Suppose that $a_i\not\in \Z_{\leq0}$ for all $i$.
Let us write $\cF^{(\sigma)}_{\underline{a}}(t)=G_{\ul a}(t)/F_{\ul a}(t)$.
If $c\in 1+2pW$, then for all $n\geq 1$
\begin{equation}\label{cong-thm-eq0}
\cF^{(\sigma)}_{\underline{a}}(t)\equiv
\frac{[G_{\underline{a}}(t)]_{<p^n}}{[F_{\underline{a}}(t)]_{<p^n}}\mod p^nW[[t]].
\end{equation}
If $p=2$ and $c\in 1+2W$ (not necessarily $c\in 1+4W$), 
then the above holds modulo $p^{n-1}$.
\end{thm}
\begin{cor}\label{cong-cor}
Suppose that there exists an integer $r\geq 0$ such that $a_i^{(r+1)}=a_i$ for all $i$
where $(-)^{(r)}$ denotes the $r$-th Dwork prime.
Then
\[
\cF_{\ul a}^{(\sigma)}(t)\in 
W\langle t,[F_{\ul a}(t)]_{<p}^{-1}, \ldots,[F_{\ul a^{(r)}}(t)]_{<p}^{-1}\rangle:=
\varprojlim_n (W/p^n[t,[F_{\ul a}(t)]_{<p}^{-1}, \ldots,[F_{\ul a^{(r)}}(t)]_{<p}^{-1}])
\]
is a convergent function.
For $\alpha\in W$ such that $[F_{\ul a^{(i)}}(\alpha)]_{<p}\not
\equiv 0$ mod $p$ for all $i$, the special value of $\cF_{\ul a}^{(\sigma)}(t)$ at $t=\alpha$
is defined, and it is explicitly given by
\[
\cF_{\ul a}^{(\sigma)}(\alpha)=\lim_{n\to\infty}
\frac{[G_{\underline{a}}(\alpha)]_{<p^n}}{[F_{\underline{a}}(\alpha)]_{<p^n}}.
\]
\end{cor}
\subsection{Proof of Congruence relations : Reduction to the case $c=1$}\label{cong-sect2}
Throughout the sections \ref{cong-sect2}, \ref{cong-sect3} and \ref{cong-sect4},
we use the following notation. Fix $s\geq 1$ and $\ul a=(a_1,\ldots,a_s)$ with $a_i\not\in \Z_{\leq0}$. Let $\sigma(t)=ct^p$ be the Frobenius. 
Put
\begin{equation}\label{tilde-defn}
A_n:=\frac{(a_1)_n}{n!}\cdots\frac{(a_s)_n}{n!},\quad
\wA_n:=\frac{(a^{(1)}_1)_n}{n!}\cdots\frac{(a_s^{(1)})_n}{n!}
\end{equation}
for $n\geq 0$.
Let $B_n$ be defined by
$G_{\ul a}(t)=\sum_{n=0}^\infty B_nt^n$, or explicitly
\begin{equation}\label{cong-sect1-eq1-0}
B_0=\psi_p(a_1)+\cdots+\psi_p(a_s)+s\gamma_p,
\end{equation}
\begin{equation}\label{cong-sect1-eq1}
B_n=\frac{A_n}{n},\,(p\nmid n),\quad
B_{mp^k}=\frac{A_{mp^k}-c^{mp^{k-1}}\wA_{mp^{k-1}}}{mp^k},\,(m,k\geq1).
\end{equation}

\begin{lem}\label{cong-lem0}
The proof of Theorem \ref{cong-thm} is reduced to the case $\sigma(t)=t^p$ (i.e. $c=1$).
\end{lem}
\begin{pf}
Write $[f(t)]_{\geq m}:=f(t)-[f(t)]_{<m}$. 
Put $n^*:=n$ if $c\in 1+2pW$ and $n^*=n-1$ if $p=2$ and $c\not
\in 1+4W$.
Theorem \ref{cong-thm}
is equivalent to saying
\[
F_{\ul a}(t)\cdot[G_{\ul a}(t)]_{\geq p^n}
\equiv [F_{\ul a}(t)]_{\geq p^n}\cdot G_{\ul a}(t)\mod p^{n^*}W[[t]],
\]
namely
\[
\sum_{i+j=m,\,i,j\geq0}A_{i+p^n}B_j-A_iB_{j+p^n}\equiv 0\mod p^{n^*}
\]
for all $m\geq 0$.
Suppose that this is true when $c=1$, namely
\begin{equation}\label{cong-lem0-eq1}
\sum_{i+j=m}A_{i+p^n}B^\circ_j-A_iB^\circ_{j+p^n}\equiv 0\mod p^{n^*}
\end{equation}
where $B^\circ_i$ are the coefficients \eqref{cong-sect1-eq1-0}
or \eqref{cong-sect1-eq1} when $c=1$.
Suppose that $c\in 1+pW$ is an arbitrary element, and
let $B_i$ be as in \eqref{cong-sect1-eq1-0}
or \eqref{cong-sect1-eq1}.
We then want to show
\begin{equation}\label{cong-lem0-eq2}
\sum_{i+j=m}A_{i+p^n}(B^\circ_j-B_j)-A_i(B^\circ_{j+p^n}-B_{j+p^n})\equiv 0\mod p^{n^*}.
\end{equation}
Let $c=1+pe$ with $e\ne0$ (if $e=0$, there is nothing to prove). 
Then
\begin{align*}
\sum_{i+j=m}A_{i+p^n}(B^\circ_j-B_j)
&=A_{m+p^n}p^{-1}\log(c)+\sum_{1\leq j\leq m} p^{-1}\frac{(c^{j/p}-1)A_{m+p^n-j}
\wA_{j/p}}{j/p}\\
&=A_{m+p^n}\sum_{i=1}^\infty\frac{(-1)^{i+1}}{i}p^{i-1}e^i+\sum_{1\leq j\leq m}(j/p)^{-1}
\sum_{i=1}^\infty\binom{j/p}{i}p^{i-1}e^iA_{m+p^n-j}\wA_{j/p}\\
&=\sum_{i=1}^\infty\left(A_{m+p^n}\frac{(-1)^{i+1}}{i}+\sum_{1\leq j\leq m}(j/p)^{-1}
\binom{j/p}{i}A_{m+p^n-j}\wA_{j/p}\right)p^{i-1}e^i\\
&=\sum_{i=1}^\infty\left(A_{m+p^n}\frac{(-1)^{i+1}}{i}+\sum_{1\leq j\leq m}i^{-1}
\binom{j/p-1}{i-1}A_{m+p^n-j}\wA_{j/p}\right)p^{i-1}e^i\\
&=\sum_{i=1}^\infty\left(\sum_{0\leq j\leq m}i^{-1}
\binom{j/p-1}{i-1}A_{m+p^n-j}\wA_{j/p}\right)p^{i-1}e^i
\end{align*}
where we always mean $A_{j/p}=\wA_{j/p}=0$ unless $p|j$.
Similarly
\[
\sum_{i+j=m}A_i(B^\circ_{j+p^n}-B_{j+p^n})
=\sum_{i=1}^\infty\left(\sum_{0\leq j\leq m}i^{-1}
\binom{(m+p^n-j)/p-1}{i-1}A_{j}\wA_{(m+p^n-j)/p}\right)p^{i-1}e^i.
\]
Therefore it is enough to show that
\[
\frac{p^{i-1}e^i}{i}\sum_{0\leq j\leq m}
\binom{j/p-1}{i-1}A_{m+p^n-j}\wA_{j/p}
\equiv\frac{p^{i-1}e^i}{i}\sum_{0\leq j\leq m}
\binom{(m+p^n-j)/p-1}{i-1}A_{j}\wA_{(m+p^n-j)/p}\mod p^{n^*}
\]
equivalently
\begin{equation}\label{cong-lem0-eq3}
\sum_{0\leq j\leq m}
(1-j/p)_{i-1}
A_{m+p^n-j}\wA_{j/p}
\equiv\sum_{0\leq j\leq m}
(1-(m+p^n-j)/p)_{i-1}A_{j}\wA_{(m+p^n-j)/p}\mod p^{n^*-i+1}i!e^{-i}
\end{equation}
for all $i\geq 1$ and $m\geq 0$. 
Recall the Dwork congruence 
\begin{equation}\label{Dwork.Cong}
\frac{F_{\ul a^{(1)}}(t^p)}{F_{\ul a}(t)}\equiv
\frac{[F_{\ul a^{(1)}}(t^p)]_{<p^m}}{[F_{\ul a}(t)]_{<p^m}}\mod p^l\Z_p[[t]],\quad m\geq l
\end{equation}
from \cite[p.37, Thm. 2, p.45]{Dwork-p-cycle}.
This immediately imples \eqref{cong-lem0-eq3} in case $i=1$.
Suppose $i\geq 2$.
To show \eqref{cong-lem0-eq3}, it is enough to show 
\begin{equation}\label{cong-lem0-eq4}
\sum_{0\leq j\leq m}
(j/p)^k
A_{m+p^n-j}\wA_{j/p}
\equiv\sum_{0\leq j\leq m}
((m+p^n-j)/p)^kA_{j}\wA_{(m+p^n-j)/p}\mod p^{n^*-i+1}i!e^{-i}
\end{equation}
for each $k\geq 0$.
Let $D=t\frac{d}{dt}$, and put $F^*(t):=D^kF_{\ul a^{(1)}}(t)=\sum_{j=0}^\infty j^k\wA_jt^j$.
Then \eqref{cong-lem0-eq4} is equivalent to saying 
\begin{equation}\label{cong-lem0-eq5}
[F_{\ul a}(t)]_{<p^n}F^*(t^p)\equiv
F_{\ul a}(t)[F^*(t^p)]_{<p^n}\mod p^{n^*-i+1}i!e^{-i}\Z_p[[t]].
\end{equation}
We show \eqref{cong-lem0-eq5}, which finishes the proof of Lemma \ref{cong-lem0}.
Using the Dwork congruence \eqref{Dwork.Cong}, one can show
\[
\frac{\frac{d^k}{dt^k}F_{\ul a^{(1)}}(t)}{F_{\ul a^{(1)}}(t)}
\equiv
\frac{\frac{d^k}{dt^k}([F_{\ul a^{(1)}}(t)]_{<p^m})}{[F_{\ul a^{(1)}}(t)]_{<p^m}}
\mod p^m\Z_p[[t]]
\]
in the same way as the proof of \cite[p.45, (3.14)]{Dwork-p-cycle}.
Hence
\[
\frac{F^*(t)}{F_{\ul a^{(1)}}(t)}\equiv
\frac{[F^*(t)]_{<p^m}}{[F_{\ul a^{(1)}}(t)]_{<p^m}}\mod p^l\Z_p[[t]],\quad m\geq l
\]
and this implies
\[
\frac{F^*(t^p)}{F_{\ul a^{(1)}}(t^p)}\equiv
\frac{[F^*(t^p)]_{<p^n}}{[F_{\ul a^{(1)}}(t^p)]_{<p^n}}\mod p^{n-1}\Z_p[[t]].
\]
Therefore we have
\[
\frac{F^*(t^p)}{F_{\ul a}(t)}=\frac{F_{\ul a^{(1)}}(t^p)}{F_{\ul a}(t)}
\frac{F^*(t^p)}{F_{\ul a^{(1)}}(t^p)}\equiv
\frac{[F_{\ul a^{(1)}}(t^p)]_{<p^n}}{[F_{\ul a}(t)]_{<p^n}}
\frac{[F^*(t^p)]_{<p^n}}{[F_{\ul a^{(1)}}(t^p)]_{<p^n}}
=\frac{[F^*(t^p)]_{<p^n}}{[F_{\ul a}(t)]_{<p^n}}
\mod p^{n-1}\Z_p[[t]].
\]
If $p\geq 3$, then $\ord_p(p^{n^*-i+1}i!)=\ord_p(p^{n-i+1}i!)\leq n-1$ for any $i\geq2$, and hence
\eqref{cong-lem0-eq5} follows.
If $p=2$, then $\ord_p(p^{n-i+1}i!)\leq n$ but not necessarily $\ord_p(p^{n-i+1}i!)\leq n-1$.
If $e\in 2W\setminus\{0\}$, then $\ord_p(p^{n^*-i+1}i!e^{-i})=\ord_p(p^{n-i+1}i!e^{-i})\leq n-i<n-1$,
and hence
\eqref{cong-lem0-eq5} follows.
If $e$ is a unit, then 
$\ord_p(p^{n^*-i+1}i!e^{-i})=\ord_p(p^{n-i}i!)\leq n-1$ for any $i\geq2$, as required again.
This completes the proof.
\end{pf}
\subsection{Proof of Congruence relations : Preliminary lemmas}\label{cong-sect3}
Until the end of \S \ref{cong-sect4}, let $\sigma$ be the Frobenius given by $\sigma(t)=t^p$
(i.e. $c=1$).
Therefore
\begin{equation}\label{cong-B}
B_0=\psi_p(a_1)+\cdots+\psi_p(a_s)+s\gamma_p,\quad
B_i=\frac{A_i-\wA_{i/p}}{i},\quad i\in\Z_{\geq1}
\end{equation}
where the notation is as in \eqref{tilde-defn} and  we 
always mean $A_{i/p}=\wA_{i/p}=0$ unless $p| i$.
\begin{lem}\label{prod-lem}
For an $p$-adic integer $\alpha\in \Z_p$ and $n\in\Z_{\geq1}$, we define 
\[
\{\alpha\}_n:=\prod_{\underset{p\nmid (a+i-1)}{1\leq i\leq n}}(\alpha+i-1),
\]
and $\{\alpha\}_0:=1$.
Let $a\in\Z_p\setminus\Z_{\leq 0}$, and 
$l\in \{0,1,\ldots,p-1\}$ the integer such that $a\equiv -l$ mod $p$.
Then for any $m\in \Z_{\geq0}$, we have
\[
m\equiv 0,1,\ldots,l\text{ \rm mod }p\quad\Longrightarrow\quad
\frac{(a)_m}{m!}\left(\frac{(a')_{\lfloor m/p\rfloor}}
{\lfloor m/p\rfloor !}\right)^{-1}=\frac{\{a\}_{m}}{\{1\}_{m}}\in\Z_p^\times,
\]
\[
m\equiv l+1,\ldots,p-1\text{ \rm mod }p\quad\Longrightarrow\quad
\frac{(a)_m}{m!}\left(\frac{(a')_{\lfloor m/p\rfloor}}
{\lfloor m/p\rfloor !}\right)^{-1}=\left(a+l+p\lfloor\frac{m}{p}\rfloor\right)
\frac{\{a\}_m}{\{1\}_m}
\]
where $a'=a^{(1)}$ is the Dwork prime.
\end{lem}
\begin{pf}
Straightforward.
\end{pf}
\begin{lem}[Dwork]\label{cong-lem5}
For any $m\in\Z_{\geq0}$, $A_m/\wA_{\lfloor m/p\rfloor}$ are $p$-adic integers, and
\[
m_1\equiv m_2\mod p^n\quad\Longrightarrow\quad
A_{m_1}/\wA_{\lfloor m_1/p\rfloor}\equiv
A_{m_2}/\wA_{\lfloor m_2/p\rfloor}\mod p^n.
\]
\end{lem}
\begin{pf}
This is \cite[p.36, Cor. 1]{Dwork-p-cycle}, or one can easily show this
by using Lemma \ref{prod-lem} on noticing the fact that
\[
\{\alpha\}_{p^n}
\equiv \prod_{i\in (\Z/p^n\Z)^\times}i
\equiv\begin{cases}
1&p=2,\,n\ne 2\\
-1&\text{otherwise}
\end{cases}
\mod p^n.
\]
\end{pf}

\begin{lem}\label{b0-lem}
Let $a\in\Z_p\setminus\Z_{\leq 0}$ and $m,n\in \Z_{\geq1}$. Then
\begin{equation}\label{b0-lem-eq1}
1-\frac{(a')_{mp^{n-1}}}{(mp^{n-1})!}
\left(\frac{(a)_{mp^n}}{(mp^n)!}\right)^{-1}\equiv mp^n
(\psi_p(a)+\gamma_p)\mod p^{2n}.
\end{equation}
Moreover $\wA_{mp^{n-1}}/A_{mp^n}$
and $B_k/A_k$ 
are $p$-adic integers for all $k,m\geq 0$, $n\geq 1$, and
\begin{equation}\label{b0-lem-eq2}
\frac{\wA_{mp^{n-1}}}{A_{mp^n}}\equiv1-mp^n(\psi_p(a_1)+\cdots+\psi_p(a_s)+s\gamma_p)
\mod p^{2n},
\end{equation}
\begin{equation}\label{b0-lem-eq3}
p\nmid m\quad \Longrightarrow\quad\frac{B_{mp^n}}{A_{mp^n}}=
\frac{1-\wA_{mp^{n-1}}/A_{mp^n}}{mp^n}\equiv B_0\mod p^n.
\end{equation}
\end{lem}
\begin{pf}
We already see that $\wA_{mp^{n-1}}/A_{mp^n}\in \Z_p$ in Lemma \ref{prod-lem}.
It is enough to show \eqref{b0-lem-eq1}.
Indeed \eqref{b0-lem-eq2} is immediate from \eqref{b0-lem-eq1}, and
\eqref{b0-lem-eq3} is immediate from \eqref{b0-lem-eq2}.
Moreover \eqref{b0-lem-eq2} also implies that $B_k/A_k\in \Z_p$ for any $k\in\Z_{\geq0}$.

Let us show \eqref{b0-lem-eq1}.
Let $a=-l+p^nb$ with $l\in \{0,1,\ldots,p^n-1\}$.
Then
\[
\frac{(a')_{mp^{n-1}}}{(mp^{n-1})!}
\left(\frac{(a)_{mp^n}}{(mp^n)!}\right)^{-1}
=\frac{\{1\}_{mp^n}}{\{a\}_{mp^n}}=\prod_{\underset{k-l\not\equiv0\text{ mod }p}{l<k<mp^n}} \frac{k-l}{k-l+p^nb}
\times\prod_{\underset{k-l\not\equiv0\text{ mod }p}{0\leq k<l}} \frac{k-l+mp^n}{k-l+p^nb}
\]
by Lemma \ref{prod-lem}.
If $(p,n)\ne(2,1)$, we have
\begin{align*}
\frac{\{1\}_{mp^n}}{\{a\}_{mp^n}}
&\equiv
\prod_{\underset{p\nmid k-l}{l<k<mp^n}}\left( 1-\frac{p^nb}{k-l}\right)
\prod_{\underset{p\nmid k-l}{0\leq k<l}}\left( 1-\frac{p^n(b-m)}{k-l}\right)
\equiv
1-p^n\left(\sum_{\underset{p\nmid k-l}{l<k<mp^n}}\frac{b}{k-l}
+\sum_{\underset{p\nmid k-l}{0\leq k<l}}\frac{b-m}{k-l}\right)\\
&\os{(\ast)}{\equiv}
1-mp^n\sum_{\underset{p\nmid k-l}{l<k<mp^n}}\frac{1}{k-l}
=
1-mp^n\sum_{1\leq k<mp^n-l,p \nmid k}\frac{1}{k}
\os{\eqref{equiv-psi}}{\equiv}1-mp^n(\psi_p(a)+\gamma_p)
\end{align*}
modulo $p^{2n}$, where $(\ast)$ follows from Lemma \ref{lemma.equiv}.
This completes the proof of \eqref{b0-lem-eq1} in case $(p,n)\ne(2,1)$.
In case $(p,n)=(2,1)$, we need another observation
(since the equivalence $(\ast)$ breaks down).
In this case
we have
\begin{align*}
\frac{\{1\}_{2m}}{\{a\}_{2m}}
&\equiv
1-2
\left(\sum_{\underset{2\nmid  k-l}{l<k<2m}}\frac{b}{k-l}
+\sum_{\underset{2\nmid k-l}{0\leq k<l}}\frac{b-m}{k-l}\right)\mod 4\\
&=
1-2
\left(\sum_{\underset{2\nmid  k-l}{l<k<2m}}\frac{m}{k-l}
+\sum_{\underset{2\nmid k-l}{0\leq k<2m}}\frac{b-m}{k-l}\right)\\
&
\equiv
1-2m
\left(\sum_{0<k<2m-l,\,2\nmid k}\frac{1}{k}
+b-m\right)\mod 4
\end{align*}
and
\[
\psi_2(a)+\gamma_2\equiv \sum_{0<k<L,\,2 \nmid k}\frac{1}{k}\mod 2
\]
where $L\in\{0,1,2,3\}$ such that $a=-l+2b\equiv L$ mod $4$.
Therefore \eqref{b0-lem-eq1} is equivalent to 
\[
m
\left(\sum_{0<k<2m-l,\,2 \nmid k}\frac{1}{k}
-\sum_{0<k<L,\,2 \nmid k}\frac{1}{k}
+b-m\right)\equiv 0\mod 2.
\]
We may assume that $m>0$ is odd and $b=0,1$ (hence $a=0,\pm 1,2$). 
Then one can check this
on a case-by-case analysis.
\end{pf}
\begin{lem}\label{cong-lem3}
For any $m_1, m_2\in\Z_{\geq 0}$ and $n\in \Z_{\geq 1}$, we have \[
m_1\equiv m_2\mod p^n\quad\Longrightarrow\quad
\frac{B_{m_1}}{A_{m_1}}\equiv
\frac{B_{m_2}}{A_{m_2}}\mod p^n.
\]
\end{lem}
\begin{pf}
If $p{\not|} m_i$, then $B_{m_i}/A_{m_i}=1/m_i$ and hence the assertion is obvious.
Let $m_1=kp^i$ with $i\geq 1$ and $p{\not|}k$. 
It is enough to show the assertion in case $m_2=m_1+p^n$.
If $n\leq i$, then
\[
\frac{B_{m_1}}{A_{m_1}}\equiv \frac{B_{m_2}}{A_{m_2}}\equiv B_0\mod p^n
\]
by \eqref{b0-lem-eq3}. Suppose $n\geq i$.
Notice that for $m\in p\Z_{\geq0}$
\[
1-m\frac{B_m}{A_m}=\frac{\wA_{m/p}}{A_m}=\prod_{r=1}^s \frac{\{1\}_m}{\{a_r\}_m}
\]
by \eqref{cong-B} and Lemma \ref{prod-lem}.
We have
\begin{align*}
1-m_2\frac{B_{m_2}}{A_{m_2}}
&=\prod_r \frac{\{1\}_{kp^i+p^n}}{\{a_r\}_{kp^i+p^n}}
=\prod_r \frac{\{1\}_{kp^i}}{\{a_r\}_{kp^i}}\frac{\{1+kp^i\}_{p^n}}{\{a_r+kp^i\}_{p^n}}
=\left(1-m_1\frac{B_{m_1}}{A_{m_1}}\right)\prod_r 
\frac{\{1+kp^i\}_{p^n}}{\{a_r+kp^i\}_{p^n}}\\
&=\left(1-m_1\frac{B_{m_1}}{A_{m_1}}\right)\prod_r 
\frac{\{1\}_{p^n}}{\{a_r+kp^i\}_{p^n}}
\frac{\{1+kp^i\}_{p^n}}{\{1\}_{p^n}}\\
&\os{(\ast)}{\equiv}\left(1-m_1\frac{B_{m_1}}{A_{m_1}}\right)\prod_r 
(1-p^n(\psi_p(a_r+kp^i)-\psi_p(1+kp^i)))\mod p^{2n}\\
&\os{(\ast\ast)}{\equiv}\left(1-m_1\frac{B_{m_1}}{A_{m_1}}\right)
(1-p^nB_0)\mod p^{n+i}.
\end{align*}
Here 
$(\ast)$ follows from Lemmas \ref{prod-lem} and \ref{b0-lem}.
The equivalence 
$(\ast\ast)$ follows from Theorem \ref{polygamma-thm1} (1) and 
\eqref{equiv-psi} in case $(p,i)\ne(2,1)$, and in case $(p,i)=(2,1)$,
it does from the fact that
\[
\psi_2(z+2)-\psi_2(z)\equiv1\mod 2,\quad z\in \Z_2.
\]
Therefore we have
\[
kp^i\left(\frac{B_{m_2}}{A_{m_2}}-\frac{B_{m_1}}{A_{m_1}}\right)
\equiv p^n\left(B_0-\frac{B_{m_2}}{A_{m_2}}\right)
\mod p^{i+n}.
\]
By \eqref{b0-lem-eq3}, the right hand side vanishes.
This is the desired assertion.
\end{pf}
\begin{lem}\label{cong-lem4}
Put $S_m:=\sum_{i+j=m}A_{i+p^n}B_{j}-A_iB_{j+p^n}$ for $m\in \Z_{\geq0}$.
Then \[
S_m\equiv\sum_{i+j=m}(A_{i+p^n}A_{j}-A_iA_{j+p^n})\frac{B_{j}}{A_{j}}\mod p^n.
\]
\end{lem}
\begin{pf}
\begin{align*}
S_m
&=\sum_{i+j=m}A_{i+p^n}B_{j}-A_iA_{j+p^n}\frac{B_{j+p^n}}{A_{j+p^n}}\\
&\equiv\sum_{i+j=m}A_{i+p^n}B_{j}-A_iA_{j+p^n}\frac{B_{j}}{A_{j}}\mod p^n\quad(\mbox{Lemma \ref{cong-lem3}})\\
&=\sum_{i+j=m}(A_{i+p^n}A_{j}-A_iA_{j+p^n})\frac{B_{j}}{A_{j}}
\end{align*}
as required.
\end{pf}
\begin{lem}\label{cong-lem6}
\[
S_m\equiv\sum_{i+j=m}
(\wA_{\lfloor j/p\rfloor}\wA_{\lfloor i/p\rfloor+p^{n-1}}
-\wA_{\lfloor i/p\rfloor}\wA_{\lfloor j/p\rfloor+p^{n-1}}
)\frac{A_i}{\wA_{\lfloor i/p\rfloor}}\frac{A_j}{\wA_{\lfloor j/p\rfloor}}
\frac{B_j}{A_j}\mod p^n.
\]
\end{lem}
\begin{pf}
This follows from Lemma \ref{cong-lem4} and Lemma \ref{cong-lem5}.
\end{pf}

\begin{lem}\label{cong-lem7}
For all $m,k,s\in \Z_{\geq0}$ and $0\leq l\leq n$, 
we have
\begin{equation}\label{cong-lem7-eq1}
\sum_{\underset{i\equiv k\text{ mod }p^{n-l}}{i+j=m}}A_iA_{j+p^{n-1}}-A_jA_{i+p^{n-1}}
\equiv0
\mod p^l.
\end{equation}
\end{lem}
\begin{pf}
There is nothing to prove in case $l=0$. If $l=n$, then \eqref{cong-lem7-eq1} is obvious as
\[
\mbox{LHS}=\sum_{i+j=m}A_iA_{j+p^{n-1}}-A_jA_{i+p^{n-1}}=0.
\]
Suppose $1\leq l\leq n-1$.
We simply write
\[
F^{(r)}(t)=F_{\ul a^{(r)}}(t)=\sum_{i=0}^\infty \frac{(a^{(r)}_1)_i}{i!}\cdots\frac{(a_s^{(r)})_i}{i!}t^i,
\quad F(t)=F^{(0)}(t).
\]
For $k\in \Z_{\geq0}$, we put
\begin{equation}\label{cong-lem7-eq3}
F^{(r)}_k(t):=\sum_{i\equiv k\text{ mod }p^{n-l}}\frac{(a^{(r)}_1)_i}{i!}
\cdots\frac{(a_s^{(r)})_i}{i!}t^i
=p^{-n+l}\sum_{s=0}^{p^{n-l}-1}\zeta^{-sk}F(\zeta^st)
\end{equation}
where $\zeta$ is a primitive $p^{n-l}$-th root of unity.
Then \eqref{cong-lem7-eq1} is equivalent to
\begin{equation}\label{cong-lem7-eq2}
F_k(t)\cdot[F_{m-k}(t)]_{<p^{n-1}}\equiv [F_k(t)]_{<p^{n-1}}\cdot F_{m-k}(t) \mod p^l.
\end{equation}
It follows from the Dwork congruence
\cite[p.37, Thm. 2]{Dwork-p-cycle} that one has
\[
\frac{F^{(i)}(t)}{F^{(i+1)}(t^p)}\equiv 
\frac{[F^{(i)}(t)]_{<p^m}}{[F^{(i+1)}(t^p)]_{<p^m}}\mod p^n
\]
for any $m\geq n\geq1$.
This implies
\[
\frac{F^{(i)}(t^p)}{F^{(i+1)}(t^{p^2})}\equiv 
\frac{[F^{(i)}(t^p)]_{<p^{n+1}}}{[F^{(i+1)}(t^{p^2})]_{<p^{n+1}}}\mod p^n,\quad
\frac{F^{(i)}(t^{p^2})}{F^{(i+1)}(t^{p^3})}\equiv 
\frac{[F^{(i)}(t^{p^2})]_{<p^{n+2}}}{[F^{(i+1)}(t^{p^3})]_{<p^{n+2}}}\mod p^n,\ldots.
\]
Hence we have 
\begin{align*}
\frac{F(t)}{F^{(n-l)}(t^{p^{n-l}})}
&=\frac{F(t)}{F^{(1)}(t^p)}
\frac{F^{(1)}(t^p)}{F^{(2)}(t^{p^2})}\cdots
\frac{F^{(n-l-1)}(t^{p^{n-l-1}})}{F^{(n-l)}(t^{p^{n-l}})}\\
&\equiv\frac{[F(t)]_{<p^d}}{[F^{(1)}(t^p)]_{<p^d}}\frac{[F(t^p)]_{<p^d}}{[F^{(1)}(t^{p^2})]_{<p^d}}\cdots
\frac{[F^{(n-l-1)}(t^{p^{n-l-1}})]_{<p^d}}{[F^{(n-l)}(t^{p^{n-l}})]_{<p^d}}
\mod p^{d-n+l+1} \Z_p[[t]]\\
&=\frac{[F(t)]_{<p^d}}{[F^{(n-l)}(t^{p^{n-l}})]_{<p^d}},
\end{align*}
namely there are $a_i\in \Z_p$ such that
\[
\frac{F(t)}{F^{(n-l)}(t^{p^{n-l}})}=
\frac{[F(t)]_{<p^d}}{[F^{(n-l)}(t^{p^{n-l}})]_{<p^d}}
+p^{d-n+l+1}\sum_{i}a_it^i.
\]
Substitute $t$ for $\zeta^s t$ in the above and 
multiply it by
\[
\left(\frac{F(t)}{F^{(n-l)}(t^{p^{n-l}})}\right)^{-1}=
\left(\frac{[F(t)]_{<p^d}}{[F^{(n-l)}(t^{p^{n-l}})]_{<p^d}}
+p^{d-n+l+1}\sum_{i}a_it^i\right)^{-1}.
\]
Then we have
\[
F(\zeta^st)\cdot [F(t)]_{<p^d}-[F(\zeta^st)]_{<p^d}\cdot F(t)=p^{d-n+l+1}\sum_{i=0}^\infty b_i(\zeta^s)t^i
\]
where $b_i(x)\in \Z_p[x]$ are polynomials which do not depend on $s$.
Applying $\sum_{s=0}^{p^{n-l}-1}\zeta^{-sk}(-)$ on both side,
one has
\[
p^{n-l}F_k(t)[F(t)]_{<p^d}-p^{n-l}[F_k(t)]_{<p^d} F(t)=p^{d-n+l+1}\sum_{i=0}^\infty 
\sum_{s=0}^{p^{n-l}-1}\zeta^{-sk}b_i(\zeta^s)t^i
\]
by \eqref{cong-lem7-eq3}. Since $\sum_{s=0}^{p^{n-l}-1}\zeta^{sj}=0$ or $p^{n-l}$,
the right hand side is zero modulo $p^{d+1}$.
Therefore
\[
\frac{F_k(t)}{F(t)}\equiv 
\frac{[F_k(t)]_{<p^d}}{[F(t)]_{<p^d}}
\mod p^{d-n+l+1}\Z_p[[t]].
\]
This implies
\[
\frac{F_k(t)[F_j(t)]_{<p^d}-[F_k(t)]_{<p^d}F_j(t)}{F(t)}\equiv
\frac{[F_k(t)]_{<p^d}[F_j(t)]_{<p^d}-[F_k(t)]_{<p^d}[F_j(t)]_{<p^d}}{[F(t)]_{<p^d}}=0\mod p^{d-n+l+1}.
\]
Now \eqref{cong-lem7-eq2} is the case $(d,j)=(n-1,s-k)$.
\end{pf}

\subsection{Proof of Congruence relations : End of proof}\label{cong-sect4}
We finish the proof of Theorem \ref{cong-thm}.
Let $S_m$ be as in Lemma \ref{cong-lem4}.
The goal is to show
\[
S_m\equiv 0\mod p^n,\quad \forall\, m\geq 0.
\]
Let us put
\[
q_i:=A_i/\wA_{\lfloor i/p\rfloor},\quad
A(i,j):=\wA_i\wA_j,\quad A^*(i,j):=A(j,i+p^{n-1})-A(i,j+p^{n-1})
\]
\[
B(i,j):=A^*(\lfloor i/p\rfloor,\lfloor j/p\rfloor).
\]
Then 
\[
S_m\equiv\sum_{i+j=m}B(i,j)q_iq_j\frac{B_j}{A_j}\mod p^n
\]
by Lemma \ref{cong-lem6}.
It follows from Lemma \ref{cong-lem3} and Lemma \ref{cong-lem5} that we have
\begin{equation}\label{pf-eq1}
k\equiv k'\mod p^i\quad\Longrightarrow\quad \frac{B_k}{A_k}\equiv\frac{B_{k'}}{A_{k'}},
\,q_k\equiv q_{k'}
\mod p^{i}.
\end{equation}
By Lemma \ref{cong-lem7}, we have
\begin{equation}\label{pf-eq2}
\sum_{\underset{i\equiv k\text{ mod }p^{n-l}}{i+j=s}}A^*(i,j)\equiv0 \mod p^l,\quad 0\leq l\leq n
\end{equation}
for all $s\geq0$.
Let
$m=l+sp$ with $l\in\{0,1,\ldots,p-1\}$.
Note
\[
B(i,m-i)=\begin{cases}
A^*(k,s-k)&kp\leq i\leq kp+l\\
A^*(k,s-k-1)&kp+l< i\leq (k+1)p-1.
\end{cases}
\]
Therefore
\begin{align*}
S_{m}&\equiv \sum_{i+j=m}B(i,j)q_iq_j\frac{B_j}{A_j}\mod p^n\\
&=
\sum_{i=0}^{p-1}\sum_{k=0}^{\lfloor(m-i)/p\rfloor}B(i+kp,m-(i+kp))q_{i+kp}q_{m-(i+kp)}
\frac{B_{m-(i+kp)}}{A_{m-(i+kp)}}
\\
&=
\sum_{k=0}^s
B(i+kp,m-(i+kp))
\sum_{i=0}^{l}q_{i+kp}q_{m-(i+kp)}
\frac{B_{m-(i+kp)}}{A_{m-(i+kp)}}
\\
&\quad+\sum_{k=0}^{s-1}
B(i+kp,m-(i+kp))
\sum_{i=l+1}^{p-1}q_{i+kp}q_{m-(i+kp)}
\frac{B_{m-(i+kp)}}{A_{m-(i+kp)}}
\\
&=
\sum_{k=0}^s
A^*(k,s-k)
\overbrace{\left(\sum_{i=0}^{l}q_{i+kp}q_{m-(i+kp)}
\frac{B_{m-(i+kp)}}{A_{m-(i+kp)}}\right)}^{P_k}
\\
&\quad+\sum_{k=0}^{s-1}
A^*(k,s-k-1)
\underbrace{\left(\sum_{i=l+1}^{p-1}q_{i+kp}q_{m-(i+kp)}
\frac{B_{m-(i+kp)}}{A_{m-(i+kp)}}\right)}_{Q_k}.
\end{align*}
We show that the first term vanishes modulo $p^n$.
It follows from \eqref{pf-eq1} that we have
\begin{equation}\label{pf-eq3}
k\equiv k'\mod p^i\quad\Longrightarrow\quad 
P_k\equiv P_{k'}
\mod p^{i+1}.
\end{equation}
Therefore one can write
\[
\sum_{k=0}^{s}A^*(k,s-k)P_k\equiv\sum_{i=0}^{p^{n-1}-1}P_i
\overbrace{\left(\sum_{k\equiv i\text{ mod } p^{n-1}}A^*(k,s-k)
\right)}^{(\ast)}\mod p^n.
\]
It follows from \eqref{pf-eq2} that $(\ast)$ is zero modulo $p$.
Therefore, again by \eqref{pf-eq3}, one can rewrite
\[
\sum_{k=0}^{s}A^*(k,s-k)P_k\equiv\sum_{i=0}^{p^{n-2}-1}P_i
\overbrace{\left(\sum_{k\equiv i\text{ mod } p^{n-2}}A^*(k,s-k)
\right)}^{(\ast\ast)}\mod p^n.
\]
It follows from \eqref{pf-eq2} that $(\ast\ast)$ is zero modulo $p^2$, so that one has
\[
\sum_{k=0}^{s}A^*(k,s-k)P_k\equiv\sum_{i=0}^{p^{n-3}-1}P_i
\left(\sum_{k\equiv i\text{ mod } p^{n-3}}A^*(k,s-k)
\right)\mod p^n
\]
by \eqref{pf-eq3}.
Continuing the same discussion, one finally obtains
\[
\sum_{k=0}^{s}A^*(k,s-k)P_k\equiv
\sum_{k=0}^sA^*(k,s-k)=0\mod p^n
\]
the vanishing of the first term.
In the same way one can show the vanishing of the second term,
\[
\sum_{k=0}^{s}A^*(k,s-1-k)Q_k\equiv0\mod p^n.
\]
We thus have $S_m\equiv 0$ mod $p^n$.
This completes the proof of Theorem \ref{cong-thm}.
\section{Geometric aspect of $p$-adic hypergeometric functions of logarithmic type}
\label{reg-sect}
We mean by a {\it fibration} over a ring $R$ a projective flat morphism 
of quasi-projective smooth $R$-schemes.
Let $X$ be a smooth $R$-scheme. We mean by a relative {\it normal crossing divisor}
(abbreviated to NCD) in $X$ over $R$
a divisor in $X$ which is locally defined by an equation $x^{r_1}_1\cdots x_s^{r_s}$
where $r_i>0$ are integers and $(x_1,\ldots,x_n)$ is a local coordinate of $X/R$.
We say a divisor $D$ simple if $D$ is a union of $R$-smooth divisors.
\subsection{Hypergeometric Curves}\label{HG-sect}
\label{fermat-sect}
Let $A$ be a commutative ring.
Let $\P^1_A(Z_0,Z_1)$ denote the projective line over $A$ with homogeneous coordinate
$(Z_0:Z_1)$.
Let $N,M\geq 2$ be integers which are invertible in $A$.
Let $t\in A$ such that $t(1-t)\in A^\times$.
Define $X$ to be a projective scheme over $A$ defined by a bihomogeneous equation
\begin{equation}\label{HGC}
(X_0^N-X_1^N)(Y_0^M-Y_1^M)=tX^N_0Y_0^M
\end{equation}
in $\P^1_A(X_0,X_1)\times\P^1_A(Y_0,Y_1)$.
We call it
a {\it hypergeometric curve} over $A$.
The morphism $X\to\Spec A$ is smooth projective
with connected fibers of relative dimension one, and 
the genus of a geometric fiber is $(N-1)(M-1)$ (e.g. Hurwitz formula).
We put $x:=X_1/X_0$ and $y:=Y_1/Y_0$, and often refer to an affine equation
\[
(1-x^N)(1-y^M)=t.
\]
In what follows, we only consider the case $A=W[t,(t-t^2)^{-1}]$ where
$W$ is a commutative ring in which $NM$ is invertible, and $t$ is an indeterminate.

\begin{lem}\label{int.model-lem}
Then the morphism $X\to \Spec W[t,(t-t^2)^{-1}]$ extends to a projective flat morphism
\[
f:Y\lra \P^1_W=\P^1_W(T_0,T_1)
\]
of smooth projective $W$-schemes satisfying the following conditions. 
\begin{enumerate}
\item[$(1)$]
$f$ has a semistable reduction at $t=0$.
The fiber $D:=f^{-1}(t=0)$ is a relative simple NCD over $W$, and 
the multiplicity of each component is one.
\item[$(2)$]
The fiber $f^{-1}(t=1)$ is a relative simple NCD over $W$, and 
the multiplicity of each component is either of $1,iN,jM$ with $i\in \{1,\ldots,M\}$, $\leq NM$.
\item[$(3)$]
The fiber $f^{-1}(t=\infty)$ is a relative simple NCD over $W$, and 
the multiplicity of each component is $\leq \max(N,M)$.
\end{enumerate}
\end{lem}
\begin{pf}
See \cite[\S 3.1]{A2} where the explicit construction of resolution of singularities
is done in \cite[Appendix B]{A2}.
\end{pf}

\subsection{Gauss-Manin connection}
In this section we assume that $W$ is an integral domain of characteristic zero.
Let $K=\Frac W$ be the fractional field.
Let $A=W[t,(t-t^2)^{-1}]$ and $S=\Spec A$.
For a $W$-scheme $T$ and a $W$-algebra $R$, we write $T_R=T\times_W R$.
The group $\mu_N\times\mu_M=\mu_N(\ol K)\times \mu_M(\ol K)$ acts on 
$X_{\ol K}$ in the following way
\begin{equation}\label{fermat-ss}
[\zeta, \nu]\cdot(x,y,t)=(\zeta x,\nu y,t),\quad(\zeta ,\nu )\in \mu_N\times \mu_M.
\end{equation}
For a $\ol K$-module $V$ with an action of $\mu_N\times\mu_M$,
let $V(i,j)$ denote the submodule on which $(\zeta ,\nu )$ acts by multiplication
by $\zeta^i\nu^j$ for all $(\zeta ,\nu )\in\mu_N\times\mu_M$.
Then one has the eigen decomposition
\[
H^1_\dR(X_{\ol K}/ S_{\ol K})=\bigoplus_{i=1}^{N-1}\bigoplus_{j=1}^{M-1}
H^1_\dR(X_{\ol K}/S_{\ol K})(i,j),
\]
and each eigenspace $H^1_\dR(X_{\ol K}/ S_{\ol K})(i,j)$ is free of rank $2$ over 
$\O(S_{\ol K})$
(\cite[Lemma 2.2]{A}).
Put
\begin{equation}\label{fermat-form-ab}
a_i:=1-\frac{i}{N},\quad 
b_j:=1-\frac{j}{M}.
\end{equation}
Let
\begin{equation}\label{fermat-form}
\omega_{i,j}:=N\,\frac{x^{i-1}y^{j-M}}{1-x^N}dx
=-M\,\frac{x^{i-N}y^{j-1}}{1-y^M}dy,
\end{equation}
\begin{equation}\label{fermat-form-2}
\eta_{i,j}:
=\frac{1}{x^N-1+t}\omega_{i,j}=Mt^{-1}x^{i-N}y^{j-M-1}dy
\end{equation}
for integers $i,j$ such that $1\leq i\leq N-1,\,1\leq j\leq M-1$.
Then $\omega_{i,j}$ is the 1st kind, and $\eta_{i,j}$ is the 2nd kind.
They form a $\O(S_{\ol K})$-free basis of $H^1_\dR(X_{\ol K}/S_{\ol K})(i,j)$.
According to this, we put
\[H_\dR^1(X_K/S_K)(i,j):=\O(S_K)\omega_{i,j}+\O(S_K)\eta_{i,j}\subset H^1_\dR(X_K/S_K).
\]
Let 
\[
F_{a,b}(t):={}_2F_1\left({a,b\atop 1};t\right)=\sum_{i=0}^\infty\frac{(a)_i}{i!}\frac{(b)_i}{i!}
t^i\in K[[t]]
\]
be the hypergeometric series.
Put
\begin{equation}\label{fermat-form-wt}
\wt\omega_{i,j}:=\frac{1}{F_{a_i,b_j}(t)}\omega_{i,j},\quad
\wt\eta_{i,j}:=-t(1-t)^{a_i+b_j}(F'_{a_i,b_j}(t)\omega_{i,j}+b_jF_{a_i,b_j}(t)\eta_{i,j}).
\end{equation}
For the later use, we give notation $V(i,j),\omega_{i,j},\eta_{i,j},\wt\omega_{i,j},\wt\eta_{i,j}$
for $(i,j)$ not necessarily a pair of integers.
Let $(i,j)=(q/r,q'/r')\in \Q^2$ such that $\gcd(r,N)=\gcd(r',M)=1$ and
$N\not| q$ and $M\not| q'$.
Let $i_0$, $j_0$ be the unique integers such that $i_0\equiv i$ mod $N$,
$j_0\equiv j$ mod $M$ and $1\leq i_0<N$, $1\leq j_0<M$. Then we define
\begin{equation}\label{rule}
V(i,j)=V(i_0,j_0),\quad \omega_{i,j}=\omega_{i_0,j_0},\quad\ldots\quad
\wt\eta_{i,j}=\wt\eta_{i_0,j_0}.
\end{equation}
\begin{prop}\label{fermat-GM}
Let $\nabla:H^1_\dR(X_K/S_K)\to\O(S_K)dt\ot H^1_\dR(X_K/S_K)$
be the Gauss-Manin connection. It naturally extends on $K((t))\ot_{\O(S)}H^1_\dR(X_K/S_K)$
which we also write by $\nabla$. Then
\[
\begin{pmatrix}
\nabla(\omega_{i,j})&\nabla(\eta_{i,j})
\end{pmatrix}
=dt\ot\begin{pmatrix}
\omega_{i,j}&\eta_{i,j}
\end{pmatrix}
\begin{pmatrix}
0&-a_i(t-t^2)^{-1}\\
-b_j&(-1+(1+a_i+b_j)t)(t-t^2)^{-1}
\end{pmatrix},
\]
\[
\begin{pmatrix}
\nabla(\wt\omega_{i,j})&\nabla(\wt\eta_{i,j})
\end{pmatrix}
=dt\ot\begin{pmatrix}
\wt\omega_{i,j}&\wt\eta_{i,j}
\end{pmatrix}
\begin{pmatrix}
0&0\\
t^{-1}(1-t)^{-a_i-b_j}F_{a_i,b_j}(t)^{-2}&0
\end{pmatrix}.
\]
\end{prop}
\begin{pf}
We may replace the base field with $\C$.
Since $\nabla$ is commutative with the action of $\mu_N(\C)\times\mu_M(\C)$, 
$\nabla$ preserves the eigen components $H^1_\dR(X/S)(i,j)$.
We think $X$ and $S$ of being complex manifolds.
For $\alpha\in \C\setminus\{0,1\}$ we write $X_\alpha=f^{-1}(t=\alpha)$.
Then there is a homology cycle $\delta_\alpha\in H_1(X_\alpha,\Q)$
such that
\[
\int_{\delta_\alpha}\omega_{i,j}=2\pi\sqrt{-1}\,{}_2F_1\left({a_i,b_j\atop 1};\alpha\right)
\]
(\cite[Lemma 2.3]{A}).
Let $\partial_t=\nabla_{\frac{d}{dt}}$ be the differential operator on 
$\O(S)^{an}\ot_{\O(S)} H^1_\dR(X/S)$.
Put $D=t\partial_t$ and $P_{\mathrm{HG}}=D^2-t(D+a_i)(D+b_j)$ the hypergeometric
differential operator.
Since the series ${}_2F_1\left({a_i,b_j\atop 1};t\right)$ is annihilated by $P_{\mathrm{HG}}$,
we have
\[
\int_{\delta_\alpha}P_{\mathrm{HG}}(\omega_{i,j})=0.
\]
Since $H_1(X_\alpha,\C)(i,j)$ is a 2-dimensional irreducible $\pi_1(S,\alpha)$-module,
we have $\int_{\gamma}P_{\mathrm{HG}}(\omega_{i,j})=0$ for all
$\gamma\in\pi_1(S,\alpha)$, which means
\begin{equation}\label{fermat-GM-eq1}
P_{\mathrm{HG}}(\omega_{i,j})=(D^2-t(D+a_i)(D+b_j))(\omega_{i,j})=0.
\end{equation}
Next we show
\begin{equation}\label{fermat-GM-eq2}
\partial_t(\omega_{i,j})=-b_j\eta_{i,j}.
\end{equation}
Put
$\phi:=N\,\frac{x^{i-1}}{1-x^N}dx\in \vg(U,\Omega^1_{X/\C})$
with $U=\{x\ne\infty,y\ne\infty\}\subset X$, and 
$\Omega_{i,j}:=y^{j-M}\phi$ a lifting of $\omega_{i,j}$.
Since $\phi$ is a linear combination of $dx/(x-\nu)$'s,
one has $d(\phi)=0$.
Therefore
\[
d(\Omega_{i,j})=
d(y^{j-M})\wedge\phi=(j-M)y^{j-M-1}dy\wedge\phi\in \vg(U,\Omega^2_{X/\C}).
\]
Taking $\wedge\phi$ on both sides of 
\[
\frac{dt}{t}=\frac{Nx^{N-1}dx}{x^N-1}+\frac{My^{M-1}dy}{y^M-1},
\]
we have
\[
\frac{dt}{t}\wedge\phi=\frac{My^{M-1}dy}{y^M-1}\wedge\phi
\iff
\frac{y^M-1}{My^{M-1}}\frac{dt}{t}\wedge\phi=dy\wedge\phi,
\]
and hence
\[
d(\Omega_{i,j})
=-(1-j/M)y^{j-M}
\frac{y^M-1}{y^{M}}\frac{dt}{t}\wedge\phi
=-(1-j/M)\frac{y^M-1}{ty^{M}}dt\wedge\Omega_{i,j}
=-b_jdt\wedge\frac{1}{x^N-1+t}\Omega_{i,j}.
\]
Since $(x^N-1+t)^{-1}\Omega_{i,j}$ is a lifting of $\eta_{i,j}$, 
this shows $\nabla(\omega_{i,j})=-b_jdt\ot\eta_{i,j}$. This completes the proof of 
\eqref{fermat-GM-eq2}.
Now all the formulas on $\nabla$ follow from \eqref{fermat-GM-eq1} and
\eqref{fermat-GM-eq2}.
\end{pf}
The following is straightforward from Proposition \ref{fermat-GM}.
\begin{cor}\label{fermat-GM-cor1}
Let $\nabla_{i,j}$ be the connection on the eigen component
$H_{i,j}:=K((t))\ot_{\O(S)}H^1_\dR(X_K/S_K)(i,j)$. Then
$\ker\nabla_{i,j}=K\wt{\eta}_{i,j}$.
Moreover let $\ol\nabla_{i,j}$ be the connection on $M_{i,j}:=H_{i,j}/K((t))\wt\eta_{i,j}$
induced from $\nabla_{i,j}$.
Then $\ker\ol\nabla_{i,j}=K\wt{\omega}_{i,j}$.
\end{cor}
We mean by
a semistable family
$g:\cX\to \Spec R[[t]]$ over a commutative ring $R$ 
that $g$ is a proper flat morphism, smooth
over $\Spec R((t))$ and it is locally described by
\[
g:\Spec R[[x_1,\ldots,x_n]]\lra \Spec R[[t],\quad g^*(t)=x_1\cdots x_r
\]
in each formal neighborhood.
Let $D$ be the fiber at $\Spec R[[t]]/(t)$, which is a relative NCD in $\cX$ over $R$
with no multiplicities.
Let ${\mathscr U}:=\cX\setminus D$.
We define the log de Rham complex $\omega^\bullet_{\cX/R[[t]]}$
to be the subcomplex of $\Omega^\bullet_{{\mathscr U}/R((t))}$ generated by $dx_i/x_i$
($1\leq i\leq r$) and $dx_j$ ($j>r$) over $\O_\cX$.
Equivalently,
\begin{equation}\label{log-dR-cpx}
\omega^\bullet_{\cX/R[[t]]}:=\Coker\left[
\frac{dt}{t}\ot\Omega^{\bullet-1}_{\cX/R}(\log D)
\to \Omega^{\bullet}_{\cX/R}(\log D)\right]
\end{equation}
where $\Omega^{\bullet}_{\cX/R}(\log D)$ denotes the $t$-adic completion of
the complex of the 
algebraic K\"ahler differentials, 
\[
\Omega^{\bullet}_{\cX/R}(\log D):=
\varprojlim_{n\geq1}\big(\Omega^{\bullet,\text{alg}}_{\cX/R}(\log D)/t^n\Omega^{\bullet,\text{alg}}_{\cX/R}(\log D)\big)
\]
\begin{cor}\label{fermat-GM-cor2}
Let $f:Y\to\P^1_W$ be the morphism of projective smooth $W$-schemes in Lemma \ref{int.model-lem}.
Let $Y_K:=Y\times_WK$.
Let $\Delta_K:=\Spec K[[t]]\hra \P^1_K$ be the formal neighborhood, and put
$\cY_K:=f^{-1}(\Delta_K)$. Let $D_K\subset\cY_K$ be the fiber at $t=0$.
Let
\[
\xymatrix{
D_K\ar[r]\ar[d]&\cY_K\ar[d]\ar[r]&Y_K\ar[d]\\
\{0\}\ar[r]&\Delta_K\ar[r]&\P^1_K
}\]
Put $H_K:=H^1_\zar(\cY_K,\omega^\bullet_{\cY_K/K[[t]]})$.
It follows from \cite[(2.18)--(2.20)]{steenbrink} or
\cite[(17)]{zucker} that $H_K\to K((t))\ot_{\O(S)} H^1_\dR(X/S)$ is injective.
We identify $H_K$ with its image.
Then the eigen component
$H_K(i,j)$ is a free $K[[t]]$-module with basis
$\{\wt\omega_{i,j},\wt\eta_{i,j}\}$.
\end{cor}
\begin{pf}
$H_K$ is called {\it Deligne's canonical extension}, and 
is characterized by the following conditions (\cite[(17)]{zucker}).
\begin{description}
\item[(D1)] $H_K$ is a free $K[[\l]]$-module such that $K((t))\ot H_K=K((t))\ot_{\O(S)} H^1_\dR(X/S)$,
\item[(D2)]
the connection extends to have log pole, 
$\nabla:H_K\to\frac{dt}{t}\ot H_K$,
\item[(D3)]
each eigenvalue $\alpha$ of $\Res(\nabla)$ satisfies $0\leq \mathrm{Re}(\alpha)<1$, 
where $\Res(\nabla)$ is the $K$-linear endomorphism defined by a commutative diagram
\[
\xymatrix{
H_K\ar[r]^\nabla\ar[d]&
\frac{dt}{t}\ot H_K\ar[d]^{\Res\ot1}\\
H_K/tH_K\ar[r]^{\Res(\nabla)}&H_K/tH_K.
}
\]
\end{description}
Put $H^0_K:=\bigoplus_{i,j}K[[t]]\wt\omega_{i,j}+K[[t]]\wt\eta_{i,j}$.
We can directly check
that $H^0_K$ satisfies {\bf(D1)}--{\bf(D3)} by Propositioon \ref{fermat-GM}.
We then conclude $H_K=H^0_K$ thanks to the uniqueness of Deligne's canonical extension.
\end{pf}
\subsection{Rigid cohomology and a category $\FilFMIC(S)$}
In what follows, let the base ring $W$ be the Witt ring $W(\ol\F_p)$ of the algebraic closure
$\ol\F_p$ with $p\nmid NM$.
Let $F$ be the $p$-th Frobenius on $W$, and $K:=\Frac W$ the fractional field.

\medskip

Let $f:X\to S$ be the hypergeometric curve as before.
Write $X_{\ol\F_p}:=X\times_{W}\ol\F_p$ and $S_{\ol\F_p}:=S\times_{W}\ol\F_p$.
Let $\sigma$ be a $F$-linear $p$-th Frobenius on $W[t,(t-t^2)^{-1}]^\dag$ the ring of
overconvergent power series, which naturally extends on $K[t,(t-t^2)^{-1}]^\dag:=
K\ot W[t,(t-t^2)^{-1}]^\dag$.
Then the $i$-th {\it rigid cohomology} group
\[
H^i_\rig(X_{\ol\F_p}/S_{\ol\F_p}) 
:=\Gamma(S_K^{\mathrm{an}}, R^if_{\rig}j^{\dag}_X\O_{X_K^{\mathrm{an}}}),
\]
is defined where $R^i f_{\rig}j^{\dag}_X\O_{X_K^{\mathrm{an}}}$
is the $i$-th relative rigid cohomology sheaf (cf. \cite[Definition 2.12]{AM}
for the notation and remark on the definition).
The required properties in below is the following (loc.cit).
\begin{itemize}
\item
$H^\bullet_\rig(X_{\ol\F_p}/S_{\ol\F_p})$ is a finitely generated 
$\O(S)^\dag=K[t,(t-t^2)^{-1}]^\dag$-module.
\item(Frobenius)
The $p$-th Frobenius $\Phi_{X/S}$ on $H^\bullet_\rig(X_{\ol\F_p}/S_{\ol\F_p})$ (depending on $\sigma$)
is defined in a natural way. This is a $\sigma$-linear endomorphism :
\[
\Phi_{X/S}(h(t)x)=\sigma(h(t))\Phi_{X/S}(x),\quad \mbox{for }x\in 
H^\bullet_\rig(X_{\ol\F_p}/S_{\ol\F_p}),\, 
h(t)\in \O(S)^\dag.
\]
\item(Comparison)
There is the comparison isomorphism with the algebraic de Rham cohomology,
\[
c:H^\bullet_\rig(X_{\ol\F_p}/S_{\ol\F_p})\cong H^\bullet_\dR(X/S)
\ot_{\O(S)} \O(S)^\dag.
\]
\end{itemize}
In \cite[\S 2.1]{AM} we introduce
a category $\FilFMIC(S)=\FilFMIC(S,\sigma)$.
It consists of collections of datum $(H_{\dR}, H_{\rig}, c, \Phi, \nabla, \Fil^{\bullet})$ such that
    \begin{itemize}
            \setlength{\itemsep}{0pt}
        \item $H_{\dR}$ is a finitely generated $\O(S)$-module,
        \item $H_{\rig}$ is a finitely generated $\O(S)^\dag$-module,
        \item $c\colon H_\rig\cong H_{\dR}\otimes_{\O(S)}
        \O(S)^\dag$, the comparison 
        \item $\Phi\colon \sigma^{\ast}H_{\rig}\xrightarrow{\,\,\cong\,\,} H_{\rig}$ is an isomorphism of $\O(S)^\dag$-module,
        \item $\nabla\colon H_{\dR}\to \Omega_{S/\Q_p}^1\otimes H_{\dR}$ is an integrable connection
    that satisfies $\Phi\nabla=\nabla\Phi$.
        \item $\Fil^{\bullet}$ is a finite descending filtration on $H_{\dR}$ of locally free 
        $\O(S)$-module (i.e. each graded piece is locally free),
    that satisfies $\nabla(\Fil^i)\subset \Omega^1\ot \Fil^{i-1}$.
    \end{itemize}

Let $\mathrm{Fil}^\bullet$ denote the Hodge filtration on the de Rham cohomology,
and $\nabla$ the Gauss-Manin connection.
Let
\[
H^i(X/S):=(H^i_\dR(X/S),H^i_\rig(X_{\ol\F_p}/S_{\ol\F_p}),c,\Phi_{X/S},\nabla,\mathrm{Fil}^\bullet)
\]
be an object of $\FilFMIC(S)$. 

For an integer $r$, the Tate object $\O_S(r)\in \FilFMIC(S)$ 
is defined in a customary way (loc.cit.).
We simply write
\[
M(r)=M\ot \O_S(r)
\]
for an object $M\in \FilFMIC(S)$.

\begin{lem}\label{frobenius-lem}
Suppose that $\sigma$ is given by $\sigma(t)=ct^p$ with some $c\in 1+pW$.
Then, with the notation in Corollary \ref{fermat-GM-cor2},
the Frobenius $\Phi_{X/S}$ induces the action on 
$H_K$ in a natural way.
\end{lem}
\begin{pf}
Let $W((t))^\wedge$ be the $p$-adic completion and write $K((t))^\wedge:=K\ot_WW((t))^\wedge$ on which $\sigma$ extends as $\sigma(t)=ct^p$.
The Frobenius $\Phi_{X/S}$ on $H^1_\dR(X/S)\ot_{\O(S)}\O(S)^\dag$ naturally extends
on $H^1_\dR(X/S)\ot_{\O(S)}K((t))^\wedge$ via the homomorphism
$\O(S)^\dag\to K((t))^\wedge$. We show that the action of $\Phi_{X/S}$ preserves
the subspace $H_K$

Let $f:Y\to\P^1_W$ be the morphism of projective smooth $W$-schemes in Lemma \ref{int.model-lem}.
Let
\[
\xymatrix{
D_W\ar[r]\ar[d]&\cY\ar[d]\ar[r]&Y\ar[d]\\
0\ar[r]&\Delta_W\ar[r]&\P^1_W
}\]
where $\Delta_W:=\Spec W[[t]]\hra \P^1_W$ is the formal neighborhood, and 
$0=\Spec W[[t]]/(t)$.
Note that $D_W$ is reduced, namely $\cY/\Delta_W$ has a semistable reduction.
Write $\cY_{\ol\F_p}:=\cY\times_W\ol\F_p$ etc.
We employ the log-crystalline cohomology
\begin{equation}\label{log-crys}
H^\bullet_{\text{log-crys}}((\cY_{\ol\F_p},D_{\ol\F_p})/(\Delta_W,0))
\end{equation}
where $(\cX,\cD)$ denotes the log scheme with log structure
induced by the divisor $\cD$.
There is the comparison theorem by Kato \cite[Theorem 6.4]{Ka-log},
\begin{equation}\label{log-crys-isom}
H^\bullet_{\text{log-crys}}((\cY_{\ol\F_p},D_{\ol\F_p})/(\Delta_W,0))
\cong 
H^\bullet(\cY,\omega^\bullet_{\cY/W[[t]]})
\end{equation}
(see \eqref{log-dR-cpx} for the complex $\omega^\bullet_{\cY/W[[t]]}$).
The log-crystalline cohomology is endowed with
the $p$-th Frobenius $\Phi_{(\cY,D_W)}$
which is compatible with $\Phi_{X/S}$ under the map
\begin{align*}
H^1(\cY,\omega^\bullet_{\cY/W[[t]]})&\to
H^1(\cY_K,\omega^\bullet_{\cY_K/K[[t]]}))\\
&\hra
H^1_\dR(X/S)\ot_{\O(S)} K((t))\\
&\hra H^1_\dR(X/S)\ot_{\O(S)} K((t))^\wedge
\end{align*}
where $\cY_K:=\cY\times_{W[[t]]}K[[t]]$ and $D_K:=D_W\times_WK$.
Thus the assertion follows.
\end{pf}
\begin{prop}\label{frobenius-thm}
Let $\wt\omega_{i,j},\wt\eta_{i,j}$ be as in \eqref{fermat-form} and \eqref{fermat-form-2}.
Suppose that $\sigma$ is given by $\sigma(t)=ct^p$ with some $c\in 1+pW$.
Then
\[
\Phi_{X/S}(\wt\eta_{p^{-1}i,p^{-1}j})\in K\wt\eta_{i,j},\quad
\Phi_{X/S}(\wt\omega_{p^{-1}i,p^{-1}j})\equiv p\wt\omega_{i,j}\mod K((t))\wt\eta_{i,j}
\]
where we use the notation \eqref{rule}.
\end{prop}
\begin{pf}
Let $\nabla$ be the Gauss-Manin connection on $H^1_\dR(X/S)\ot_{\O(S)}K((t))$.
Since $\Phi_{X/S}\nabla=\nabla\Phi_{X/S}$, we have 
$\Phi_{X/S}\ker(\nabla)\subset \ker(\nabla)$.
Moreover,
$\Phi_{X/S}$ sends the eigencomponents $H_{i,j}:=H^1_\dR(X/S)(i,j)\ot_{\O(S)}K((t))$
onto the component $H_{pi,pj}$ as 
$\Phi_{X/S}[\zeta,\nu]=[\zeta,\nu]\Phi_{X/S}$.
Therefore we have
\[
\Phi_{X/S}(\wt\eta_{p^{-1}i,p^{-1}j})\in K\wt\eta_{i,j}
\]
by Corollary \ref{fermat-GM-cor1}.
We show the latter.
By Lemma \ref{frobenius-lem} together with Corollary \ref{fermat-GM-cor2}, 
there are $f_{i,j}(t), g_{i,j}(t)\in K[[t]]$ such that
$\Phi_{X/S}(\wt\omega_{p^{-1}i,p^{-1}j})=f_{i,j}(t)\wt\omega_{i,j}+g_{i,j}(t)\wt\eta_{i,j}$.
Put $M_{i,j}:=H_{i,j}/K((t))\wt\eta_{i,j}$.
Then $\Phi_{X/S}(M_{p^{-1}i,p^{-1}j})\subset M_{i,j}$ and $\Phi_{X/S}$ is
commutative with the connection $\ol\nabla$ on $M_{i,j}$.
Therefore $f_{i,j}(t)=C_{i,j}$ is a constant as $\ker(\ol\nabla)=K\wt\omega_{i,j}$
by Corollary \ref{fermat-GM-cor1},
\begin{equation}\label{frob-lem-eq1}
\Phi_{X/S}(\wt\omega_{p^{-1}i,p^{-1}j})=C_{i,j}\wt\omega_{i,j}+g_{i,j}(t)\wt\eta_{i,j}.
\end{equation}
We want to show $C_{i,j}=p$.
To do this, we recall the log-crystalline cohomology \eqref{log-crys}
\[
H^\bullet_{\text{log-crys}}((\cY_{\ol\F_p},D_{\ol\F_p})/(\Delta_W,0))
\cong H^\bullet(\cY,\omega^\bullet_{\cY/W[[t]]}).
\]
where we keep the notation in the proof of Lemma \ref{frobenius-lem}.
Let $Z_W$ be the intersection locus of $D_W$.
This is a disjoint union of  $NM$-copies of $\Spec W$.
More precisely, let $P_{\zeta,\nu}$ be the point of $Z_W$ defined by $x=\zeta$ and
$y=\nu$.
Then $Z_W=\{P_{\zeta,\nu}\mid(\zeta,\nu)\in\mu_N\times\mu_M\}$.
We consider the composition of morphisms 
\[
\omega^\bullet_{\cY/W[[t]]}\os{\wedge\frac{dt}{t}}{\lra}
\Omega^{\bullet+1}_{\cY/W}(\log D)\os{\Res}{\lra}
\O_{Z_W}[-1]
\]
of complexes where $\Res$ is the Poincare residue.
This gives rise to the map
\[
R:H^1(\cY_K,\Omega^\bullet_{\cY_K/K[[t]]}(\log D_K))\lra 
H^0(Z_K,\O_{Z_K})=\bigoplus_{\zeta,\nu}K\cdot P_{\zeta,\nu}
\]
which is compatible with respect to the Frobenius $\Phi_{X/S}$ 
on the left
and the Frobenius $\Phi_Z$ on the right in the sense that
\begin{equation}\label{frobenius-thm-eq0}
R\circ \Phi_{X/S}=p\Phi_Z\circ R.
\end{equation}
Notice that $\Phi_Z$ is a $F$-linear map such that
$\Phi_Z(P_{\zeta,\nu})=P_{\zeta,\nu}$ where $F$ is the Frobenius
on $W$.
We claim
\begin{equation}\label{frobenius-thm-eq1}
R(\wt\eta_{i,j})=0,
\end{equation}
and
\begin{equation}\label{frobenius-thm-eq2}
R(\wt\omega_{i,j})=\sum_{\zeta,\nu}\zeta^i\nu^jP_{\zeta,\nu}.
\end{equation}
To show \eqref{frobenius-thm-eq1}, we recall the definition \eqref{fermat-form-wt}.
Since $R(t\wt\omega_{i,j})=0$ obviously, 
it is enough to show $R(t\eta_{i,j})=0$.
However since $t\eta_{i,j}=Mx^{i-N}y^{j-M}dy$ by
\eqref{fermat-form-2}, we have
\[
R(t\eta_{i,j})=\Res\left(Mx^{i-N}y^{j-M}dy\frac{dt}{t}\right)=0
\]
as required. One can show \eqref{frobenius-thm-eq2} as follows.
\begin{align*}
R(\wt\omega_{i,j})=R(\omega_{i,j})
&=\Res\left(
M\,\frac{x^{i-N}y^{j-1}}{y^M-1}dy\wedge\frac{dt}{t}\right)\\
&=\Res\left(
\frac{Mx^{i-N}y^{j-1}}{y^M-1}dy\wedge\frac{Nx^{N-1}}{x^N-1}dx\right)\\
&=\sum_{\zeta,\nu}\zeta^i\nu^jP_{\zeta,\nu}.
\end{align*}
We turn to the proof of $C_{i,j}=p$ in \eqref{frob-lem-eq1}.
Apply $R$ on the both side of \eqref{frob-lem-eq1}.
By \eqref{frobenius-thm-eq1}, the right hand side is $\alpha_{i,j}R(\wt\omega_{i,j})$,
and 
the left hand side is $p\Phi_Z\circ R(\wt\omega_{p^{-1}i,p^{-1}j})$ 
by \eqref{frobenius-thm-eq0},
\[
C_{i,j}R(\wt\omega_{i,j})=p\Phi_Z\circ R(\wt\omega_{p^{-1}i,p^{-1}j}).
\]
Apply \eqref{frobenius-thm-eq1} to the above. We have
\[
C_{i,j}\left(\sum_{\zeta,\nu}\zeta^i\nu^jP_{\zeta,\nu}\right)=
p\Phi_Z\left(\sum_{\zeta,\nu}\zeta^{p^{-1}i}\nu^{p^{-1}j}\cdot P_{\zeta,\nu}\right)
=p\left(\sum_{\zeta,\nu}\zeta^{i}\nu^{j}P_{\zeta,\nu}\right)
\]
and hence $C_{i,j}=p$ as required.
\end{pf}

\begin{thm}[Unit root formula]\label{uroot-thm}
Suppose $\sigma(t)=t^p$.
Let $(i,j)$ be a pair of integers $(i,j)$ with $0< i<N$ and $0< j<M$.
Put
\[
e^{\mathrm{unit}}_{i,j}:=(1-t)^{-a_i-b_j}F_{a_i,b_j}(t)^{-1}\wt\eta_{i,j}.
\]
Let $s\geq 0$ be the minimal
integer such that $a_i^{(s+1)}=a_i$ and $b_j^{(s+1)}=b_j$
\footnote{For any $a\in \Z_{(p)}$ with $0<a<1$, 
there exists $i>0$ such that $a^{(i)}=a$.
Let $n>0$ be an integer prime to $p$.
 For $a=l/n$ with $0<l<n$, the Dwork prime $a^{(1)}=k/n$ is characterized by
 $pk\equiv l$ mod $n$ and $0<k<n$.
 Therefore, the map $a\mapsto a^{(i)}$ induces a bijection on the set
 $\{1/n,2/n,\ldots,(n-1)/n\}$, and it is identity if and only if $p^i\equiv1$ mod $n$.
}.
Put $h(t):=\prod_{m=0}^s [F_{a_i^{(m)},b_j^{(m)}}(t)]_{<p}$. Then
\begin{equation}\label{uroot-thm-1}
e^{\mathrm{unit}}_{i,j}\in 
H^1_\dR(X/S)\ot_{\O(S)}K\langle t,(t-t^2)^{-1},h(t)^{-1}\rangle
\end{equation}
and
\begin{equation}\label{uroot-thm-2}
\Phi_{X/S}(e^{\mathrm{unit}}_{p^{-1}i,p^{-1}j})=
\frac{(1-t)^{a_i+b_j}}{(1-t^p)^{a^{(1)}_i+b^{(1)}_j}}
\cF^{\mathrm{Dw}}_{a_i,b_j}(t)e^{\mathrm{unit}}_{i,j}
\end{equation}
where $\cF^{\mathrm{Dw}}_{a_i,b_j}(t)$ is the Dwork $p$-adic hypergeometric function
\eqref{Dwork},
and we apply the convention \eqref{rule} to the notation
$e^{\mathrm{unit}}_{p^{-1}i,p^{-1}j}$.
In particular
$e^{\mathrm{unit}}_{p^{-1}i,p^{-1}j}$ is the eigen vector of 
$\Phi_{X/S}^{s+1}=\overbrace{\Phi_{X/S}\circ\cdots\circ\Phi_{X/S}}^{s+1}$, and
\begin{equation}\label{uroot-thm-3}
\Phi_{X/S}^{s+1}(e^{\mathrm{unit}}_{i,j})
=\left(\prod_{m=0}^s
\frac{(1-t^{p^m})^{a^{(m)}_i+b^{(m)}_j}}{(1-t^{p^{m+1}})^{a^{(m+1)}_i+b^{(m+1)}_j}}
\cF^{\mathrm{Dw}}_{a^{(m)}_i,b^{(m)}_j}(t^{p^m})\right)e^{\mathrm{unit}}_{i,j}.
\end{equation}
\end{thm}
Notice that $(1-t)^{a_i+b_j}\not\in \Z_p\langle t,(1-t)^{-1}\rangle$
but $(1-t)^{a_i+b_j}/(1-t^p)^{a_i^{(1)}+b_j^{(1)}}\in \Z_p\langle t,(1-t)^{-1}\rangle$.
\begin{pf}
Since
\begin{equation}\label{uroot-thm-eq1}
\frac{F'_{a_i,b_j}(t)}{F_{a_i,b_j}(t)}\in \Z_p\langle t,(t-t^2)^{-1},h(t)^{-1}\rangle
\end{equation}
by \cite[p.45, Lem. 3.4 ]{Dwork-p-cycle}, \eqref{uroot-thm-1} follows.
We show \eqref{uroot-thm-2}, which is equivalent to
\begin{equation}\label{uroot-thm-eq2}
\Phi_{X/S}(\wt\eta_{p^{-1}i,p^{-1}j})=\wt\eta_{i,j}\in H^1_\dR(X/S)\ot_{\O(S)}K((t)).
\end{equation}
Let
\[
Q:H^1_\dR(X_K/S_K)\ot H^1_\dR(X_K/S_K)\lra H^2_\dR(X_K/S_K)\cong \O(S_K)
\]
be the cup-product pairing which is anti-symmetric and non-degenerate.
This extends on $H^1_\dR(X/S)\ot_{\O(S)}K((t))$ which we also write by $Q$.
Then the following is satisfied.
\begin{description}
\item[(Q1)]
$Q(\Phi_{X/S}(x),\Phi_{X/S}(y))=pQ(x,y)^\sigma$ for $x,y\in H^1_\dR(X/S)\ot_{\O(S)}K((t))$,
\item[(Q2)]
$Q(gx,gy)=Q(x,y)$ for $g=(\zeta,\nu)\in \mu_N\times \mu_M$,
\item[(Q3)]
$Q(F^1,F^1)=0$ where $F^1=\vg(X_K,\Omega^1_{X_K/S_K})$ is the Hodge filtration, 
\item[(Q4)]
$Q(\nabla(x),y)+Q(x,\nabla(y))=dQ(x,y)$.
\end{description}
Put $H_{i,j}:=H^1_\dR(X/S)(i,j)\ot_{\O(S)}K((t))$ eigen components.
By {\bf(Q2)}, $Q$ induces a perfect pairing $H_{i,j}\ot H_{N-i,M-j}\to\O(S)$.
Therefore $Q(\wt\omega_{i,j},\wt\eta_{N-i,M-j})\ne0$ by {\bf(Q3)}.
We claim 
\begin{equation}\label{uroot-thm-eq3}
Q(\wt\eta_{i,j},\wt\eta_{N-i,M-j})=0,
\end{equation}
\begin{equation}\label{uroot-thm-eq4}
Q(\wt\omega_{i,j},\wt\eta_{N-i,M-j})\in \Q^\times.
\end{equation}
To show \eqref{uroot-thm-eq3}, we recall $H_K$ in Corollary \ref{fermat-GM-cor2}.
Since $Q$ is the cup-product pairing, this induces a pairing $H_K\ot H_K\to K[[t]]$, and hence 
\[
\ol Q:H_K/tH_K\ot_K H_K/tH_K\lra K.
\]
Since $\nabla(\wt\eta_{i,j})=0$, $Q(\wt\eta_{i,j},\wt\eta_{N-i,M-j})$ is a constant by {\bf(Q4)}.
Therefore if one can show $\ol Q(\wt\eta_{i,j},\wt\eta_{N-i,M-j})=0$, 
then \eqref{uroot-thm-eq3} follows. It follows from {\bf(Q4)} that
\[
\ol Q(\Res(\nabla)(x),y)+\ol Q(x,\Res(\nabla)(y))=0,\quad \forall\,x,y\in H_K/tH_K
\]
where $\Res(\nabla)$ is as in the proof of Corollary \ref{fermat-GM-cor2}.
Since $\Res(\nabla)(\wt\omega_{i,j})=\wt\eta_{i,j}$ 
and $\Res(\nabla)\wt\eta_{i,j}=0$ by Proposition \ref{fermat-GM},
one has
\[
\ol Q(\wt\eta_{i,j},\wt\eta_{N-i,M-j})
=\ol Q(\Res(\nabla)\wt\omega_{i,j},\wt\eta_{N-i,M-j})
=-\ol Q(\wt\omega_{i,j},\Res(\nabla)\wt\eta_{N-i,M-j})=0
\]
as required. We show \eqref{uroot-thm-eq4}.
Since $\nabla(\wt\omega_{i,j})
\in K((t))\wt\eta_{i,j}$, we have
\[
dQ(\wt\omega_{i,j},\wt\eta_{N-i,M-j})=
Q(\nabla(\wt\omega_{i,j}),\wt\eta_{N-i,M-j})=0
\]
by \eqref{uroot-thm-eq3} which means that $Q(\wt\omega_{i,j},\wt\eta_{N-i,M-j})$ is a constant.
Since $X/S$, $Q$ and $\wt\omega_{i,j},\wt\eta_{i,j}$ are defined over $\Q((t))$, 
the constant 
should belong to $\Q^\times$.
This completes the proof of \eqref{uroot-thm-eq4}.

\medskip
We turn to the proof of \eqref{uroot-thm-eq2}.
By Proposition \ref{frobenius-thm}, there is a constant $\alpha\in K$ such that
$\Phi_{X/S}(\wt\eta_{p^{-1}i,p^{-1}j})=\alpha\wt\eta_{i,j}$.
Put $c:=Q(\wt\omega_{i,j},\wt\eta_{N-i,M-j})$ which belongs to $\Q^\times$
by \eqref{uroot-thm-eq4}.
By {\bf(Q1)}, we have
\[Q(\Phi_{X/S}(\wt\omega_{i,j}),\Phi_{X/S}(\wt\eta_{N-i,M-j}))=pQ(\wt\omega_{i,j},\wt\eta_{N-i,M-j})^\sigma=pc,\]
and hence \[\alpha Q(\Phi_{X/S}(\wt\omega_{i,j}),\wt\eta_{N-i,M-j})=pc.\]
It follows from \eqref{uroot-thm-eq3} and Proposition \ref{frobenius-thm} 
that the left hand side is 
\[
p\alpha Q(\wt\omega_{i,j},\wt\eta_{N-i,M-j})=p\alpha c.
\]
Therefore $\alpha=1$. This completes the proof.
\end{pf}

\subsection{Syntomic Regulators of hypergeometric curves} \label{syn-reg-sect}
Let $f_R:Y_R\to \P^1_R$ be the fibration of hypergeometric curves 
in Lemma \ref{int.model-lem}
over a ring $R:=\Z[\zeta_N,\zeta_M,(NM)^{-1}]$ 
where $\zeta_n$ is a primitve $n$-th root of unity in $\ol\Q$.
Let $S_R:=\Spec R[t,(t-t^2)^{-1}]$ and $X_R:=f_R^{-1}(S_R)$ as before.
Let $\ol X_R:=f_R^{-1}(\Spec R[t,t^{-1}])$
and $\ol U_R:=\Spec R[x,y,t,t^{-1}]/((1-x^N)(1-y^M)-t)$.
Put $Z_R:=\ol X_R\setminus \ol U_R$. 
By the construction in \S \ref{HG-sect},
$Z_R$ consists of disjoint $(N+M)$-components, and every components
are isomorphic to $\bG_{m,R}=\Spec R[t,t^{-1}]$.
For $(\nu_1,\nu_2)\in\mu_N(R)\times \mu_M(R)$, let
\begin{equation}\label{m-fermat-eq1}
\xi=\xi(\nu_1,\nu_2)=\left\{
\frac{x-1}{x-\nu_1},\frac{y-1}{y-\nu_2}
\right\}\in K^M_2(\O(\ol U_R))
\end{equation}
be a Milnor symbol in $K_2$.
Let $K_2(\ol U_R)^{(2)}\os{\partial}{\lra} K_1(Z_R)^{(1)}=(R[t,t^{-1}]^\times)^\op\ot\Q$
be the boundary map where $K_i(-)^{(j)}$ denotes the Adams weight piece, 
which is explicitly described by
\[
\{f,g\}\longmapsto (-1)^{\ord_P(f)\ord_P(g)}\frac{g^{\ord_P(f)}}{f^{\ord_P(g)}}\bigg|_P,\quad 
P\in Z_R.
\]
It is a simple exercise to show that
$\partial(\xi)=0$ and hence $\xi$ lies in the image of $K_2(\ol X_R)^{(2)}$.
Since $K_2(\ol X_R)^{(2)}\to K_2(\ol U_R)^{(2)}$ is injective as $K_2(\bG_{m,R})^{(1)}=0$, 
we have an element in $K_2(\ol X_R)^{(2)}$ which we also write by $\xi$.
Let $\mathrm{dlog}:K_2(\ol X_R)^{(2)}\to \vg(\ol X,\Omega^2_{\ol X_R/R})\ot\Q$
be the dlog map which is given by $\{f,g\}\mapsto\frac{df}{f}\wedge\frac{dg}{g}$.
One immediately has
\begin{equation}\label{m-fermat-eq2}
\mathrm{dlog}(\xi)=N^{-1}M^{-1}
\sum_{i=1}^{N-1}\sum_{j=1}^{M-1}(1-\nu^{-i}_1)(1-\nu^{-j}_2)\frac{dt}{t}\wedge\omega_{i,j}.
\end{equation}

\medskip

Let $p>\max(N,M)$ be a prime. Let 
$W=W(\ol\F_p)$ be the Witt ring and $K:=\Frac W$ the fractional field. 
Fix an embedding $R\hra W$.
Write $Y:=Y\times_RW$, $\ol X:=\ol X_R\times _RW$ etc.
Let $D_i=f^{-1}(t=i)$ be the fiber at $t=i$ for 
$i\in \{0,1,\infty\}$.
By Lemma \ref{int.model-lem}, the morphism
\[
f:(Y,D_0\cup (D_1)_{\mathrm{red}}\cup (D_{\infty})_{\mathrm{red}})\lra (\P^1,\{0,1,\infty\})
\]
of log schemes is smooth where $(-)_{\mathrm{red}}$ denotes the reduced part,
and a pair $(V,D)$ denotes the log scheme whose log structure is defined by the divisor $D$.
The boundary $Z\subset \ol X$ consists of sections $\{X_0-\zeta X_1=Y_0=0\}$ and 
$\{Y_0-\zeta' Y_1=X_0=0\}$ in the equation \eqref{HGC}.
One easily sees that the closure $\ol Z\subset Y$ in $Y$ also consists of sections
which are disjoint, and each section intersects
with regular reduced locus of $D_i$ transversally for every $i\in \{0,1,\infty\}$.
For $c\in 1+pW$, there is a $p$-th Frobenius $\sigma$ on the weak completion 
$\O(S)^\dag$ given by $\sigma(t)=ct^p$ compatible with the Frobenius on $W$. 
This setting is under the setting in \cite[\S 4.1]{AM}, so that
the comparison map
\[
\O(S)^\dag_K\ot_{\O(S_K)}H^i_\dR(U_K/S_K)\lra  H^i_\rig(U_{\ol\F_p}/S_{\ol\F_p})
\]
is bijective, and the symbol map
\[
[-]_{U/S}:K_2^M(\O(U))\lra \Ext^1_{\FilFMIC(S,\sigma)}(\O_S,H^1(U/S)(2))
\]
is defined.
Since the Milnor symbol $\xi\in K_2^M(\O(U))$ has no boundary at $X\setminus U$,
the symbol map $[-]_{U/S}$ also defines a $1$-extension
\begin{equation}\label{syn-reg-ext}
0\lra H^1(X/S)(2)\lra M_\xi(X/S)\lra \O_S\lra 0
\end{equation}
in the exact category $\FilFMIC(S,\sigma)$ (\cite[Prop.4.3]{AM}).
Let $e_\xi\in \mathrm{Fil}^0M_\xi(X/S)_\dR$ be the unique lifting of $1\in \O_S(S)$.
Let $\ve_k^{(i,j)}(t)$ and $E_k^{(i,j)}(t)$ be defined by
\begin{align}
e_\xi-\Phi(e_\xi)
&=N^{-1}M^{-1}
\sum_{i=1}^{N-1}\sum_{j=1}^{M-1}(1-\nu^{-i}_1)(1-\nu^{-j}_2)[\ve^{(i,j)}_1(t)\omega_{i,j}
+\ve^{(i,j)}_2(t)\eta_{i,j}]\label{fermat-e-eq1}\\
&=N^{-1}M^{-1}
\sum_{i=1}^{N-1}\sum_{j=1}^{M-1}(1-\nu^{-i}_1)(1-\nu^{-j}_2)[E^{(i,j)}_1(t)\wt\omega_{i,j}
+E^{(i,j)}_2(t)\wt\eta_{i,j}].\label{fermat-e-eq2}
\end{align}
Notice that $\ve_k^{(i,j)}(t)$ and $E^{(i,j)}_k(t)$ 
depend on the choice of the Frobenius $\sigma$.
The relation between $\ve_k^{(i,j)}(t)$ and $E^{(i,j)}_k(t)$ is explicitly given by
\begin{align}
\ve^{(i,j)}_1(t)&=E^{(i,j)}_1(t)F_{a_i,b_j}(t)^{-1}-t(1-t)^{a_i+b_j}F'_{a_i,b_j}(t)E_2^{(i,j)}(t)
\label{syn-reg-eq3}
\\
\ve^{(i,j)}_2(t)&=-b_jt(1-t)^{a_i+b_j}F_{a_i,b_j}(t)E_2^{(i,j)}(t)\label{syn-reg-eq4}.
\end{align}
By the definition $\ve^{(i,j)}_k(t)$ are automatically overconvergent functions,
\[\ve^{(i,j)}_k(t)\in K[t,(t-t^2)^{-1}]^\dag.\]
On the other hand, since $F'_{a_i,b_j}(t)/F_{a_i,b_j}(t)$ is a convergent function 
(cf. \eqref{uroot-thm-eq1}), so is
$E^{(i,j)}_1(t)/F_{a_i,b_j}(t)$,
\begin{equation}\label{syn-reg-eq5}
\frac{E^{(i,j)}_1(t)}{F_{a_i,b_j}(t)}\in K\langle
t,(t-t^2)^{-1},h(t)^{-1}\rangle,\quad h(t):=\prod_{m=0}^s [F_{a_i^{(m)},b_j^{(m)}}(t)]_{<p}
\end{equation}
where $s\geq 0$ is the minimal
integer such that $a_i^{(s+1)}=a_i$ and $b_j^{(s+1)}=b_j$.

\medskip

The following is the main theorem in this paper, which provides
a geometric aspect of 
$\cF^{(\sigma)}_{a,b}(t)$ the $p$-adic hypergeometric function
of logarithmic type defined in \S \ref{pHGlog-defn}.

\begin{thm}\label{fermat-main1}
Suppose $p>\max(N,M)$.
We have 
\begin{equation}\label{main-thm1-eq1}
\frac{E_1^{(i,j)}(t)}{F_{a_i,b_j}(t)}=-\cF^{(\sigma)}_{a_i,b_j}(t).
\end{equation}
Hence
\[
e_\xi-\Phi(e_\xi)\equiv 
\sum_{i=1}^{N-1}\sum_{j=1}^{M-1}\frac{(1-\nu^{-i}_1)(1-\nu^{-j}_2)}{NM}\cF^{(\sigma)}_{a_i,b_j}(t)
\omega_{i,j}
\mod \sum_{i,j}K\langle t,(t-t^2),h(t)^{-1}\rangle e^{\mathrm{unit}}_{i,j}.\]
\end{thm}
\begin{pf}
The Frobenius  $\sigma$ extends on $K((t))$, and 
$\Phi$ also extends on \[K((t))\ot H^1_\dR(X/S)\cong
H^1_{\text{log-crys}}((\cY_{\ol\F_p},D_{\ol\F_p})/(\Delta_W,0))\ot_{W[[t]]}K((t))
\] in the natural way where the isomorphism follows from \eqref{log-crys-isom}.
Apply the Gauss-Manin connection $\nabla$ on \eqref{fermat-e-eq2}.
Since $\nabla\Phi=\Phi\nabla$ and $\nabla(e_\xi)=-\mathrm{dlog}\xi$ (\cite[(2.30)]{AM}),
we have
\begin{align}\label{main-thm1-eq2-left}
&(1-\Phi)\left(-N^{-1}M^{-1}\sum_{i=1}^{N-1}\sum_{j=1}^{M-1}(1-\nu^{-i}_1)(1-\nu^{-j}_2)
\frac{dt}{t}\omega_{i,j}\right)\\
=&N^{-1}M^{-1}\sum_{i=1}^{N-1}\sum_{j=1}^{M-1}(1-\nu^{-i}_1)(1-\nu^{-j}_2)
\nabla(E^{(i,j)}_1(t)\wt\omega_{i,j}+E^{(i,j)}_2(t)\wt\eta_{i,j})\label{main-thm1-eq2-right}
\end{align}
by \eqref{m-fermat-eq2}.
Let $\Phi_{X/S}$ denote the $p$-th Frobenius on $H^1_\rig(X_{\ol\F_p}/S_{\ol\F_p})$.
Then the $\Phi$ on $H^1_\rig(X/S)(2)$ agrees with $p^{-2}\Phi_{X/S}$ by definition of Tate twists. 
It follows from Proposition \ref{frobenius-thm} that we have
\[
\Phi_{X/S}(\wt\omega_{p^{-1}i,p^{-1}j})\equiv p\wt\omega_{i,j}\mod K((t))\wt\eta_{i,j}
\]
where $m\in \{1,\ldots,N-1\}$ with $pm\equiv n$ mod $N$.
Therefore
\[
\eqref{main-thm1-eq2-left}\equiv-
N^{-1}M^{-1}\sum_{i=1}^{N-1}\sum_{j=1}^{M-1}(1-\nu^{-i}_1)(1-\nu^{-j}_2)
(F_{a_i,b_j}(t)-F_{a_i^{(1)},b_j^{(1)}}(t^\sigma))\frac{dt}{t}\wt\omega_{i,j}
\]
modulo $\sum_{i,j} K((t))\wt\eta_{i,j}$.
On the other hand, 
\[
\eqref{main-thm1-eq2-right}\equiv
N^{-1}M^{-1}\sum_{i=1}^{N-1}\sum_{j=1}^{M-1}(1-\nu^{-i}_1)(1-\nu^{-j}_2)
t\frac{d}{dt}(E^{(i,j)}_1(t))\cdot\frac{dt}{t}\wt\omega_{i,j}
\mod \sum_{i,j}K((t))\wt\eta_{i,j}
\]
by Proposition \ref{fermat-GM}.
We thus have
\begin{equation}\label{main-thm1-eq3}
t\frac{d}{dt}E_1^{(i,j)}(t)=-F_{a_i,b_j}(t)+F_{a_i^{(1)},b_j^{(1)}}(t^\sigma),
\end{equation}
and hence 
\[
E_1^{(i,j)}(t)=-\left(C+\int_0^tF_{a_i,b_j}(t)-F_{a_i^{(1)},b_j^{(1)}}(t^\sigma)\frac{dt}{t}\right)
\]
for some constant $C\in K$.
We determine the constant $C$ in the following way.
Firstly $E_1^{(i,j)}(t)/F_{a_i,b_j}(t)$ is a convergent function by \eqref{syn-reg-eq5}.
If $C=\psi_p(a_i)+\psi_p(b_j)+2\gamma_p-p^{-1}\log(c)$, then 
$E_1^{(i,j)}(t)/F_{a_i,b_j}(t)=\cF^{(\sigma)}_{a_i,b_j}(t)$ belongs to
$K\langle t,(t-t^2)^{-1},h(t)^{-1}\rangle$ by Corollary \ref{cong-cor}.
If there is another $C'$ such that 
$E_1^{(i,j)}(t)/F_{a_i,b_j}(t)\in K\langle t,(t-t^2)^{-1},h(t)^{-1}\rangle$,
then it follows
\[
\frac{C-C'}{F_{a_i,b_j}(t)}\in K\langle t,(t-t^2)^{-1},h(t)^{-1}\rangle.
\]
We show that this is impossible.
We recall from Theorem \ref{uroot-thm}
the formula \eqref{uroot-thm-3}
\begin{align*}
\Phi_{X/S}^{s+1}(e^{\mathrm{unit}}_{i,j})
&=\left(\prod_{m=0}^s
\frac{(1-t^{p^m})^{a^{(m)}_i+b^{(m)}_j}}{(1-t^{p^{m+1}})^{a^{(m+1)}_i+b^{(m+1)}_j}}
\cF^{\mathrm{Dw}}_{a^{(m)}_i,b^{(m)}_j}(t^{p^m})\right)e^{\mathrm{unit}}_{i,j}\\
&=\left(\frac{(1-t)^{a_i+b_j}}{(1-t^{p^{s+1}})^{a_i+b_j}}\right)
\frac{F_{a_i,b_j}(t)}{F_{a_i,b_j}(t^{p^{s+1}})}e^{\mathrm{unit}}_{i,j}.
\end{align*}
Iterating $\Phi_{X/S}^{s+1}$ to the above, we have
\[
(\Phi_{X/S}^{s+1})^n(e^{\mathrm{unit}}_{i,j})
=\overbrace{\left(\frac{(1-t)^{a_i+b_j}}{(1-t^{p^{n(s+1)}})^{a_i+b_j}}\right)}^{\mu(t)}
\frac{F_{a_i,b_j}(t)}{F_{a_i,b_j}(t^{p^{n(s+1)}})}e^{\mathrm{unit}}_{i,j}.
\]
Put $q:=p^{n(s+1)}$.
Let $\alpha\in W$ satisfy $\alpha^q=\alpha$ and 
$(\alpha-\alpha^2)h(\alpha)\not\equiv0$ mod $p$.
Then the evaluation $\mu(\alpha)$ is a root of unity. 
Suppose $g(t):=F_{a_i,b_j}(t)^{-1}\in K\langle t,(t-t^2)^{-1},h(t)^{-1}\rangle$.
Let $g(t)=(t-\alpha)^kg_0(t)$ with $g_0(t)\in 
K\langle t,(t-t^2)^{-1},h(t)^{-1}\rangle$ such that $g_0(\alpha)\ne0$.
Then we have 
\[
\frac{F_{a_i,b_j}(t)}{F_{a_i,b_j}(t^q)}\bigg|_{t=\alpha}
=\frac{(t^q-\alpha)^kg_0(t^q)}{(t-\alpha)^kg_0(t)}\bigg|_{t=\alpha}
=(q\alpha^{q-1})^k=q^k.
\]
Since the first evaluation is a unit in $W$, we have $k=0$.
Thus
an eigen value of $(\Phi_{X/S}^{s+1})^n|_{t=\alpha}$ is a root of unity.
This contradicts with the Weil-Riemann hypothesis.
\end{pf}

\begin{thm}[Syntomic Regulator Formula]\label{main-thm3}
Let $\alpha\in W$ such that $\alpha\not\equiv0,1$ mod $p$.
Let $\sigma_{\alpha}$ be the Frobenius given by $t^\sigma=F(\alpha)\alpha^{-p}t^p$
where $F$ is the Frobenius on $W$.
Let $X_\alpha$ be the fiber at $t=\alpha$ ($\Leftrightarrow$ $\l=1-\alpha$), 
which is a smooth projective variety over $W$ of relative dimension one.
Let
\[
\reg_\syn:K_2(X_\alpha)\lra H^2_\syn(X_\alpha,\Q_p(2))\cong H^1_\dR(X_\alpha/K)
\]
be the syntomic regulator map.
Then
\[
\reg_\syn(\xi|_{X_\alpha})=-
N^{-1}M^{-1}
\sum_{i=1}^{N-1}\sum_{j=1}^{M-1}(1-\nu^{-i}_1)(1-\nu^{-j}_2)
[\ve^{(i,j)}_1(\alpha)\omega_{i,j}
+\ve^{(i,j)}_2(\alpha)\eta_{i,j}].
\]
\end{thm}
\begin{pf}
This follows from \cite[Theorem 4.4]{AM}.
\end{pf}
\begin{cor}\label{main-thm4}
Let the notation and assumption be as in Theorem \ref{main-thm3}. Suppose further that
$h(\alpha)\not\equiv0$ mod $p$ where $h(t)$ is as in \eqref{syn-reg-eq5}.
Let $e_{N-i,M-j}^{\text{\rm unit}}$ be as in Theorem \ref{uroot-thm}, and
$Q: H^1_\dR(X_\alpha/K)\ot H^1_\dR(X_\alpha/K)\to H^2_\dR(X_\alpha/K)\cong K$ 
the cup-product pairing.
Then we have
\[
Q(\reg_\syn(\xi|_{X_\alpha}), e_{N-i,M-j}^{\text{\rm unit}})
=-N^{-1}M^{-1}
(1-\nu^{-i}_1)(1-\nu^{-j}_2)
\cF_{a_i,b_j}^{(\sigma_\alpha)}(\alpha)
Q(\omega_{i,j}, e_{N-i,M-j}^{\text{\rm unit}}).
\]
\end{cor}
\begin{pf}
Noticing $Q(e_{i,j}^{\text{\rm unit}},e_{N-i,M-j}^{\text{\rm unit}})=0$ by \eqref{uroot-thm-eq3},
this is immediate from Theorem \ref{main-thm3}.
\end{pf}

\subsection{Syntomic regulator of the Ross symbols of Fermat curves}
\label{fermatcurve-sect}
We apply Theorem \ref{fermat-main1} to the study of the syntomic regulator 
of the {\it Ross symbol} \cite{ross2}
\[
\{1-z,1-w\}\in K_2(F)\ot\Q
\]
of the (projective smooth) Fermat curve $F$ defined by an affine equation
$z^N+w^M=1$ over a field $K$ of characteristic zero.
The group $\mu_N\times \mu_M$ acts on $F$ by $(\ve_1,\ve_2)\cdot (z,w)=(\ve_1z,\ve_2w)$.
Let $H^1_\dR(F/K)(i,j)$ denote the subspace on which $(\ve_1,\ve_2)$ acts by multiplication
by $\ve_1^i\ve_2^j$.
Let
\[
I=\left\{(i,j)\in\Z^2\mid1\leq i\leq N-1,1\leq j\leq M-1,\,\frac{i}{N}+\frac{j}{M}\ne1\right\},
\]
then 
\[
H^1_\dR(F/K)=\bigoplus_{(i,j)\in I}H^1_\dR(F/K)(i,j).
\]
Each eigen space $H^1_\dR(F/K)(i,j)$ 
is one-dimensional with basis 
$z^{i-1}w^{j-M}dz=-N^{-1}Mz^{i-N}w^{j-1}dw$, and 
\[
H^1_\dR(F/K)(i,j)\subset \vg(F,\Omega^1_{F/K})\quad\Longleftrightarrow\quad
\frac{i}{N}+\frac{j}{M}<1
\]
(e.g. \cite[\S 2]{gross}).
In particular, the genus of $F$ is $1+\frac{1}{2}(NM-N-M-\gcd(N,M))$.

\begin{thm}\label{fermat-main2}
Let $p>\max(N,M)$ be a prime and $W=W(\ol \F_p)$ the Witt ring and $K=\Frac(W)$.
Let $F$ be the Fermat curve defined by an affine equation $z^N+w^M=1$
that is smooth and projective over $W$.  
Let
\[
\reg_\syn:K_2(F)\ot\Q\lra H^2_\syn(F,\Q_p(2))\cong H^1_\dR(F/K)
\]
be the syntomic regulator map and let $A^{(i,j)}\in K$ be defined by
\[
\reg_\syn(\{1-z,1-w\})=
\sum_{(i,j)\in I} A^{(i,j)}M^{-1}z^{i-1}w^{j-M}dz.
\]
Suppose that $(i,j)\in I$ satisfies
\begin{equation}\label{fermat-main2-eq1}
\mbox{\rm(i) }
\frac{i}{N}+\frac{j}{M}<1,\quad
\mbox{\rm(ii) }
[F_{\frac{i}{N},\frac{j}{M}}(t)]_{<p^n}|_{t=1}\not\equiv0\mod p,\quad\forall n\geq1,
\end{equation}
where $f(t)|_{t=a}$ denotes the evaluation $f(a)$ at $t=a$.
Then we have
\begin{equation}\label{fermat-main2-eq2}
A^{(i,j)}=\cF^{(\sigma)}_{\frac{i}{N},\frac{j}{M}}(1)
\end{equation}
where $\sigma=\sigma_1$ (i.e. $\sigma(t)=t^p$).
\end{thm}
The following lemma gives a sufficient condition for that 
the conditions \eqref{fermat-main2-eq1} are satisfied.
\begin{lem}\label{fermat-main2-lem1}
\begin{enumerate}
\item[$(1)$]
Let $a,b\in \Z_p$. Then
$[F_{a,b}(t)]_{<p^n}|_{t=1}\not\equiv0$ mod $p$ for all $n\geq1$ if and only if
$[F_{a^{(k)},b^{(k)}}(t)]_{<p}|_{t=1}\not\equiv0$ mod $p$ for all $k\geq0$ where $a^{(k)}$ denotes the
Dwork $k$-th prime.
\item[$(2)$]
Let $a_0,b_0\in\{0,1,,\ldots,p-1\}$ satisfy $a\equiv -a_0$ and $b\equiv -b_0$ mod $p$.
Then 
\[
[F_{a,b}(t)]_{<p}|_{t=1}\equiv \frac{\Gamma(1+a_0+b_0)}{\Gamma(1+a_0)\Gamma(1+b_0)}
=\frac{(a_0+b_0)!}{a_0!b_0!}
\mod p.
\]
In particular
\[
[F_{a,b}(t)]_{<p}|_{t=1}\not\equiv 0
\quad\Longleftrightarrow\quad 
a_0+b_0\leq  p-1.
\]
\item[$(3)$]
Suppose that $N|(p-1)$ and $M|(p-1)$. Then for any 
$(i,j)$ such that $0<i<N$ and $0<j<M$ and $i/N+j/M<1$, 
the conditions \eqref{fermat-main2-eq1} hold.
\end{enumerate}
\end{lem}
\begin{pf}
(1) is a consequence of the Dwork congruence \eqref{Dwork-congruence}.
We show (2). Obviously $[F_{a,b}(t)]_{<p}\equiv [F_{-a_0,-b_0}(t)]_{<p}$ mod $p\Z_p[t]$, and
$[F_{-a_0,-b_0}(t)]_{<p}=F_{-a_0,-b_0}(t)$ as $a_0$ and
$b_0$ are non-positive integers greater than $-p$.
Then apply Gauss' formula (e.g. \cite{NIST} 15.4.20)
\[
{}_2F_1\left({a,b\atop c};1\right)=
\frac{\Gamma(c)\Gamma(c-a-b)}{\Gamma(c-a)\Gamma(c-b)},\quad \mathrm{Re}(c-a-b)>0.
\]
To see (3), letting $a=i/N$ and $b=j/M$, we note that 
$a^{(k)}=a$, $b^{(k)}=b$ and $a_0=i(p-1)/N$, $b_0=j(p-1)/M$.
Then the condition \eqref{fermat-main2-eq1} (ii) follows by (1) and (2).
\end{pf}

\bigskip

\noindent{\it Proof of Theorem \ref{fermat-main2}}.
We show that the theorem is reduced to the case $M=N$. Let $L$ be the least common multiple of $M,N$,
and $F_1$ the Fermat curve defined by an affine equation $F_1:z^L_1+w_1^L=1$.
There is a finite surjective map $\rho:F_1\to F$ given by $
\rho^*(z)=z_1^A$ $\rho^*(w)=w_1^B$ where $AN=BM=L$.
There is a commutative diagram
\[
\xymatrix{
K_2(F_1)\ot\Q\ar[r]\ar[d]_{\rho_*}&H^2_\syn(F_1,\Q_p(2))\ar[d]^{\rho_*}\\
K_2(F)\ot\Q\ar[r]&H^2_\syn(F,\Q_p(2))
}
\]
with surjective vertical arrows. It is a simple exercise to show that $\rho_*\{1-z_1,
1-w_1\}=\{1-z.1-w\}$ and $\rho_*(L^{-1}z_1^{i-1}w_1^{j-L}dz_1)
=M^{-1}z^{i'-1}w^{j'-M}dz$ if $(i,j)=(i'A,j'B)$ and $=0$ otherwise.
Thus the theorem for $F$ can be deduced from the theorem for $F_1$.

\medskip

We assume $N=M$ until the end of the proof.
Let $f:Y_s\to\P^1_W$ be the curve \eqref{HGC},
which has bad fibers at $t=0,1,\infty$.
Let $\l:=1-t$ be a new parameter, and let $\l_0^N=\l$.
Let $\ol S_0:=\Spec W[\l_0,(1-\l^N_0)^{-1}]\to \P^1_W$
and $S_0:=\Spec W[\l_0,\l_0^{-1},(1-\l^N_0)^{-1}]\subset \ol S_0$.
Let $\ol X_s:=Y_s\times_{\P^1_W}\ol S$. Then
$\ol X_s$ has a unique singular point $(x,y,\l_0)=(0,0,0)$ in an affine open set
\[U_s=\Spec W[x,y,\l_0,(1-\l^N_0)^{-1}]/(x^Ny^N-x^N-y^N-\l_0^N)\subset \ol X_s.\]
Let $\ol X_0\to \ol X_s$ be the blow-up at $(x,y,\l_0)=(0,0,0)$. Then $\ol X_0\to \Spec W$ is smooth, and the morphism
\begin{equation}\label{fermat-main2-lem2}
\xymatrix{
f_0:\ol X_0\ar[r]& \ol S_0
}
\end{equation}
is projective flat
such that $X_0:=f_0^{-1}(S_0)\to S_0$ smooth and $f_0$ has a 
semistable reduction at $\l_0=0$.
The fiber $Z:=f^{-1}(\l_0=0)$ is a reduced divisor with two irreducible components
$F$ and $E$ where $F$ is the proper transform of the curve $x^Ny^N-x^N-y^N=0$
$\Leftrightarrow$ $z^N+w^N=1$ $(z:=x^{-1},w:=y^{-1})$, and $E$ is the exceptional curve. Both curves are
isomorphic to the Fermat curve $u^N+v^N=1$. 
Moreover $E$ and $F$ intersects transversally at $N$-points.

We recall the $K_2$-symbols $\xi(\nu_1,\nu_2)$ in \eqref{m-fermat-eq1}.
We think them to be elements of $K_2(\ol X_0)\ot\Q$, and
put
\[
\Xi:=\sum_{(\nu_1,\nu_2)\in \mu_N\times \mu_N}\xi(\nu_1,\nu_2)
=
\left\{\frac{(x-1)^N}{x^N-1},\frac{(y-1)^N}{y^N-1}\right\}\in K_2(\ol X_0)\ot\Q.
\]
Then the restriction of $\Xi$ on $F$ is 
\begin{align*}
\Xi|_F
&=\{(1-z)^N,(1-w)^N\}-\{(1-z)^N,1-w^N\}-\{1-z^N,(1-w)^N\}+\{1-z^N,1-w^N\}\notag\\
&=\{(1-z)^N,(1-w)^N\}-\{(1-z)^N,z^N\}-\{w^N,(1-w)^N\}+\{w^N,1-w^N\}\notag\\
&=N^2\{1-z,1-w\}.
\end{align*}
This is the Ross symbol.
We thus have
\begin{equation}\label{Xi-eq1}
N^2\reg_\syn(\{1-z,1-w\})=
\sum_{(\nu_1,\nu_2)\in \mu_N\times \mu_N}\reg_\syn(\xi(\nu_1,\nu_2)|_F).
\end{equation}
Write $\xi=\xi(\nu_1,\nu_2)$ simply.
Let $\sigma$ be the $p$-th Frobenius on $W[\l_0,\l_0^{-1}]$ such that $\sigma(t)=t^p$
$\Leftrightarrow$ $\sigma(\l_0)=(1-(1-\l_0^N)^p)^{\frac1N}$.
Recall \eqref{fermat-e-eq1} and \eqref{fermat-e-eq2},
\begin{align*}
e_\xi-\Phi_\sigma(e_\xi)
&=N^{-2}\sum_{1\leq i,j\leq N-1}(1-\nu^{-i}_1)(1-\nu^{-j}_2)[
\ve^{(i,j)}_{1,\sigma}(t)\omega_{i,j}
+\ve^{(i,j)}_{2,\sigma}(t)\eta_{i,j}]\\
&=N^{-2}\sum_{1\leq i,j\leq N-1}
(1-\nu^{-i}_1)(1-\nu^{-j}_2)[
E^{(i,j)}_{1,\sigma}(t)\wt\omega_{i,j}
+E^{(i,j)}_{2,\sigma}(t)\wt\eta_{i,j}]
\end{align*}
where we write ``$(-)_\sigma$'' 
to emphasize that they depend on $\sigma$.
Let $\tau$ be the $p$-th Frobenius on $W[[\l_0]]$ such that $\tau(\l_0)=\l_0^p$.
Let
\begin{equation}\label{fermat-main2-eq5}
e_\xi-\Phi_\tau(e_\xi)
=N^{-2}\sum_{1\leq i,j\leq N-1}
(1-\nu^{-i}_1)(1-\nu^{-j}_2)[\ve^{(i,j)}_{1,\tau}(\l)\omega_{i,j}
+\ve^{(i,j)}_{2,\tau}(\l)\eta_{i,j}]
\end{equation}
be defined in the same way.
This is related to \eqref{Xi-eq1} in the following way.
Let
$\Delta:=\Spec W[[\l_0]]\to \ol S_0$, and  $\cX:=f^{-1}(\Delta)$.
We have the syntomic regulator
\[
\reg_\syn(\xi)\in H^2_\syn(\cX,\Z_p(2))
\]
in the syntomic cohomology group. 
We endow the log structure on $\Delta$ (resp. $\cX$)
defined by the divisor $O=\Spec W[[\l_0]]/(\l_0)$ (resp. $E+F$)
which is denoted by the same notation $O$ (resp. $E+F$).
Let $\omega_{\cX/\Delta}$ be the log de Rham complex for $(\cX,E+F)/(\Delta,O)$.
Recall the log syntomic cohomology groups (e.g. \cite[\S 2]{T})
\[
H^i_\syn((X,M),\Z_p(j))
\]
of a log scheme $(X,M)$ satisfying several conditions 
(all log schemes appearing in this proof satisfy them).
Moreover one can further define the syntomic cohomology groups
$H^i_\syn((\cX,E+F)/(\Delta,O,\tau),\Z_p(j))$ following the construction in \cite[\S 3.1]{AM}, 
where we note that 
$\tau$ induces the $p$-th Frobenius on $(\Delta,O)$ (while so does not $\sigma$).
Let 
\[
\rho_\tau:H^2_\syn(\cX,\Z_p(2))
\lra H^2_\syn((\cX,E+F)/(\Delta,O,\tau),\Z_p(2))
\os{\cong}{\longleftarrow}H^1_\zar(\cX,\omega^\bullet_{\cX/\Delta})
\]
be the composition of natural maps.
We endow the log structure on $F$ defined by the divisor $T:=E\cap F$ 
which is denoted by $T$. Put $U:=F\setminus T$.
Let $\omega^\bullet_F:=\Omega^\bullet_{F/W}(\log T)$ 
the log de Rham complex for $(F,T)/W$.
Let $\iota:F\hra\cX$ be the closed immersion.
Then there is a commutative diagram
\[
\xymatrix{
H^2_\syn(\cX,\Z_p(2))\ar[d]_{\rho_\tau}\ar[r]& H^2_\syn((F,T),\Z_p(2))\\
H^1_\zar(\cX,\omega^\bullet_{\cX/\Delta})\ar[r]^-{\iota^*}\ar[d]_\pi
&H^1_\zar(F,\omega^\bullet_{F})\ar[u]_\cong\ar[r]^-{\subset}&
H^1_\dR(U/K)\\
W((\l_0))\ot H^1_\dR(X_0/S_0)
}
\]
and we have
\begin{equation}\label{1-ext-thm1-eq1}
(\iota^*\circ\rho_\tau)(\reg_\syn(\xi))=\reg_\syn(\xi|_F)\in H^1_\dR(U/K).
\end{equation}
Moreover it follows from \cite[Theorem 4.5]{AM} that
\begin{equation}\label{1-ext-thm1-eq2}
(\pi\circ\rho_\tau)(\reg_\syn(\xi))=\Phi(e_\xi)-e_\xi
\in 
W((\l_0))\ot H^1_\dR(X_0/S_0).
\end{equation}
Note that $H^1_\dR(F/K)\to H^1_\dR(U/K)$ is injective, and the above element belongs to
the image of $H^1_\dR(F/K)$, so that we may replace $H^1_\dR(U/K)$ with $H^1_\dR(F/K)$
in \eqref{1-ext-thm1-eq1}.

\medskip

Write $X_{0,K}:=X_0\times_WK$ etc.
Let $\Delta_K:=\Spec K[[\l_0]]$ and $\cX_K:=\ol X_0\times_{\ol S_0} \Delta_K$.
Put 
\[
H_K:=H^1(\cX_K,\omega^\bullet_{\cX_K/\Delta_K})\hra K((\l_0))\ot H^1_\dR(X_0/S_0).
\]

\begin{lem}\label{fermat-main2-lem3}
Put $s:=(a_i+b_j)N$ which is a positive integer.
If $a_i+b_j<1$, then the eigencomponent $H_K(i,j)$ is a free $K[[\l_0]]$-module
of rank two
with a basis $\{\omega_{i,j},\l_0^{s}\eta_{i,j}\}$.
\end{lem}
\begin{pf}
This is proven in the same way as the proof of
Corollary \ref{fermat-GM-cor2}.
\end{pf}

\begin{lem}\label{fermat-main2-lem0}
Let $1\leq i,j\leq N-1$ be integers, and put
$a_i:=1-i/N$ and $b_j:=1-j/N$.
Put
\[
f_n(t)=f_{n,i,j}(t):=-\frac{(1-\nu_1^{-i})(1-\nu_2^{-j})}{N^2}\frac{1}{F_{a_i,b_j}(t)}
\left(\frac{d^{n-1}}{dt^{n-1}}\left(\frac{F_{a^{(1)}_i,b^{(1)}_j}(t)}{t}\right)\right)^\sigma
\]
for $n\in\Z_{\geq 1}$. 
Then
\[
\ve^{(i,j)}_{1,\tau}(\l)-\cF^{(\sigma)}_{a_i,b_j}(t)
=\sum_{n=1}^\infty\frac{(t^\tau-t^\sigma)^n}{n!}
p^{-1}f_n(t)
+b_j^{-1}\frac{F'_{a_i,b_j}(t)}{F_{a_i,b_j}(t)}
\ve^{(i,j)}_{2,\tau}(\l).
\]
\end{lem}
Notice that $f_n(t)$ is a convergent function on the rigion 
$\{[F_{a_i,b_j}(t)]_{<p^n}\not\equiv 0\}$ by
\cite[p.37, Thm. 2, p.45 Lem. 3.4]{Dwork-p-cycle}
\begin{pf}
The relation between $\ve^{(i,j)}_{k,\sigma}(t)$ and $\ve^{(i,j)}_{k,\tau}(t)$ is the following (e.g. \cite[6.1]{EK}, 
\cite[17.3.1]{Ke})
\begin{equation}\label{fermat-main2-eq6}
\Phi_\tau(e_\xi)-\Phi_\sigma(e_\xi)=
\sum_{n=1}^\infty\frac{(t^\tau-t^\sigma)^n}{n!}
\Phi_\sigma\partial^n_t e_\xi
\end{equation}
where $\partial_t=\nabla_{\frac{d}{dt}}$ is the differential operator on $M_\xi(X/S)_\dR$.
By \eqref{m-fermat-eq2},
\[
\partial_t(e_\xi)
=-\sum_{1\leq i,j\leq N-1}\frac{(1-\nu^{-i}_1)(1-\nu^{-j}_2)}{N^2}\frac{1}{t}\omega_{i,j}
=-\sum_{1\leq i,j\leq N-1}\frac{(1-\nu^{-i}_1)(1-\nu^{-j}_2)}{N^2}\frac{F_{a_i,b_j}(t)}{t}\wt\omega_{i,j}.
\] 
Let $\eta_{i,j}^*:=(1-t)^{-a_i-b_j}F_{a_i,b_j}(t)^{-1}\wt\eta_{i,j}\in H^1_\dR(X/S)\ot
K\langle t,(t-t^2)^{-1},h(t)^{-1}\rangle$ where $h(t)=\prod_{m=0}^N
F_{a^{(m)}_i,b^{(m)}_j}(t)$ with $N\gg0$.
By Proposition \ref{fermat-GM}, 
\[
\partial^n_t(e_\xi)
=-\sum_{1\leq i,j\leq N-1}\frac{(1-\nu^{-i}_1)(1-\nu^{-j}_2)}{N^2}
\frac{d^{n-1}}{dt^{n-1}}\left(\frac{F_{a_i,b_j}(t)}{t}\right)\wt\omega_{i,j}+(\cdots)\wt\eta_{i,j}
\]
and hence
\[
\Phi_\sigma\partial^n_t(e_\xi)
\equiv \sum_{1\leq i,j\leq N-1}
p^{-1}f_{n,i,j}(t)\mod K\langle t,(t-t^2)^{-1},h(t)^{-1}\rangle\eta^*_{i,j}\]
by Proposition \ref{frobenius-thm}.
Take the reduction of the both side of \eqref{fermat-main2-eq6} modulo 
$K\langle t,(t-t^2)^{-1},h(t)^{-1}\rangle\eta^*_{i,j}$.
We then have
\[
\ve^{(i,j)}_{1,\tau}(\l)-\ve^{(i,j)}_{1,\sigma}(t)
-b_j^{-1}\frac{F'_{a_i,b_j}(t)}{F_{a_i,b_j}(t)}
(\ve^{(i,j)}_{2,\tau}(\l)-\ve^{(i,j)}_{2,\sigma}(t))
=\sum_{n=1}^\infty\frac{(t^\tau-t^\sigma)^n}{n!}
p^{-1}f_n(t).
\]
On the other hand,
\[
\ve^{(i,j)}_{1,\sigma}(t)=\cF^{(\sigma)}_{a_i,b_j}(t)+b_j^{-1}\frac{F'_{a_i,b_j}(t)}{F_{a_i,b_j}(t)}
\ve^{(i,j)}_{2,\sigma}(t)
\]
by \eqref{syn-reg-eq3}, \eqref{syn-reg-eq4} and Theorem \ref{fermat-main1}.
Hence
\[
\ve^{(i,j)}_{1,\tau}(\l)-\cF^{(\sigma)}_{a_i,b_j}(t)
=\sum_{n=1}^\infty\frac{(t^\tau-t^\sigma)^n}{n!}
p^{-1}f_n(t)
+b_j^{-1}\frac{F'_{a_i,b_j}(t)}{F_{a_i,b_j}(t)}
\ve^{(i,j)}_{2,\tau}(\l)
\]
as required.
\end{pf}
\begin{lem}\label{fermat-main2-lem4}
If $a_i+b_j<1$, then $\ord_{\l=0}(\ve_{1,\tau}^{(i,j)}(\l))\geq 0$ and 
$\ord_{\l=0}(\ve_{2,\tau}^{(i,j)}(\l))\geq 1$.
\end{lem}
\begin{pf}
Since $e_\xi-\Phi_\tau(e_\xi)\in H_K$,
we have
\[
\ve_{1,\tau}^{(i,j)}(\l)\omega_{i,j}+\ve_{2,\tau}^{(i,j)}(\l)\eta_{i,j}\in H_K(i,j).
\]
If $a_i+b_j<1$, then this means
\[
\ve_{1,\tau}^{(i,j)}(\l_0^N)\, , \l_0^{-s}\ve_{2,\tau}^{(i,j)}(\l_0^N)\in K[[\l_0]].
\]
by Lemma \ref{fermat-main2-lem3}. Since $s=(a_i+b_j)N<N$, the assertion follows.
\end{pf}
\begin{lem}\label{fermat-main2-lem5}
If $a_i+b_j<1$ and $F_{a_i,b_j}(1)_{<p^n}\not\equiv 0$ mod $p$ for all $n\geq 1$, then 
\[
\ve_{1,\tau}^{(i,j)}(0)=\cF^{(\sigma)}_{a_i,b_j}(1)
\]
where the left hand side denotes the evaluation at $\l=0$ ($\Leftrightarrow$ $t=1$)
and the right hand side denotes the evaluation at $t=1$.
Note that the left value is defined by Lemma \ref{fermat-main2-lem4}.
\end{lem}
\begin{pf}
This is straightforward from Lemma \ref{fermat-main2-lem0} on noticing that 
$F'_{a_i,b_j}(t)/F_{a_i,b_j}(t)$ and $f_n(t)$ are convergent at $t=1$ by 
\cite[p.45, Lem. 3.4 ]{Dwork-p-cycle} under the condition that
$[F_{a_i,b_j}(t)]_{<p^n}|_{t=1}\not\equiv 0$ mod $p$ for all $n\geq 1$.
\end{pf}

\medskip

We finish the proof of Theorem \ref{fermat-main2}.
Let $(i,j)$ satisfy $a_i+b_j<1$.
Let 
\[\rho_\tau(\reg_\syn(\xi))(i,j)\in H_K(i,j)\] be the eigencomponent of
$\psi_\tau(\reg_\syn(\xi))$, which agrees with
\[
-N^{-2}(1-\nu^{-i}_1)(1-\nu^{-j}_2)[
\ve^{(i,j)}_{1,\tau}(t)\omega_{i,j}
+\ve^{(i,j)}_{2,\tau}(t)\eta_{i,j}]
\]
by \eqref{fermat-main2-eq5} and \eqref{1-ext-thm1-eq2}.
It is straightforward to see $\iota^*(\omega_{i,j})=Nz^{N-i-1}w^{-j}dz.$
Then 
\begin{align*}
&\iota^*[\ve^{(i,j)}_{1,\tau}(\l)\omega_{i,j}+\ve^{(i,j)}_{2,\tau}(\l)\eta_{i,j}]\\
=&
\ve_{1,\tau}^{(i,j)}(\l_0^L)|_{\l_0=0}\cdot \iota^*(\omega_{i,j})+
(\l_0^{-s}\ve_{2,\tau}^{(i,j)}(\l_0^L))|_{\l_0=0}\cdot \iota^*(\l_0^s\eta_{i,j})]
&\text{(Lemma \ref{fermat-main2-lem3})}\\
=&
\ve_{1,\tau}^{(i,j)}(0)\cdot Nz^{N-i-1}w^{-j}dz
&\text{(Lemma \ref{fermat-main2-lem4} and $s<L$)}\\
=&
\cF_{a_i,b_j}^{(\sigma)}(1)\cdot Nz^{N-i-1}w^{-j}dz
& \text{(Lemma \ref{fermat-main2-lem5})}.
\end{align*}
Therefore 
\[\reg_\syn(\xi|_F)
=-N^{-2}
\sum_{a_i+b_j<1}(1-\nu^{-i}_1)(1-\nu^{-j}_2)
\cF_{a_i,b_j}^{(\sigma)}(1)Nz^{N-i-1}w^{-j}dz+\sum_{a_i+a_j> 1}(-)
\]
by \eqref{1-ext-thm1-eq1}.
Taking the summation over $(\nu_1,\nu_2)\in \mu_N\times \mu_M$, we have
\[
N^2\reg_{\syn}(\{1-z,1-w\})=-\sum_{a_i+a_j< 1}\cF_{a_i,b_j}^{(\sigma)}(1)Nz^{N-i-1}w^{-j}dz
+\sum_{a_i+a_j> 1}(-)
\]
by \eqref{Xi-eq1}. This finishes the proof of Theorem \ref{fermat-main2}.
 
 \bigskip

In \cite{ross2}, Ross showed the non-vanishing of the Beilinson regulator
\[
\reg_B\{1-z,1-w\}\in H^2_\cD(F,\R(2))\cong H^1_B(F,\R)^{F_\infty=-1}
\]          
of his element in the Deligne-Beilinson cohomology group.
We expect the non-vanishing also in the $p$-adic situation.
\begin{conj}
Under the condition \eqref{fermat-main2-eq1},
$\cF_{\frac{i}{N},\frac{j}{M}}^{(\sigma)}(1)\ne0$.
\end{conj}
By the congruence relation for $\cF^{(\sigma)}_{\ul a}(t)$ (Theorem \ref{cong-thm}), 
the non-vanishing $\cF_{\frac{i}{N},\frac{j}{M}}^{(\sigma)}(1)\ne0$ is equivalent to 
\[
[G_{\frac{i}{N},\frac{j}{M}}^{(\sigma)}(t)]_{<p^n}|_{t=1}\not\equiv 0\mod p^n
\]
for some $n\geq 1$.
A number of computations by computer indicate that this holds
(possibly $n\ne1$). 
Moreover if the Fermat curve has a quotient to an elliptic curve over $\Q$, one can expect
that the syntomic regulator agrees with the special value of the $p$-adic $L$-function
according to the $p$-adic Beilinson conjecture by Perrin-Riou \cite[4.2.2]{Perrin-Riou}.
See Conjecture \ref{ell6-conj} below for detail.
\subsection{Syntomic Regulators of Hypergeometric curves of Gauss type}
\label{gauss-sect}
Let $W=W(\ol\F_p)$ and $K=\Frac W$. 
Let $N\geq 2$, $A,B>0$ be integers such that $0<A,B<N$ and $\gcd(N,A)=\gcd(N,B)=1$.
Let $X_{\gauss,K}\to\Spec K[\l,(\l-\l^2)^{-1}]$ a smooth projective morphism
of relative dimension one whose generic fiber is 
defined from an affine equation
\[
v^N=u^A(1-u)^B(1-\l u)^{N-B}.
\]
We call $X_{\gauss,K}$ the {\it hypergeometric curve of Gauss type} (\cite[\S 2.3]{A}).
The genus of a smooth fiber is $N-1$.
Let $X$ be the hypergeometric curve in \S \ref{fermat-sect} defined by an
affine equation $(1-x^N)(1-y^N)=t$.
Then there is a finite cyclic covering 
\begin{equation}
\rho:X\times_WK\lra X_{\gauss,K},\quad
\begin{cases}
\rho^*(u)=x^{-N}\\
\rho^*(v)=x^{-A}(1-x^{-N})y^{N-B}\\
\rho^*(\l)=1-t
\end{cases}
\end{equation}
of degree $N$ whose Galois group is generated by an automorphism 
$g_{A,B}:=[\zeta_N^B,\zeta_N^{-A}]\in \mu_N(K)\times \mu_N(K)$ 
(see \eqref{fermat-ss} for the notation) where $\zeta_N$ is
a fixed primitive $N$-th root of unity
($\rho$ is a generalization of the Fermat quotient, e.g. \cite[p.211]{gross}).

\medskip

Suppose that $N$ is prime to $p$.
We construct an integral model $X_\gauss$ over $W$ in the following way.
The cyclic group $\langle g_{A,B}\rangle$ generated by $g_{A,B}$ acts on $X$,
and it is a free action.
We define
$X_\gauss$ to be the quotient
\[
f:X_\gauss\os{\text{def}}{=}X/\langle g_{A,B}\rangle\lra S=\Spec W[t,(t-t^2)^{-1}].
\]
of $X$ by the cyclic group $\langle g_{A,B}\rangle$. Since it is a free action, $X_\gauss$
is smooth over $W$.
The cyclic group $\mu_N(K)\times \mu_N(K)/\langle g_{A,B}\rangle$ acts on $X_{\gauss,K}$,
and it is generated by an automorphism $h$ given by $(u,v)\mapsto (u,\zeta^{-1}_N v)$.
For an integer $n$, let $V(n)$ be the eigenspace on which $h$ acts by
multiplication by $\zeta_N^n$ for all $\zeta_N\in \mu_N(K)$.
Then the pull-back $\rho^*$ satisfies
\begin{equation}\label{gauss-eq1}
\rho^*(H^1_\dR(X_{\gauss,K}/S_K)(n))=H^1_\dR(X_K/S_K)(nA,nB),\quad
0< n<N
\end{equation}
and the push-forward $\rho_*$ satisfies
\begin{equation}\label{gauss-eq2}
\rho_*(H^1_\dR(X_K/S_K)(i,j))=\begin{cases}
H^1_\dR(X_{\gauss,K}/S_K)(n)&(i,j)\equiv (nA,nB)\text{ mod }N\\
0&\text{otherwise}
\end{cases}
\end{equation}
for $0<i,j<N$.
Put
\begin{equation}\label{gauss-eq3}
\omega_n:=\rho_*(\omega_{nA,nB}),\quad 
\eta_n:=\rho_*(\eta_{nA,nB})
\end{equation}
a basis of $H^1_\dR(X_{\gauss,K}/S_K)(n)$ (see \eqref{rule} for the notation $(-)_{nA,nB}$).
Recall $e^{\mathrm{unit}}_{i,j}$ from 
Theorem \ref{uroot-thm}.
Put
\[
e^{\mathrm{unit}}_n:=\rho_*e^{\mathrm{unit}}_{nA,nB}\in 
H^1_\dR(X_{\gauss,K}/S_K)(n)
\ot K\langle t,(t-t^2)^{-1},h(t)^{-1}\rangle
\] 
for $0<n<N$. 
Notice that
$\rho^*(\omega_n)=N\omega_{nA,nB}$, 
$\rho^*(\eta_n)=N\eta_{nA,nB}$ and 
$\rho^*e^{\mathrm{unit}}_n=Ne^{\mathrm{unit}}_{nA,nB}$ by \eqref{gauss-eq1} and
\eqref{gauss-eq2} together with the fact that $\rho_*\rho^*=N$.
We put
\begin{equation}\label{gauss-eq4}
\xi_\gauss=\xi_\gauss(\nu_1,\nu_2):=\rho_*\xi(\nu_1,\nu_2)\in K_2(X_\gauss)^{(2)}\subset
K_2(X_\gauss)\ot\Q
\end{equation}
where $\xi(\nu_1,\nu_2)$ is as in the beginning of \S \ref{syn-reg-sect}.
Let $\sigma_{\alpha}$ be the Frobenius given by $t^\sigma=ct^p$ with $c\in 1+pW$.
Taking the fixed part of \eqref{syn-reg-ext} by $\langle g_{A,B}\rangle$,
we have a $1$-extension
\[
0\lra H^1(X_\gauss/S)(2)\lra M_{\xi_\gauss}(X_\gauss/S)\lra \O_S\lra 0
\]
in the exact category $\FilFMIC(S)$.
Let $e_{\xi_\gauss}\in \Fil^0M_{\xi_\gauss}(X_\gauss/S)_\dR$ 
be the unique lifting of $1\in \O(S)$. 
\begin{thm}\label{gauss-main1}
Put $a_n:=-nA/N-\lfloor -nA/N\rfloor$ and $b_n:=-nB/N-\lfloor -nB/N\rfloor$.
Let $h(t)=\prod_{m=0}^s[F_{a_n^{(m)},b_n^{(m)}}(t)]_{<p}$ where $s$ is the minimal
integer such that $(a_n^{(s+1)},b_n^{(s+1)})=(a_n,b_n)$ for all $n\in\{1,2,\ldots,N-1\}$.
Then 
\[
e_{\xi_\gauss}-\Phi(e_{\xi_\gauss})\equiv -
\sum_{n=1}^{N-1}\frac{(1-\nu^{-nA}_1)(1-\nu^{-nB}_2)}{N^2}\cF^{(\sigma)}_{a_n,b_n}(t)
\omega_n\]
modulo $\sum_{n=1}^{N-1}K\langle t,(t-t^2),h(t)^{-1}\rangle e^{\mathrm{unit}}_n.$
\end{thm}
\begin{pf}
This is immediate by applying $\rho_*$ on the formula in Theorem \ref{fermat-main1}.
\end{pf}
\begin{cor}\label{gauss-main2}
Suppose $p>N$.
Let $a\in W$ such that $a\not\equiv0,1$ mod $p$.
Let $\sigma_a$ be the Frobenius given by $t^\sigma=F(a)a^{-p}t^p$
where $F$ is the Frobenius on $W$.
Let $X_{\gauss,a}$ be the fiber at $t=a$ ($\Leftrightarrow$ $\l=1-a$), 
which is a smooth projective variety over $W$ of relative dimension one.
Let
\[
\reg_\syn:K_2(X_{\gauss,a})\lra 
H^2_\syn(X_{\gauss,a},\Q_p(2))\cong H^1_\dR(X_{\gauss,a}/K)
\]
be the syntomic regulator map. 
Let $Q$ be the cup-product pairing on $H^1_\dR(X_{\gauss,a}/K)$.
Then
\[
Q(\reg_\syn(\xi_\gauss|_{X_a}),e^{\mathrm{unit}}_{n})=
N^{-2}
(1-\nu^{-nA}_1)(1-\nu^{-nB}_2)
\cF_{a_n,b_n}^{(\sigma_a)}(a)
Q(\omega_n,e^{\mathrm{unit}}_{n}).
\]
\end{cor}
\begin{pf}
This is immediate from Theorem \ref{gauss-main1}
on noticing $Q(e^{\mathrm{unit}}_n,e^{\mathrm{unit}}_m)=0$ for any $n,m$, cf. Corollary \ref{main-thm4}.
\end{pf}

\subsection{Syntomic Regulators of elliptic curves}\label{elliptic-sect}
The methods in \S \ref{syn-reg-sect} work not only for the hypergeometric curves but
also for the elliptic fibrations listed in \cite[\S 5]{A}.
We here give the results together with a sketch of the proof because the discussion
is similar to before.
\begin{thm}\label{elliptic-thm1}
Let $p>5$ be a prime number. Let $W=W(\ol\F_p)$ be the Witt ring and $F:=\Frac(W)$
the fractional field.
Let $f:X\to\P^1\setminus\{0,1,\infty\}$ be the elliptic fibration over $W$
defined by a Weierstrass equation
$3y^2=2x^3-3x^2+1-t$ over $W$. Put $\omega=dx/y$.
Let 
\[\xi:=\left\{
\frac{y-x+1}{y+x-1},
\frac{t}{2(x-1)^3}
\right\}\in K_2(X).
\]
Let $a\in W$ satisfy that
$a\not\equiv0,1$ mod $p$ and $X_a$
has a good ordinary reduction where $X_a$ is the fiber at 
$\Spec W[t]/(t-a)$.
Let $e_{\text{\rm unit}}\in H^1_\dR(X_a/K)$ be the eigen vector of
the unit root with respect to the $p$-th Frobenius $\Phi$.
Let $\sigma_a$ denote the $p$-th Frobenius given by $\sigma_a(t)=F(a)
a^{-p}t^p$.
Then, we have
\[
Q(\reg_\syn(\xi|_{X_a}), e_{\text{\rm unit}})
=\cF_{\frac{1}{6},\frac{5}{6}}^{(\sigma_a)}(a)
Q(\omega,e_{\text{\rm unit}})
\]
\end{thm}
\begin{pf} (sketch).
Let $U\subset X$ be the complement of divisors
$\{y=\pm(x-1)\}$, $\{x=1\}$ and $\{x=\infty\}$ so that
the symbol $\xi$ lies in the image of $K^M_2(\O(U))$.
It is not hard to construct an elliptic fibration $f:Y\to \P^1_W$ such that the 
union the closure of $X\setminus U$ 
and singular fibers of $f$ is a relative simple NCD over $W$ and that
the multiplicity of any component of the singular fibers is at most $6$.
Therefore this setting is under the setting in \cite[\S 4.1]{AM}.
Let $\cE$ be the fiber over the formal neighborhood $\Spec W[[t]]\hra \P^1_W$.
Let $\rho:\bG_m\to \cE$ be the uniformization, and $u$ the uniformizer of $\bG_m$.
Then we have 
\[
\rho^*\omega=F(t)\frac{du}{u}
\]
and a formal power series $F(t)\in W[[t]]$ can be computed by
the usual method of the Picard-Fuchs equation. One sees that $F(t)$ is a solution 
of the differential equation
\[
(t-t^2)\frac{d^2y}{dt^2}+(1-2t)\frac{dy}{dt}-\frac{5}{36}y=0,
\]
and therefore it agrees with the hypergeometric power series 
\[
F_{\frac{1}{6},\frac{5}{6}}(t)={}_2F_1\left({\frac{1}{6},\frac{5}{6}\atop 1};t\right)
\]
up to scalar.
Looking at the residue of $\omega$ at the point $(x,y,t)=(1,0,0)$, one finds that
the constant term of $F(t)$ is $1$, and hence we have
\[
\rho^*\omega=F_{\frac{1}{6},\frac{5}{6}}(t)\frac{du}{u}.
\]
It is straightforward to show
\[
\dlog(\xi)=\frac{dx}{y}\frac{dt}{t}=\omega\wedge\frac{dt}{t}.
\]
Then the rest of the proof goes in the same way as the proof of Theorem \ref{fermat-main1}.
\end{pf}

The following theorems are proven by the same argument 
as in the proof of Theorem \ref{elliptic-thm1}. 
\begin{thm}\label{elliptic-thm2}
Let $p>3$ be a prime and $W=W(\ol \F_p)$ the Witt ring.
Let $f:X\to\P^1\setminus\{0,1,\infty\}$ 
be the elliptic fibration over $W$ defined by a Weierstrass equation
$y^2=x^3+(3x+4t)^2$, and 
\[\xi:=\left\{
\frac{y-3x-4t}{-8t},
\frac{y+3x+4t}{8t}
\right\}.
\]
Then, under the same notation in Theorem \ref{elliptic-thm1}, we have
\[
Q(\reg_\syn(\xi|_{X_a}), e_{\text{\rm unit}})
=3\cF_{\frac{1}{3},\frac{2}{3}}^{(\sigma_a)}(a)Q(
\omega,e_{\text{\rm unit}}).
\]
\end{thm}


\begin{thm}\label{elliptic-thm3}
Let $p>3$ be a prime and $W=W(\ol \F_p)$ the Witt ring.
Let $f:X\to\P^1\setminus\{0,1,\infty\}$ 
be the elliptic fibration defined by a Weierstrass equation
$y^2=x^3-2x^2+(1-t)x$, and 
\[\xi:=\left\{\frac{y-(x-1)}{y+(x-1)},\frac{-t x}{(x-1)^3}\right\}.
\]
Then, under the same notation in Theorem \ref{elliptic-thm1}, we have
\[
Q(\reg_\syn(\xi|_{X_a}), e_{\text{\rm unit}})
=\cF_{\frac{1}{4},\frac{3}{4}}^{(\sigma_a)}(a)
Q(\omega,e_{\text{\rm unit}}).
\]
\end{thm}
\begin{rem}\label{rem.miya}
The syntomic regulator in Theorem \ref{elliptic-thm2}
is also considered in \cite[Theorem 4.8]{AM}, where the authors give a full description
of the element 
$\reg_\syn(\xi|_{X_a})$ in $H^1_\dR(X_a/K)$ (not only the cup-product with a unit root vector).
\end{rem}
\section{$p$-adic Beilinson conjecture for elliptic curves over $\Q$}\label{weakB-sect}
\subsection{Statement}
The Beilinson regulator is a generalization of Dirichlet's regulators of number fields.
in higher $K$-groups of varieties.
He conjectured the formulas on the regulators and special values of motivic $L$-functions
which generalize the analytic class number formula.
For an elliptic curve $E$ over $\Q$, Beilinson proved 
that there is an integral symbol $\xi\in K_2(E)$
such that the real regulator $\reg_\R(\xi)\in H^2_\cD(E,\R(2))\cong \R$ 
agrees with the special value of the $L$-function $L(E,s)$ of $E$,
\[
\reg_\R(\xi)\sim_{\Q^\times}L'(E,0)\quad
(\Longleftrightarrow\,
\reg_\R(\xi)/L'(E,0)\in \Q^\times)
\]
where
$L'(E,s)=\frac{d}{ds}L(E,s)$ (\cite[Theorem 1.3]{B-reg}).

The $p$-adic counterpart of the Beilinson conjecture
was formulated by Perrin-Riou \cite[4.2.2]{Perrin-Riou},
which we call the {\it $p$-adic Beilinson conjecture}.
See also \cite{Colmez} for a general survey.
Her conjecture is formulated in terms of the $p$-adic etale regulators (which
are compatible with the syntomic regulators thanks to Besser's theorem) and the conjectural
$p$-adic measures which provide the $p$-adic $L$-functions of motives.
However
there are only a few of literatures
due to the extremal difficulty of the statement.
In a joint paper \cite{ABC} with Chida, 
we give a concise statement of the $p$-adic Beilinson
conjecture by restricting ourselves to the case of elliptic curves over $\Q$.

\begin{conj}[Weak $p$-adic Beilinson conjecture]\label{weakB-conj}
Let $E$ be an elliptuic curve over $\Q$.
Let $p>2$ be a prime at which $E$ has a good ordinary reduction.
Let $E_{\Q_p}:=E\times_\Q\Q_p$ and 
let $e_{\text{\rm unit}}\in H^1_\dR(X_\alpha/K)$ be the eigen vector for
the unit root $\alpha_{E,p}$ of the $p$-th Frobenius $\Phi$.
Let $L_p(E,\chi,s)$ be the $p$-adic $L$-function by Mazur and Swinnerton-Dyer
\cite{MS}.
Let $Q:H^1_\dR(E_{\Q_p}/\Q_p)^{\ot2}\to H^2_\dR(E_{\Q_p}/\Q_p)\cong\Q_p$ be the cup-product pairing.
Let 
\[
\reg_\syn:K_2(E)\to H^2_\syn(E,\Q_p(2))\cong H^1_\dR(E_{\Q_p}/\Q_p)
\]
be the syntomic regulator map. 
Fix a regular 1-form $\omega_E\in \vg(E,\Omega^1_{E/\Q})$.
Let $\omega:\Z_p^\times\to\mu_{p-1}$ be the Teichm\"uller character.
Then there is a constant $C\in \Q^\times$ which does not depend on $p$ such that
\[
(1-p\alpha_{E,p}^{-1})
\frac{Q(\reg_\syn(\xi),e_\unit)}{Q(\omega_E,e_\unit)}
=CL_p(E,\omega^{-1},0).
\]
\end{conj}
In \cite[Conjecture 3.3]{ABC}, we give a precise description of the constant $C$
in terms of the real regulator.

\subsection{Conjecture on Rogers-Zudilin type formulas}\label{RZ-sect}
In their paper \cite{RZ}, Rogers and Zudilin give descriptions of special values of
$L$-functions of elliptic curves in terms of
the hypergeometric functions ${}_3F_2$ or ${}_4F_3$. 
Apply theorems in \S \ref{elliptic-sect} to Conjecture \ref{weakB-conj}, we have
a statement of the $p$-adic counterpart of the Rogers-Zudilin type formulas.

\medskip

The following is obtained from Conjecture \ref{weakB-conj} and 
Corollary \ref{gauss-main2}
(note that the symbol \eqref{RZ-eq2} below
agrees with $\xi_\gauss|_{X_a}$ in \S \ref{gauss-sect} up to a constant). 
\begin{conj}\label{ell1-conj}
Let $X\to \P^1$ be the elliptic fibration defined by a Weierstrass equation
$y^2=x(1-x)(1-(1-t)x)$ (i.e. the hypergeometric curve of Gauss type, cf. \S \ref{gauss-sect}).
Let $X_a$ be the fiber at $t=a\in\Q\setminus\{0,1\}$.
Suppose that the symbol
\begin{equation}\label{RZ-eq2}
\xi_a=\left\{\frac{y-1+x}{y+1-x},\frac{a x^2}{(1-x)^2}\right\}\in K_2(X_a)
\end{equation}
is integral in the sense of Scholl \cite{Scholl}.
Let $p>2$ be a prime such that
$X_a$ has a good ordinary reduction at $p$. 
Let $\sigma_a:\Z_p[[t]]\to\Z_p[[t]]$ be the $p$-th Frobenius given by
$\sigma_a(t)=a^{1-p}t^p$.
Let $\alpha_{X_a,p}$ be the unit root.
Then there is a rational number $C_a\in \Q^\times$ not depending on $p$ such that
\[
(1-p\alpha^{-1}_{X_a,p})\cF_{\frac{1}{2},\frac{1}{2}}^{(\sigma_a)}(a)=C_a
L_p(X_a,\omega^{-1},0)
\]
where $\omega$ is the Teichm\"uller character.
\end{conj}
Here are examples of $a$ such that the symbol \eqref{RZ-eq2} is integral 
(cf. \cite[5.4]{A})
\[
a=-1,\pm2,\pm4,\pm8,\pm16,
\pm\frac{1}{2},  
\pm\frac{1}{8},\pm\frac{1}{4}, \pm\frac{1}{16}.
\]
F. Brunault compared
the symbol \eqref{RZ-eq2} with the Beilinson-Kato element in case $a=4$
(\cite[Appendix B]{ABC}).
Then, thanks to the main result of his paper \cite{regexp}, 
it follows that $\reg_\syn(\xi_a)$ gives the $p$-adic $L$-value
of $X_4$ (see \cite[Theorem 5.2]{ABC} for the precise statement). 
We thus have
\begin{thm}
\label{brunault}
When $a=4$, Conjecture \ref{ell1-conj} is true and $C_4=-1$.
\end{thm}

We obtain the following statements from
Theorems \ref{elliptic-thm1}, \ref{elliptic-thm2} and \ref{elliptic-thm3}.
\begin{conj}\label{ell2-conj}
Let $a\in \Q\setminus\{0,1\}$ and 
let $X_a$ be the ellptic curve over $\Q$ defined by an affine equation
$3y^2=2x^3-3x^2+1-a$. 
Suppose that the symbol
\begin{equation}\label{RZ-eq3}
\left\{
\frac{y-x+1}{y+x-1},
\frac{a}{2(x-1)^3}
\right\}\in K_2(X_a)
\end{equation}
is integral.
Let $p>5$ be a prime such that 
$X_a$ has a good ordinary reduction at $p$. 
Then there is a rational number $C_a\in \Q^\times$ not depending on $p$ such that
\[
(1-p\alpha^{-1}_{X_a,p})\cF_{\frac{1}{6},\frac{5}{6}}^{(\sigma_a)}(a)=C_a
L_p(X_a,\omega^{-1},0).
\]
\end{conj}
There are infinitely many $a$ such that
the symbol \eqref{RZ-eq3} is integral.
For example, if $a=1/n$ with
$n\in\Z_{\geq 2}$ and $n\equiv 0,2$ mod $6$,
then 
the symbol \eqref{RZ-eq3} is integral (cf. \cite[5.4]{A}).

\begin{conj}\label{ell3-conj}
Let $a\in \Q\setminus\{0,1\}$ and 
let $X_a$ be the ellptic curve over $\Q$ defined by an affine equation
$y^2=x^3+(3x+4a)^2$.
Suppose that the symbol
\begin{equation}\label{RZ-eq4}
\left\{
\frac{y-3x-4a}{-8a},
\frac{y+3x+4a}{8a}
\right\}
\in K_2(X_a)
\end{equation}
is integral.
Let $p>3$ be a prime such that 
$X_a$ has a good ordinary reduction at $p$. 
Then there is a rational number $C_a\in \Q^\times$ not depending on $p$ such that
\[
(1-p\alpha^{-1}_{X_a,p})\cF_{\frac{1}{3},\frac{2}{3}}^{(\sigma_a)}(a)=C_a
L_p(X_a,\omega^{-1},0).
\]
\end{conj}
If $a=\frac{1}{6n}$ with
$n\in\Z_{\geq 1}$ arbitrary,
then 
the symbol \eqref{RZ-eq4} is integral (cf. \cite[5.4]{A}).
\begin{conj}\label{ell4-conj}
Let $a\in \Q\setminus\{0,1\}$ and 
let $X_a$ be the elliptic curve over $\Q$ defined by an affine equation
$y^2=x^3-2x^2+(1-a)x$.
Suppose that the symbol
\begin{equation}\label{RZ-eq5}
\left\{\frac{y-(x-1)}{y+(x-1)},\frac{-a x}{(x-1)^3}\right\}
\in K_2(X_a)
\end{equation}
is integral.
Let $p>2$ be a prime such that
$X_a$ has a good ordinary reduction at $p$. 
Then there is a rational number $C_a\in \Q^\times$ not depending on $p$ such that
\[
(1-p\alpha^{-1}_{X_a,p})\cF_{\frac{1}{4},\frac{3}{4}}^{(\sigma_a)}(a)=C_a
L_p(X_a,\omega^{-1},0).
\]
\end{conj}
If the denominator of $j(X_a)=64(1+3a)^3/(a(1-a)^2)$ is prime to
$a$ (e.g. $a=1/n$, $n\in \Z_{\geq 2}$), then 
the symbol \eqref{RZ-eq5} is integral.

\medskip

From Corollary \ref{main-thm4} and Theorem \ref{fermat-main2}, 
we have the following conjectures.

\begin{conj}\label{ell6-conj}
Let $a\in \Q\setminus\{0,1\}$ and 
let $X_a$ be the ellptic curve over $\Q$ defined by an affine equation
$(x^2-1)(y^2-1)=a$.
Suppose that the symbol
\begin{equation}\label{RZ-eq6}
\left\{\frac{x-1}{x+1},\frac{y-1}{y+1}\right\}\in K_2(X_a)
\end{equation}
is integral.
Let $p>2$ be a prime such that 
$X_a$ has a good ordinary reduction at $p$. 
Then there is a rational number $C_a\in \Q^\times$ not depending on $p$ such that
\[
(1-p\alpha^{-1}_{X_a,p})\cF_{\frac{1}{2},\frac{1}{2}}^{(\sigma_a)}(1)=C_a
L_p(X_a,\omega^{-1},0).
\]
\end{conj}
If the denominator of $j(X_a)=16(a^2-16a+16)^3/((1-a)a^4)$
is prime to $a$ (e.g. $a=\pm 2^n$, $n\in \{\pm 1,\pm 2,\pm 3\}$), then 
the symbol \eqref{RZ-eq6} is integral.
\begin{conj}\label{ell7-conj}
Let $F_{N,M}$ be the Fermat curve defined by an affine equation 
$z^N+w^M=1$, and $F^*_{2,4}$ the curve $z^2=w^4+1$. 
Let $\sigma=\sigma_1$ (i.e. $\sigma(t)=t^p$).
Then there are rational numbers $C,C^\prime,
C^{\prime\prime}\in \Q^\times$ not depending on $p$ such that
\[
(1-p\alpha^{-1}_{F_{3,3},p})\cF_{\frac{1}{3},\frac{1}{3}}^{(\sigma)}(1)=C
L_p(F_{3,3},\omega^{-1},0),
\]
\[
(1-p\alpha^{-1}_{F_{2,4},p})\cF_{\frac{1}{2},\frac{1}{4}}^{(\sigma)}(1)=C^\prime
L_p(F_{2,4},\omega^{-1},0),
\]
\[
(1-p\alpha^{-1}_{F_{2,4}^*,p})\cF_{\frac{1}{4},\frac{1}{4}}^{(\sigma)}(1)=C^{\prime\prime}
L_p(F_{2,4}^*,\omega^{-1},0).
\]
\end{conj}

\noindent
Department of Mathematics, Hokkaido University,
\par\noindent
Sapporo 060-0810,
JAPAN

\smallskip

\noindent
{\it E-mail} : \textbf{asakura@math.sci.hokudai.ac.jp}
\end{document}